\documentclass[12pt,leqno]{article}
\usepackage[mathscr]{eucal}
\usepackage{amsmath,amssymb,latexsym,theorem,bbm,amsfonts}
\usepackage{verbatim}
\usepackage{color}
\usepackage{appendix}
\usepackage{array,epsfig,theorem,bbm,graphics,booktabs,url}

\newcommand{\CC}{\mathbb{C}}
\newcommand{\DD}{\mathbb{D}}

\newcommand{\NN}{\mathbb{N}}
\newcommand{\QQ}{\mathbb{Q}}
\newcommand{\RR}{\mathbb{R}}
\newcommand{\oRR}{\overline{\RR}}
\newcommand{\SSS}{\mathbb{S}}

\newcommand{\ZZ}{\mathbb{Z}}

\newcommand{\ba}{{\boldsymbol{a}}}

\newcommand{\ta}{\widetilde{a}}

\newcommand{\bc}{{\boldsymbol{c}}}

\newcommand{\tC}{\widetilde{C}}

\newcommand{\be}{{\boldsymbol{e}}}

\newcommand{\omu}{\overline{\mu}}

\newcommand{\bv}{{\boldsymbol{v}}}

\newcommand{\bx}{{\boldsymbol{x}}}
\newcommand{\bX}{{\boldsymbol{X}}}
\newcommand{\by}{{\boldsymbol{y}}}
\newcommand{\bY}{{\boldsymbol{Y}}}

\newcommand{\btheta}{{\boldsymbol{\theta}}}

\newcommand{\bTheta}{{\boldsymbol{\Theta}}}
\newcommand{\tnu}{\widetilde{\nu}}
\newcommand{\hmu}{\widehat{\mu}}

\newcommand{\tvare}{\widetilde{\vare}}
\newcommand{\bzero}{{\boldsymbol{0}}}
\newcommand{\bone}{{\boldsymbol{1}}}

\newcommand{\cB}{{\mathcal B}}

\newcommand{\cD}{{\mathcal D}}

\newcommand{\cM}{{\mathcal M}}

\newcommand{\cU}{{\mathcal U}}

\newcommand{\bcU}{\boldsymbol{\cU}}

\newcommand{\cX}{{\mathcal X}}

\newcommand{\bcX}{\boldsymbol{\cX}}
\newcommand{\tbcX}{\widetilde{\bcX}}

\newcommand{\cY}{{\mathcal Y}}
\newcommand{\cZ}{{\mathcal Z}}

\newcommand{\bcY}{\boldsymbol{\cY}}

\newcommand{\dd}{\mathrm{d}}
\newcommand{\ee}{\mathrm{e}}
\newcommand{\ff}{\mathrm{f}}
\newcommand{\ii}{\mathrm{i}}

\newcommand{\EE}{\operatorname{\mathbb{E}}}

\newcommand{\PP}{\operatorname{\mathbb{P}}}
\newcommand{\oo}{\operatorname{o}}

\newcommand{\sign}{\operatorname{sign}}

\newcommand{\hC}{\widehat{C}}

\newcommand{\ovarrho}{\overline{\varrho}}

\newcommand{\vare}{\varepsilon}

\renewcommand{\mid}{\,|\,}

\renewcommand{\leq}{\leqslant}
\renewcommand{\geq}{\geqslant}

\newcommand{\distr}{\stackrel{\cD}{\longrightarrow}}
\newcommand{\distrf}{\stackrel{\cD_\ff}{\longrightarrow}}
\newcommand{\distre}{\stackrel{\cD}{=}}
\newcommand{\distrw}{\stackrel{w}{\longrightarrow}}
\newcommand{\distrv}{\stackrel{v}{\longrightarrow}}

\newcommand{\as}{\stackrel{{\mathrm{a.s.}}}{\longrightarrow}}

\newcommand{\bbone}{\mathbbm{1}}

\newcommand{\nt}{{\lfloor nt\rfloor}}

\newcommand{\proofend}{\hfill\mbox{$\Box$}}

\numberwithin{equation}{section}

\theoremstyle{change} \theorembodyfont{\em}
\newtheorem{Lem}{Lemma.}[section]
\newtheorem{Thm}[Lem]{Theorem.}
\newtheorem{Pro}[Lem]{Proposition.}
\newtheorem{Cor}[Lem]{Corollary.}
\newtheorem{Def}[Lem]{Definition.}

\theorembodyfont{\rm}
\newtheorem{Rem}[Lem]{Remark.}

\begin{document}

\begin{center}
 {\bfseries\Large On aggregation of subcritical Galton--Watson branching \\[2mm]
                   processes with regularly varying immigration}

\vskip0.5cm

 {\sc\large
  M\'aty\'as $\text{Barczy}^{*,\diamond}$,
  Fanni K. $\text{Ned\'enyi}^{*}$,
  Gyula $\text{Pap}^{**}$}

\end{center}

\vskip0.2cm

\noindent
 * MTA-SZTE Analysis and Stochastics Research Group,
   Bolyai Institute, University of Szeged,
   Aradi v\'ertan\'uk tere 1, H--6720 Szeged, Hungary.

\noindent
 ** Bolyai Institute, University of Szeged,
     Aradi v\'ertan\'uk tere 1, H--6720 Szeged, Hungary.

\noindent e--mails: barczy@math.u-szeged.hu (M. Barczy),
                    nfanni@math.u-szeged.hu (F. K. Ned\'enyi).

\noindent $\diamond$ Corresponding author.

\vskip0.2cm


\renewcommand{\thefootnote}{}
\vspace*{0.2cm}
\footnote{M\'aty\'as Barczy is supported by the J\'anos Bolyai Research Scholarship of the Hungarian Academy of Sciences.
Fanni K. Ned\'enyi is supported by the UNKP-18-3 New National
Excellence Program of the Ministry of Human Capacities.
Gyula Pap was supported by the Ministry for Innovation and Technology, Hungary grant TUDFO/47138-1/2019-ITM.}

\vspace*{-10mm}

\begin{abstract}
We study an iterated temporal and contemporaneous aggregation of \ $N$ \ independent copies of a strongly stationary subcritical Galton--Watson branching process with regularly varying immigration having index \ $\alpha \in (0, 2)$.
\ Limits of finite dimensional distributions of appropriately centered and scaled aggregated partial sum processes are shown to exist
 when first taking the limit as \ $N \to \infty$ \ and then the time scale \ $n \to \infty$.
\ The limit process is an \ $\alpha$-stable process if \ $\alpha \in (0, 1) \cup (1, 2)$, \ and a deterministic line with slope \ $1$ \ if \ $\alpha = 1$.
\end{abstract}

\section{Introduction}
\label{intro}

The field of temporal and contemporaneous (also called cross-sectional) aggregations of independent stationary stochastic processes is an important and very active
 research area in the empirical and theoretical statistics and in other areas as well.
Robinson \cite{Rob} and Granger \cite{Gra} started to investigate the scheme of contemporaneous aggregation of
 random-co\-ef\-fi\-cient autoregressive processes of order 1 in order to obtain the long memory phenomenon in aggregated time series.
For surveys on aggregation of different kinds of stochastic processes, see, e.g., Pilipauskait{\.e} and Surgailis \cite{PilSur},
 Jirak \cite[page 512]{Jir} or the arXiv version of Barczy et al.\ \cite{BarNedPap}.

Recently, Puplinskait{\.e} and Surgailis \cite{PupSur1, PupSur2} studied iterated aggregation of random coefficient autoregressive processes of order 1
 with common innovations and with so-called idiosyncratic innovations, respectively, belonging to the domain of attraction of an $\alpha$-stable law.
Limits of finite dimensional distributions of appropriately centered and scaled aggregated partial
 sum processes are shown to exist when first the number of copies \ $N \to \infty$ \ and
 then the time scale \ $n \to \infty$.
\ Very recently, Pilipauskait{\.e} et al.\ \cite{PilSkoSur} extended the results of Puplinskait{\.e} and Surgailis \cite{PupSur2} (idiosyncratic case)
 deriving limits of finite dimensional distributions of appropriately centered and scaled aggregated partial sum processes
 when first the time scale \ $n \to \infty$ \ and then the number of copies \ $N \to \infty$, \ and when \ $n \to \infty$ \ and
 \ $N \to \infty$ \ simultaneously with possibly different rates.

The above listed references are all about aggregation procedures for times series, mainly for randomized autoregressive processes. According to our knowledge this question has not been studied before in the literature.
The present paper investigates aggregation schemes for some branching processes with low moment condition.
Branching processes, especially Galton--Watson branching processes with immigration, have attracted a lot of attention due to the fact that they
 are widely used in mathematical biology for modelling the growth of a population in time.
In Barczy et al.\ \cite{BarNedPap2}, we started to investigate the limit behavior of temporal and contemporaneous aggregations of independent copies
 of a stationary multitype Galton--Watson branching process with immigration under third order moment conditions on the offspring and immigration
 distributions in the iterated and simultaneous cases as well.
In both cases, the limit process is a zero mean Brownian motion with the same covariance function.
As of 2020, modeling the COVID-19 contamination of the population of a certain region or country is of great importance.
Multitype Galton--Watson processes with immigration have been frequently used to model the spreading of a number of diseases, and they can be applied for this new disease as well.
For example, Yanev et al. \cite{YanStoAta} applied a two-type Galton--Watson process with immigration to model the number of detected, COVID-19-infected and undetected,
 COVID-19-infected people in a population. The temporal and contemporaneous aggregation of the first coordinate process of
 the two-type branching process in question would mean the total number of detected, infected people up to some given time point across several regions.

In this paper we study the limit behavior of temporal and contemporaneous aggregations of independent copies of
 a strongly stationary Galton--Watson branching process \ $(X_k)_{k\geq0}$ \ with regularly varying immigration having index in \ $(0, 2)$ \ (yielding infinite variance)
 in an iterated, idiosyncratic case, namely, when first the number of copies \ $N \to \infty$ \ and then the time scale \ $n \to \infty$.
\ Our results are analogous to those of Puplinskait{\.e} and Surgailis \cite{PupSur2}.

The present paper is organized as follows.
In Section \ref{results}, first we collect our assumptions
 that are valid for the whole paper, namely, we consider a sequence of independent copies of \ $(X_k)_{k\geq0}$ \ such that the expectation of
 the offspring distribution is less than \ $1$ \ (so-called subcritical case).
In case of  \ $\alpha\in[1,2)$, \ we additionally suppose the finiteness of the second moment of the offspring distribution.
Under our assumptions, by Basrak et al. \cite[Theorem 2.1.1]{BasKulPal} (see also Theorem \ref{Xtail}), the unique stationary distribution of \ $(X_k)_{k\geq0}$ \ is also regularly varying with the same index \ $\alpha$.

In Theorem \ref{simple_aggregation1_stable_fd}, we show that the appropriately centered and scaled partial sum process of finite segments of independent copies
 of \ $(X_k)_{k\geq0}$ \ converges to an \ $\alpha$-stable process.
The characteristic function of the \ $\alpha$-stable limit process is given explicitly as well.
In Remarks \ref{mu_properties1} and \ref{mu_properties2}, we collect some properties of the \ $\alpha$-stable limit process in question,
 such as the support of its L\'evy measure.
The proof of Theorem \ref{simple_aggregation1_stable_fd} is based on a slight modification of Theorem 7.1 in Resnick \cite{Res}, namely,
 on a result of weak convergence of partial sum processes towards L\'evy processes, see Theorem \ref{7.1}, where we consider a  different centering.
In the course of the proof of Theorem \ref{simple_aggregation1_stable_fd} one needs to verify that the so-called limit measures of finite segments of
 \ $(X_k)_{k\geq 0}$ \ are in fact L\'evy measures.
We determine these limit measures explicitly (see part (i) of Proposition \ref{Pro_limit_meaure}) applying an expression for the so-called tail measure of
 a strongly stationary regularly varying sequence based on the corresponding (whole) spectral tail process given
 in Planini\'c and Soulier \cite[Theorem 3.1]{PlaSou}.

While the centering in Theorem \ref{simple_aggregation1_stable_fd} is the so-called truncated mean,
 in Corollary \ref{simple_aggregation1_stable_centering_fd} we consider no-centering if \ $\alpha\in(0,1)$, \
 and centering with the mean if \ $\alpha\in(1,2)$.
\ In both cases the limit process is an \ $\alpha$-stable process, the same one as in Theorem \ref{simple_aggregation1_stable_fd} plus some deterministic drift
 depending on \ $\alpha$.
\ Theorem \ref{simple_aggregation1_stable_fd} and Corollary \ref{simple_aggregation1_stable_centering_fd} together yield the weak convergence of finite dimensional
 distributions of appropriately centered and scaled contemporaneous aggregations of independent copies of \ $(X_k)_{k\geq0}$ \
 towards the corresponding finite dimensional distributions of a strongly stationary, subcritical autoregressive process of order 1 with
  \ $\alpha$-stable innovations as the number of copies tends to infinity, see Corollary \ref{aggr_copies1} and Proposition \ref{Pro_AR1}.

Theorem \ref{iterated_aggr_1} contains our main result, namely, we determine the weak limit of appropriately centered and scaled
 finite dimensional distributions of temporal and contemporaneous aggregations of independent copies of \ $(X_k)_{k\geq0}$, \ where the limit is taken in a way that first the number of copies tends to infinity and then the time corresponding to temporal
 aggregation tends to infinity.
It turns out that the limit process is an \ $\alpha$-stable process if \ $\alpha \in (0, 1) \cup (1, 2)$, \ and a deterministic line with slope
 \ $1$ \ if \ $\alpha = 1$.
\ We consider different kinds of centerings, and we give the explicit characteristic function of the limit process as well.
In Remark \ref{Rem_char_spectral_proc}, we rewrite this characteristic function in case of \ $\alpha\in(0,1)$ \
 in terms of the spectral tail process of \ $(X_k)_{k\geq 0}$.

We close the paper with five appendices.
In Appendix \ref{App_cont_map_theorem} we recall a version of the continuous mapping theorem due to Kallenberg \cite[Theorem 3.27]{Kal}.
Appendix \ref{App_vague} is devoted to some properties of the underlying punctured space \ $\RR^d\setminus\{\bzero\}$ \ and
 vague convergence.
In Appendix \ref{App_reg_var_distr} we recall the notion of a regularly varying random vector and its limit measure, and, in Proposition \ref{Pro_mapping}, the limit measure of an appropriate positively homogeneous real-valued function
 of a regularly varying random vector.
In Appendix \ref{App_Resnick_gen} we formulate a result on weak convergence of partial sum processes towards L\'evy processes
 by slightly modifying Theorem 7.1 in Resnick \cite{Res} with a different centering.
In the end, we recall a result on the tail behavior and forward tail process of \ $(X_k)_{k\geq0}$ \ due to Basrak et al.\ \cite{BasKulPal}, and we determine the limit measures of finite segments of \ $(X_k)_{k\geq0}$, \ see Appendix \ref{App_tail}.

Finally, we summarize the novelties of the paper.
According to our knowledge, studying aggregation of regularly varying Galton--Watson branching processes with immigration
 has not been considered before.
In the proofs we make use of the explicit form of the (whole) spectral tail process
 and a very recent result of Planini\'c and Soulier \cite[Theorem 3.1]{PlaSou} about the tail measure of strongly stationary sequences.
We explicitly determine the limit measures of finite segments of \ $(X_k)_{k\geq0}$, \ see part (i) of Proposition \ref{Pro_limit_meaure}.

In a companion paper, we will study the other iterated, idiosyncratic aggregation scheme, namely, when first the time scale
 \ $n \to \infty$ \ and then the number of copies \ $N \to \infty$.

\section{Main results}
\label{results}

Let \ $\ZZ_+$, \ $\NN$, \ $\QQ$, \ $\RR$, \ $\RR_+$, \ $\RR_{++}$, \ $\RR_-$,
 \ $\RR_{--}$ \ and \ $\CC$ \ denote the set of non-negative integers, positive
 integers, rational numbers, real numbers, non-negative real numbers, positive real numbers,
 non-positive real numbers, negative real numbers and complex numbers, respectively.
For each \ $d \in \NN$, \ the natural basis in \ $\RR^d$ \ will be denoted by \ $\be_1$, \ldots, $\be_d$.
\ Put \ $\bone_d := (1, \ldots, 1)^\top$ \ and \ $\SSS^{d-1} := \{\bx \in \RR^d : \|\bx\| = 1\}$,
 \ where \ $\|\bx\|$ \ denotes the Euclidean norm of \ $\bx\in\RR^d$,
 \ and denote by \ $\cB(\SSS^{d-1})$ \ the Borel \ $\sigma$-field of \ $\SSS^{d-1}$.
\ For a probability measure \ $\mu$ \ on \ $\RR^d$, \ $\hmu$ \ will denote its characteristic function, i.e., \ $\hmu(\btheta) := \int_{\RR^d} \ee^{\ii\langle\btheta,\bx\rangle} \, \mu(\dd\bx)$ \ for \ $\btheta \in \RR^d$.
\ Convergence in distributions and almost sure convergence of random variables, and weak convergence of probability measures will be denoted by \ $\distr$, \ $\as$ \ and \ $\distrw$, \ respectively.
Equality in distribution will be denoted by \ $\distre$.
\ We will use \ $\distrf$ \ or \ $\cD_\ff\text{-}\hspace*{-1mm}\lim$ \ for weak
 convergence of finite dimensional distributions.
A function \ $f : \RR_+ \to \RR^d$ \ is called \emph{c\`adl\`ag} if it is right
 continuous with left limits.
Let \ $\DD(\RR_+, \RR^d)$ \ and \ $\CC(\RR_+, \RR^d)$ \ denote the space of
 all \ $\RR^d$-valued c\`adl\`ag and continuous functions on \ $\RR_+$,
 \ respectively.
Let \ $\cB(\DD(\RR_+, \RR^d))$ \ denote the Borel \ $\sigma$-algebra on
 \ $\DD(\RR_+, \RR^d)$ \ for the metric defined in Chapter VI, (1.26) of Jacod and Shiryaev \cite{JacShi}.
With this metric \ $\DD(\RR_+, \RR^d)$ \ is a
 complete and separable metric space and the topology induced by this metric is
 the so-called Skorokhod topology.
For \ $\RR^d$-valued stochastic processes \ $(\bcY_t)_{t \in \RR_+}$ \ and
 \ $(\bcY^{(n)}_t)_{t \in \RR_+}$, \ $n \in \NN$, \ with c\`adl\`ag paths we write
 \ $\bcY^{(n)} \distr \bcY$ \ as \ $n\to\infty$ \ if the distribution of \ $\bcY^{(n)}$ \ on the
 space \ $(\DD(\RR_+, \RR^d), \cB(\DD(\RR_+, \RR^d)))$ \ converges weakly to the
 distribution of \ $\bcY$ \ on the space
 \ $(\DD(\RR_+, \RR^d), \cB(\DD(\RR_+, \RR^d)))$ \ as \ $n \to \infty$.

Let \ $(X_k)_{k\in\ZZ_+}$ \ be a Galton--Watson branching process with immigration.
For each \ $k, j \in \ZZ_+$, \ the number of individuals in the \ $k^\mathrm{th}$
 \ generation will be denoted by \ $X_k$, \ the number of offsprings produced by
 the \ $j^\mathrm{th}$ \ individual belonging to the \ $(k-1)^\mathrm{th}$
 \ generation will be denoted by \ $\xi_{k,j}$, \ and the number of immigrants in the
 \ $k^\mathrm{th}$ \ generation will be denoted by \ $\vare_k$.
\ Then we have
 \[
   X_k = \sum_{j=1}^{X_{k-1}} \xi_{k,j} + \vare_k , \qquad k \in \NN ,
 \]
 where we define \ $\sum_{j=1}^0 := 0$.
\ Here \ $\bigl\{X_0, \, \xi_{k,j}, \, \vare_k : k, j \in \NN\bigr\}$ \ are supposed
 to be independent non-negative integer-valued random variables.
Moreover, \ $\{\xi_{k,j} : k, j \in \NN\}$ \ and
 \ $\{\vare_k : k \in \NN\}$ \ are supposed to consist of identically
 distributed random variables, respectively.
For notational convenience, let \ $\xi$ \ and \ $\vare$ \ be independent random variables such that \ $\xi \distre \xi_{1,1}$ \ and \ $\vare \distre \vare_1$.

If \ $m_\xi := \EE(\xi) \in [0, 1)$ \ and \ $\sum_{\ell=1}^\infty \log(\ell) \PP(\vare = \ell) < \infty$, \ then the Markov chain \ $(X_k)_{k\in\ZZ_+}$ \ admits a unique stationary distribution \ $\pi$, \ see, e.g., Quine \cite{Qui}.
Note that if \ $m_\xi \in [0, 1)$ \ and \ $\PP(\vare = 0) = 1$, \ then \ $\sum_{\ell=1}^\infty \log(\ell) \PP(\vare = \ell) = 0$ \ and \ $\pi$ \ is the Dirac measure \ $\delta_0$ \ concentrated at the point \ $0$.
\ In fact, \ $\pi = \delta_0$ \ if and only if \ $\PP(\vare = 0) = 1$.
\ Moreover, if \ $m_\xi = 0$ \ (which is equivalent to \ $\PP(\xi = 0) = 1$), \ then \ $\pi$ \ is the distribution of \ $\vare$.

In what follows, we formulate our assumptions valid for the whole paper.
We assume that \ $m_\xi \in [0, 1)$ \ (so-called subcritical case) and \ $\vare$ \ is regularly varying with index \ $\alpha \in (0, 2)$, \ i.e.,
 \ $\PP(\vare > x)\in\RR_{++}$ \ for all \ $x\in\RR_{++}$ \ and
 \[
   \lim_{x\to\infty} \frac{\PP(\vare > qx)}{\PP(\vare > x)} = q^{-\alpha} \qquad
   \text{for all \ $q \in \RR_{++}$.}
 \]
Then \ $\PP(\vare = 0) < 1$ \ and \ $\sum_{\ell=1}^\infty \log(\ell) \PP(\vare = \ell) < \infty$, \ see, e.g., Barczy et al.\ \cite[Lemma E.5]{BarBosPap},
 hence the Markov process \ $(X_k)_{k\in\ZZ_+}$ \ admits a unique stationary distribution \ $\pi$.
\ We suppose that \ $X_0 \distre \pi$, \ yielding that the Markov chain \ $(X_k)_{k\in\ZZ_+}$ \ is strongly stationary.
In case of \ $\alpha \in [1, 2)$, \ we suppose additionally that \ $\EE(\xi^2) < \infty$.
\ By Basrak et al. \cite[Theorem 2.1.1]{BasKulPal} (see also Theorem \ref{Xtail}), \ $X_0$ \ is regularly varying with index \ $\alpha$, \ yielding the existence of a sequence \ $(a_N)_{N\in\NN}$ \ in \ $\RR_{++}$ \ with
 \ $N \PP(X_0 > a_N) \to 1$ \ as \ $N \to \infty$, \ see, e.g., Lemma \ref{a_n}.
Let us fix an arbitrary sequence \ $(a_N)_{N\in\NN}$ \ in \ $\RR_{++}$ \ with this property.
In fact, \ $a_N = N^{\frac{1}{\alpha} } L(N)$, \ $N \in \NN$, \ for some slowly varying continuous function \ $L : \RR_{++} \to \RR_{++}$, \ see, e.g.,
 Araujo and Gin\'e \cite[Exercise 6 on page 90]{AraGin}.
Let \ $X^{(j)} = (X^{(j)}_k)_{k\in\ZZ_+}$, \ $j \in \NN$, \ be a sequence of
 independent copies of \ $(X_k)_{k\in\ZZ_+}$.
\ We mention that we consider so-called idiosyncratic immigrations, i.e., the immigrations \ $(\vare^{(j)}_k)_{k\in\NN}$, \ $j \in \NN$, \ belonging to \ $(X^{(j)}_k)_{k\in\ZZ_+}$, \ $j \in \NN$, \ are independent.
One could study the case of common immigrations as well, i.e., when \ $(\vare^{(j)}_k)_{k\in\NN} = (\vare^{(1)}_k)_{k\in\NN}$, \ $j \in \NN$.

\begin{Thm}\label{simple_aggregation1_stable_fd}
For each \ $k \in \ZZ_+$,
 \begin{equation}\label{help6_simple_aggregation1_stable_fd}
  \begin{aligned}
   &\biggl(\frac{1}{a_N}
  \sum_{j=1}^{\lfloor Nt\rfloor}
   \Bigl(X^{(j)}_0
         - \EE\bigl(X^{(j)}_0
                    \bbone_{\{X^{(j)}_0\leq a_N\}}\bigr), \ldots, X^{(j)}_k
         - \EE\bigl(X^{(j)}_k
                    \bbone_{\{X^{(j)}_k\leq a_N\}}\bigr)\Bigr)^\top\biggr)_{t\in\RR_+} \\
   &= \biggl(\frac{1}{a_N}
  \sum_{j=1}^{\lfloor Nt\rfloor}
   \bigl(X^{(j)}_0, \ldots, X^{(j)}_k\bigr)^\top
         - \frac{\lfloor Nt\rfloor}{a_N} \EE\bigl(X_0
                    \bbone_{\{X_0\leq a_N\}}\bigr) \bone_{k+1}\biggr)_{t\in\RR_+}
    \distr \bigl(\bcX_t^{(k,\alpha)}\bigr)_{t\in\RR_+}
  \end{aligned}
 \end{equation}
 as \ $N\to\infty$, \ where \ $\bigl(\bcX_t^{(k,\alpha)}\bigr)_{t\in\RR_+}$ \ is a \ $(k + 1)$-dimensional \ $\alpha$-stable process such that the characteristic function of the distribution \ $\mu_{k,\alpha}$ \ of \ $\bcX_1^{(k,\alpha)}$ \ has the form
 \begin{align*}
  &\widehat{\mu_{k,\alpha}}(\btheta) \\
  &= \exp\biggl\{(1 - m_\xi^\alpha)
                  \sum_{j=0}^k
                   \int_0^\infty
                    \biggl(\ee^{\ii\langle\btheta,\bv_j^{(k)}\rangle u} - 1
                    - \ii u \sum_{\ell=j+1}^{k+1} \langle\be_\ell, \btheta\rangle \langle\be_\ell, \bv_j^{(k)}\rangle \bbone_{(0,1]}(u\langle\be_\ell, \bv_j^{(k)}\rangle)\biggr)
                    \alpha u^{-1-\alpha} \, \dd u\biggr\}
 \end{align*}
 for \ $\btheta \in \RR^{k+1}$ \ with the \ $(k + 1)$-dimensional vectors
 \[
  \bv_0^{(k)} := (1 - m_\xi^\alpha)^{-\frac{1}{\alpha}}
           \begin{bmatrix}
            1 \\ m_\xi \\ m_\xi^2 \\ \vdots \\ m_\xi^k
           \end{bmatrix} , \quad
  \bv_1^{(k)} :=
           \begin{bmatrix}
            0 \\ 1 \\ m_\xi \\ \vdots \\ m_\xi^{k-1}
           \end{bmatrix} , \quad
  \bv_2^{(k)} :=
           \begin{bmatrix}
            0 \\ 0 \\ 1 \\ \vdots \\ m_\xi^{k-2}
           \end{bmatrix} ,
           \quad \ldots , \quad
  \bv_k^{(k)} := \begin{bmatrix}
            0 \\ 0 \\ \vdots \\ 0 \\ 1
           \end{bmatrix} .
 \]
Moreover, for \ $\btheta \in \RR^{k+1}$,
 \[
   \widehat{\mu_{k,\alpha}}(\btheta)
   = \begin{cases}
      \exp\Bigl\{- C_\alpha (1 - m_\xi^\alpha)
                  \sum_{j=0}^k
                   |\langle\btheta,\bv_j^{(k)}\rangle|^\alpha
                    \left(1 - \ii \tan\left(\frac{\pi\alpha}{2}\right)
                                  \sign(\langle\btheta,\bv_j^{(k)}\rangle)\right) \\[2mm]
      \phantom{\exp\Bigl\{}
                 - \ii \frac{\alpha}{1-\alpha} \langle\btheta,
\bone_{k+1}\rangle\Bigr\} ,
       & \text{if \ $\alpha \ne 1$,} \\[4mm]
      \exp\Bigl\{- C_1 (1 - m_\xi)
                  \sum_{j=0}^k
                   |\langle\btheta,\bv_j^{(k)}\rangle|
                   \bigl(1 + \ii \frac{2}{\pi} \sign(\langle\btheta,\bv_j^{(k)}\rangle)
                                 \log(|\langle\btheta,\bv_j^{(k)}\rangle|)\bigr) \\[2mm]
      \phantom{\exp\Bigl\{}
                 + \ii C \langle\btheta,
\bone_{k+1}\rangle \\[2mm]
      \phantom{\exp\Bigl\{}
                 + \ii (1 - m_\xi) \sum_{j=0}^k \sum_{\ell=j+1}^{k+1} \langle\be_\ell, \btheta\rangle
                   \langle\be_\ell, \bv_j^{(k)}\rangle \log(\langle\be_\ell, \bv_j^{(k)}\rangle)\Bigr\} ,
       & \text{if \ $\alpha = 1$,}
     \end{cases}
 \]
 with the convention \ $0 \log(0) := 0$,
 \[
   C_\alpha
   := \begin{cases}
       \frac{\Gamma(2 - \alpha)}{1-\alpha}
       \cos\left(\frac{\pi\alpha}{2}\right) ,
        & \text{if \ $\alpha \ne 1$,} \\
       \frac{\pi}{2} , & \text{if \ $\alpha = 1$,}
      \end{cases}
  \]
 and
 \begin{equation}\label{c}
  C := \int_1^\infty u^{-2} \sin(u) \, \dd u
       + \int_0^1 u^{-2} (\sin(u) - u) \, \dd u .
 \end{equation}
\end{Thm}

Note that \ $C$ \ exists and is finite, since \ $\int_1^\infty u^{-2} |\sin(u)| \, \dd u \leq \int_1^\infty u^{-2} \, \dd u = 1$, \ and, by L'H\^ospital's rule, \ $\lim_{u\to0} u^{-2} (\sin(u) - u) = 0$, \ hence the integrand \ $u^{-2} (\sin(u) - u)$ \ can be extended to \ $[0, 1]$ \ continuously, yielding that its integral on \ $[0, 1]$ \ is finite.

Note that the scaling and the centering in \eqref{help6_simple_aggregation1_stable_fd} do not depend on \ $j$ \ or \ $k$, \ since the copies are independent and the process \ $(X_k)_{k\in\ZZ_+}$ \ is strongly stationary, and especially, \ $\EE\bigl(X^{(j)}_k \bbone_{\{X^{(j)}_k\leq a_N\}}\bigr) = \EE(X_0 \bbone_{\{X_0\leq a_N\}})$ \ for all \ $j \in \NN$ \ and \ $k \in \ZZ_+$.

The next two remarks are devoted to the study of some properties of \ $\mu_{k,\alpha}$.

\begin{Rem}\label{mu_properties1}
By the proof of Theorem \ref{simple_aggregation1_stable_fd} (see \eqref{fint}), it turns out that the L\'evy measure of \ $\mu_{k,\alpha}$ \ is
 \[
   \nu_{k,\alpha}(B)
   = (1 - m_\xi^\alpha) \sum_{j=0}^k \|\bv_j^{(k)}\|^\alpha \int_0^\infty \bbone_B\left(u \frac{\bv_j^{(k)}}{\|\bv_j^{(k)}\|}\right) \alpha u^{-\alpha-1} \, \dd u , \qquad B \in \cB(\RR^{k+1}_0) ,
 \]
 where the space \ $\RR^{k+1}_0 := \RR^{k+1} \setminus \{\bzero\}$ \ and its topological properties are discussed in Appendix \ref{App_vague}.
The radial part of \ $ \nu_{k,\alpha}$ \ is \ $u^{-\alpha-1} \, \dd u$, \ and the spherical part of \ $\nu_{k,\alpha}$ \ is any positive constant
 multiple of the measure \ $\sum_{j=0}^k\|\bv_j^{(k)}\|^\alpha \epsilon_{\bv_j^{(k)}/\|\bv_j^{(k)}\|}$  \ on \ $\SSS^k$, \  where
 for any \ $\bx \in \RR^{k+1}$, \ $\epsilon_{\bx}$ \ denotes the Dirac measure concentrated at the point \ $\bx$.
\ Particularly, the support of \ $\nu_{k,\alpha}$ \ is \ $\cup_{j=0}^k (\RR_{++} \bv_j^{(k)})$.
\ The vectors \ $\bv_0^{(k)}$, \ldots, $\bv_k^{(k)}$ \ form a basis in \ $\RR^{k+1}$, \ hence there is no proper linear subspace \ $V$ \ of \ $\RR^{k+1}$ \ covering
the support of \ $\nu_{k,\alpha}$.
\ Consequently, \ $\mu_{k,\alpha}$ \ is a nondegenerate measure in the sense that there are no \ $\ba \in \RR^{k+1}$ \ and a proper linear subspace \ $V$ \ of \ $\RR^{k+1}$ \ such that \ $a + V$ \ covers the support of \ $\mu_{k,\alpha}$, \ see, e.g., Sato \cite[Proposition 24.17 (ii)]{Sato}.
\proofend
\end{Rem}

\begin{Rem}\label{mu_properties2}
If \ $\alpha \in (0, 1)$, \ then, for each \ $\btheta \in \RR^{k+1}$,
 \[
   \widehat{\mu_{k,\alpha}}(\btheta) \\
   = \exp\biggl\{(1 - m_\xi^\alpha)
                  \sum_{j=0}^k
                   \int_0^\infty
                    \bigl(\ee^{\ii\langle\btheta,\bv_j^{(k)}\rangle u} - 1\bigr)
                    \alpha u^{-1-\alpha} \, \dd u
                 - \ii \frac{\alpha}{1-\alpha}
                   \langle\btheta, \bone_{k+1}\rangle\biggr\} ,
 \]
 see the proof of Theorem \ref{simple_aggregation1_stable_fd}.
Consequently, the drift of \ $\mu_{k,\alpha}$ \ is \ $- \frac{\alpha}{1-\alpha} \bone_{k+1}$, \ see, e.g., Sato \cite[Remark 14.6]{Sato}.
\ This drift is nonzero, hence \ $\mu_{k,\alpha}$ \ is not strictly
 $\alpha$-stable, see, e.g., Sato \cite[Theorem 14.7 (iv) and Definition 13.2]{Sato}.

The \ $1$-stable probability measure \ $\mu_{k,1}$ \ is not strictly \ $1$-stable, since the spherical part of its nonzero L\'evy measure \ $\nu_{k,1}$ \ is concentrated on \ $\RR_+^{k+1} \cap \SSS^k$, \ and hence the condition (14.12) in Sato \cite[Theorem 14.7 (v)]{Sato} is not satisfied.

If \ $\alpha \in (1, 2)$, \ then, for each \ $\btheta \in \RR^{k+1}$,
 \[
   \widehat{\mu_{k,\alpha}}(\btheta) \\
   = \exp\biggl\{(1 - m_\xi^\alpha)
                  \sum_{j=0}^k
                   \int_0^\infty
                    \bigl(\ee^{\ii\langle\btheta,\bv_j^{(k)}\rangle u} - 1 - \ii \langle\btheta,\bv_j^{(k)}\rangle u\bigr)
                    \alpha u^{-1-\alpha} \, \dd u
                 + \ii \frac{\alpha}{\alpha-1}
                   \langle\btheta, \bone_{k+1}\rangle\biggr\} ,
 \]
 see the proof of Theorem \ref{simple_aggregation1_stable_fd}.
Consequently, the center of \ $\mu_{k,\alpha}$ \ is \ $\frac{\alpha}{\alpha-1} \bone_{k+1}$, \ which is, in fact, the expectation of \ $\mu_{k,\alpha}$, \ and it
 is nonzero, and hence \ $\mu_{k,\alpha}$ \ is not strictly stable, see, e.g., Sato
 \cite[Theorem 14.7 (vi) and Definition 13.2]{Sato}.

All in all, \ $\mu_{k,\alpha}$ \ is not strictly \ $\alpha$-stable, but \ $\alpha$-stable for any
 \ $\alpha \in (0, 2)$.
\ We also note that \ $\mu_{k,\alpha}$ \ is absolutely continuous, see, e.g., Sato \cite[Theorem 27.4 and Proposition 14.5]{Sato}.
\proofend
\end{Rem}

The centering in Theorem \ref{simple_aggregation1_stable_fd} can be simplified in case of \ $\alpha \ne 1$.
\ Namely, if \ $\alpha \in (0, 1]$, \ then for each \ $t \in \RR_{++}$, \ by Lemma \ref{truncated_moments},
 \begin{equation}\label{centering01T}
 \begin{aligned}
  \frac{\lfloor Nt\rfloor}{a_N} \EE(X_0 \bbone_{\{X_0\leq a_N\}})
  &= \frac{\lfloor Nt\rfloor}{N}
     \frac{\EE(X_0 \bbone_{\{X_0\leq a_N\}})}{a_N\PP(X_0 > a_N)}
     N \PP(X_0> a_N) \\
  &\to \begin{cases}
        \frac{\alpha}{1-\alpha} t & \text{for \ $\alpha \in (0, 1)$,} \\
        \infty & \text{for \ $\alpha = 1$}
       \end{cases}
   \qquad \text{as \ $N \to \infty$.}
 \end{aligned}
 \end{equation}
In a similar way, if \ $\alpha \in (1, 2)$, \ then for each \ $t \in \RR_{++}$,
 \[
   \frac{\lfloor Nt\rfloor}{a_N} \EE(X_0 \bbone_{\{X_0\leq a_N\}})
   = \frac{\lfloor Nt\rfloor}{a_N} \EE(X_0)
     - \frac{\lfloor Nt\rfloor}{a_N} \EE(X_0 \bbone_{\{X_0>a_N\}}) ,
 \]
 where, \ $\lim_{N\to\infty} \frac{\lfloor Nt\rfloor}{a_N}  = \lim_{N\to\infty} t N^{1-\frac{1}{\alpha}} L(N)^{-1} = \infty$, \ and,
 by Lemma \ref{truncated_moments},
 \begin{equation}\label{centering12T}
   \frac{\lfloor Nt\rfloor}{a_N} \EE(X_0 \bbone_{\{X_0>a_N\}})
   \to \frac{\alpha}{\alpha-1} t
      \qquad \text{as \ $N \to \infty$.}
 \end{equation}
This shows that in case of \ $\alpha \in (0, 1)$, \ there is no need for centering, in case of \ $\alpha \in (1, 2)$ \ one can center with the expectation as well,
 while in case of \ $\alpha = 1$, \ neither non-centering nor centering with the expectation works even if the expectation does exist.
More precisely, without centering in case of \ $\alpha\in(0,1)$ \ or with centering with the expectation in case of \ $\alpha\in(1,2)$,
 \ we have the following convergences.

\begin{Cor}\label{simple_aggregation1_stable_centering_fd}
In case of \ $\alpha \in (0, 1)$, \ for each \ $k \in \ZZ_+$, \ we have
 \[
  \biggl(\frac{1}{a_N}
  \sum_{j=1}^{\lfloor Nt\rfloor} \bigl(X^{(j)}_0, \ldots, X^{(j)}_k\bigr)^\top\biggr)_{t\in\RR_+}
  \distr \Bigl(\bcX_t^{(k,\alpha)} + \frac{\alpha}{1-\alpha} t \bone_{k+1}\Bigr)_{t\in\RR_+}
 \]
 as \ $N\to\infty$, \ and, in case of \ $\alpha \in (1, 2)$, \ for each \ $k \in \ZZ_+$, \ we have
 \begin{equation}\label{help6_simple_aggregation1_stable_centering2_fd}
  \begin{aligned}
   &\biggl(\frac{1}{a_N} \sum_{j=1}^{\lfloor Nt\rfloor}
   \bigl(X^{(j)}_0
         - \EE(X^{(j)}_0), \ldots, X^{(j)}_k
         - \EE(X^{(j)}_k)\bigr)^\top\biggr)_{t\in\RR_+} \\
   &= \biggl(\frac{1}{a_N} \sum_{j=1}^{\lfloor Nt\rfloor}
   \bigl(X^{(j)}_0, \ldots, X^{(j)}_k\bigr)^\top - \frac{\lfloor Nt\rfloor}{a_N} \EE(X_0) \bone_{k+1}\biggr)_{t\in\RR_+}
  \distr \Bigl(\bcX_t^{(k,\alpha)} + \frac{\alpha}{1-\alpha} t \bone_{k+1}\Bigr)_{t\in\RR_+}
  \end{aligned}
 \end{equation}
 as \ $N\to\infty$.
\ Moreover, \ $\bigl(\bcX_t^{(k,\alpha)} + \frac{\alpha}{1-\alpha} t \bone_{k+1}\bigr)_{t\in\RR_+}$ \ is a \ $(k + 1)$-dimensional \ $\alpha$-stable process such that the characteristic function of \ $\bcX_1^{(k,\alpha)} + \frac{\alpha}{1-\alpha} \bone_{k+1}$ \ has the form
 \begin{align*}
  &\EE\Bigl(\exp\Bigl\{\ii\Bigl\langle\btheta,\bcX_1^{(k,\alpha)} + \frac{\alpha}{1-\alpha} \bone_{k+1}\Bigr\rangle\Bigr\}\Bigr) \\
  &= \begin{cases}
      \exp\biggl\{(1 - m_\xi^\alpha)
                  \sum_{j=0}^k
                   \int_0^\infty
                    \bigl(\ee^{\ii\langle\btheta,\bv_j^{(k)}\rangle u} - 1\bigr)
                    \alpha u^{-1-\alpha} \, \dd u\biggr\} ,
       &\text{if \ $\alpha \in (0, 1)$,} \\[3mm]
      \exp\biggl\{(1 - m_\xi^\alpha)
                  \sum_{j=0}^k
                   \int_0^\infty
                    \bigl(\ee^{\ii\langle\btheta,\bv_j^{(k)}\rangle u} - 1 - \ii \langle\btheta,\bv_j^{(k)}\rangle u\bigr)
                    \alpha u^{-1-\alpha} \, \dd u\biggr\} ,
       &\text{if \ $\alpha \in (1, 2)$,}
     \end{cases} \\
  &=\exp\biggl\{- C_\alpha (1 - m_\xi^\alpha)
                  \sum_{j=0}^k
                   |\langle\btheta, \bv_j^{(k)}\rangle|^\alpha
                   \left(1 - \ii \tan\left(\frac{\pi\alpha}{2}\right)
                             \sign(\langle\btheta, \bv_j^{(k)}\rangle)\right)\biggr\} , \quad \text{if \ $\alpha \ne 1$,}
 \end{align*}
 for \ $\btheta \in \RR^{k+1}$.
\end{Cor}

Note that in case of \ $\alpha \in (1, 2)$, \ the scaling and the centering in \eqref{help6_simple_aggregation1_stable_centering2_fd} do not depend on \ $j$ \ or \ $k$, \ since the copies are independent and the process \ $(X_k)_{k\in\ZZ_+}$ \ is strongly stationary, and especially, \ $\EE\bigl(X^{(j)}_k\bigr) = \EE(X_0) = \frac{m_\vare}{1-m_\xi}$ \ for all \ $j \in \NN$ \ and \ $k \in \ZZ_+$ \ with \ $m_\vare := \EE(\vare)$,
 \ see, e.g., Barczy et al.\ \cite[formula (14)]{BarNedPap2}.

The next remark is devoted to study some distributional properties of the \ $\alpha$-stable process \ $\bigl(\bcX_t^{(k,\alpha)} + \frac{\alpha}{1-\alpha} t\bone_{k+1}\bigr)_{t\in\RR_+}$ \ in case of \ $\alpha \ne 1$.

\begin{Rem}\label{varrho_properties}
The L\'evy measure of the distribution of \ $\bcX_1^{(k,\alpha)} + \frac{\alpha}{1-\alpha} \bone_{k+1}$ \ is the same as that of \ $\bcX_1^{(k,\alpha)}$, \ namely,
 \ $\nu_{k,\alpha}$ \ given in Remark \ref{mu_properties1}.

If \ $\alpha \in (0, 1)$, \ then the drift of the distribution of \ $\bcX_1^{(k,\alpha)} + \frac{\alpha}{1-\alpha} \bone_{k+1}$ \ is \ $\bzero$, \ hence the process \ $\bigl(\bcX_t^{(k,\alpha)} + \frac{\alpha}{1-\alpha} t\bone_{k+1}\bigr)_{t\in\RR_+}$ \ is strictly \ $\alpha$-stable, see, e.g., Sato
 \cite[Theorem 14.7 (iv)]{Sato}.

If \ $\alpha \in (1, 2)$, \ then the center, i.e., the expectation of
 \ $\bcX_1^{(k,\alpha)} + \frac{\alpha}{1-\alpha} \bone_{k+1}$ \ is \ $\bzero$, \ hence the process \ $\bigl(\bcX_t^{(k,\alpha)} + \frac{\alpha}{1-\alpha} t\bone_{k+1}\bigr)_{t\in\RR_+}$ \ is strictly \ $\alpha$-stable see, e.g., Sato \cite[Theorem 14.7 (vi)]{Sato}.

All in all, \ $\bigl(\bcX_t^{(k,\alpha)} + \frac{\alpha}{1-\alpha} t\bone_{k+1}\bigr)_{t\in\RR_+}$ \ is strictly \ $\alpha$-stable for any
 \ $\alpha \ne 1$.
\ We also note that for each \ $t \in \RR_{++}$, \ the distribution of \ $\bcX_t^{(k,\alpha)} + \frac{\alpha}{1-\alpha} t \bone_{k+1}$ \ is absolutely continuous, see, e.g., Sato \cite[Theorem 27.4 and Proposition 14.5]{Sato}.
\proofend
\end{Rem}

Let \ $\bigl(\cY^{(\alpha)}_k\bigr)_{k\in\ZZ_+}$ \ be a strongly stationary process such that
 \begin{align}\label{calY_def}
    \bigl(\cY^{(\alpha)}_k\bigr)_{k\in\{0,\ldots,K\}} \distre \bcX_1^{(K,\alpha)} \qquad \text{for each \ $K \in \ZZ_+$.}
 \end{align}
The existence of \ $\bigl(\cY^{(\alpha)}_k\bigr)_{k\in\ZZ_+}$ \ follows from the Kolmogorov extension theorem.
Its strong stationarity is a consequence of \eqref{help6_simple_aggregation1_stable_fd} together with the strong stationarity of \ $(X_k)_{k\in\ZZ_+}$.
\ We note that the common distribution of \ $\cY^{(\alpha)}_k$, \ $k \in \ZZ_+$, \ depends only on \ $\alpha$, \ it does not depend on \ $m_\xi$, \ since its characteristic function has the form
 \begin{align*}
  &\EE\bigl(\ee^{\ii\vartheta\cY^{(\alpha)}_0}\bigr)
   = \EE\bigl(\ee^{\ii\vartheta\bcX_1^{(0,\alpha)}}\bigr) \\
  &= \exp\biggl\{(1 - m_\xi^\alpha)
                 \int_0^\infty
                  \Bigl(\ee^{\ii\vartheta(1-m_\xi^\alpha)^{-\frac{1}{\alpha}}u} - 1
                         - \ii u \vartheta (1-m_\xi^\alpha)^{-\frac{1}{\alpha}} \bbone_{(0,1]}(u (1-m_\xi^\alpha)^{-\frac{1}{\alpha}})\Bigr) \alpha
                 u^{-1-\alpha} \, \dd u\biggr\} \\
  &= \exp\biggl\{\int_0^\infty
                  \bigl(\ee^{\ii\vartheta v} - 1
                         - \ii \vartheta v \bbone_{(0,1]}(v)\bigr) \alpha
                 v^{-1-\alpha} \, \dd v\biggr\} , \qquad \vartheta \in \RR .
 \end{align*}

\begin{Pro}\label{Pro_AR1}
For each \ $\alpha \in (0, 2)$, \ the strongly stationary process \ $\bigl(\cY^{(\alpha)}_k\bigr)_{k\in\ZZ_+}$ \ is a subcritical
 autoregressive process of order 1 with autoregressive coefficient \ $m_\xi$ \ and with \ $\alpha$-stable innovations,
 namely,
 \[
   \cY^{(\alpha)}_k = m_\xi \cY^{(\alpha)}_{k-1} + \tvare^{(\alpha)}_k , \qquad k \in \NN ,
 \]
 where
 \[
   \tvare^{(\alpha)}_k := \cY^{(\alpha)}_k - m_\xi \cY^{(\alpha)}_{k-1}, \qquad k \in \NN,
 \]
 is a sequence of independent, identically distributed \ $\alpha$-stable random variables such that
 for all \ $k\in\NN$, \ $\tvare^{(\alpha)}_k$ \ is independent of \ $(\cY^{(\alpha)}_0, \ldots, \cY^{(\alpha)}_{k-1})^\top$.
 \ Therefore, \ $\bigl(\cY^{(\alpha)}_k\bigr)_{k\in\ZZ_+}$ \ is a strongly stationary, time homogeneous Markov process.
\end{Pro}

Theorem \ref{simple_aggregation1_stable_fd} and Corollary \ref{simple_aggregation1_stable_centering_fd} have the following consequences
 for a contemporaneous aggregation of independent copies with different centerings.

\begin{Cor}\label{aggr_copies1}
\renewcommand{\labelenumi}{{\rm(\roman{enumi})}}
 \begin{enumerate}
  \item
   For each \ $\alpha \in (0, 2)$,
   \begin{align*}
    \biggl(\frac{1}{a_N}
          \sum_{j=1}^N
           \Bigl(X^{(j)}_k
                 - \EE\bigl(X^{(j)}_k
                            \bbone_{\{X^{(j)}_k\leq a_N\}}\bigr)
           \Bigr)
   \biggr)_{k\in\ZZ_+}
   &= \biggl(\frac{1}{a_N}
          \sum_{j=1}^N
           X^{(j)}_k
           - \frac{N}{a_N} \EE\bigl(X_0
                            \bbone_{\{X_0\leq a_N\}}\bigr)
   \biggr)_{k\in\ZZ_+} \\
    &\distrf
   \bigl(\cY^{(\alpha)}_k\bigr)_{k\in\ZZ_+} \qquad \text{as \ $N \to \infty$,}
 \end{align*}
  \item
   in case of \ $\alpha \in (0, 1)$,
  \[
   \biggl(\frac{1}{a_N} \sum_{j=1}^N X^{(j)}_k \biggr)_{k\in\ZZ_+}
   \distrf
   \Bigl(\cY^{(\alpha)}_k + \frac{\alpha}{1-\alpha}\Bigr)_{k\in\ZZ_+} \qquad \text{as \ $N \to \infty$,}
 \]
  \item
   in case of \ $\alpha \in (1, 2)$,
  \begin{align*}
   \biggl(\frac{1}{a_N} \sum_{j=1}^N \bigl(X^{(j)}_k - \EE\bigl(X^{(j)}_k \bigr)\bigr)\biggr)_{k\in\ZZ_+}
   &= \biggl(\frac{1}{a_N} \sum_{j=1}^N X^{(j)}_k - \frac{N}{a_N} \EE\bigl(X_0\bigr) \biggr)_{k\in\ZZ_+} \\
   &\distrf
   \Bigl(\cY^{(\alpha)}_k + \frac{\alpha}{1-\alpha} \Bigr)_{k\in\ZZ_+} \qquad \text{as \ $N \to \infty$,}
 \end{align*}
 \end{enumerate}
 where \ $(\cY^{(k)})_{k\in\ZZ_+}$ \ is given by \eqref{calY_def}.
\end{Cor}

Limit theorems will be presented for the aggregated stochastic process \ $\bigl(\sum_{k=1}^{\lfloor nt \rfloor} \sum_{j=1}^N X^{(j)}_k\bigr)_{t\in\RR_+}$ \
 with different centerings and scalings.
We will provide limit theorems in an iterated manner such that first $N$, and then $n$ converges to infinity.

\begin{Thm}\label{iterated_aggr_1}
In case of \ $\alpha \in (0, 1)$, \ we have
 \begin{equation}\label{iterated_aggr_1_1}
 \begin{aligned}
  &\cD_\ff\text{-}\hspace*{-1mm}\lim_{n\to\infty} \,
   \cD_\ff\text{-}\hspace*{-1mm}\lim_{N\to\infty} \,
    \biggl(\frac{1}{n^{\frac{1}{\alpha}}a_N}
           \sum_{k=1}^\nt \sum_{j=1}^N
            \Bigl(X^{(j)}_k
                  - \EE\bigl(X^{(j)}_k
                             \bbone_{\{X^{(j)}_k\leq a_N\}}\bigr)
            \Bigr)
   \biggr)_{t\in\RR_+} \\
  &=\cD_\ff\text{-}\hspace*{-1mm}\lim_{n\to\infty} \,
   \cD_\ff\text{-}\hspace*{-1mm}\lim_{N\to\infty} \,
    \biggl(\frac{1}{n^{\frac{1}{\alpha}}a_N}
           \sum_{k=1}^\nt \sum_{j=1}^N X^{(j)}_k
           - \frac{\nt N}{n^{\frac{1}{\alpha}}a_N}
             \EE\bigl(X_0 \bbone_{\{X_0\leq a_N\}}\bigr)\biggr)_{t\in\RR_+} \\
  &= \Bigl(\cZ_t^{(\alpha)} + \frac{\alpha}{1-\alpha} t\Bigr)_{t\in\RR_+} ,
 \end{aligned}
 \end{equation}
 and
 \begin{equation}\label{iterated_aggr_1_2}
   \cD_\ff\text{-}\hspace*{-1mm}\lim_{n\to\infty} \,
   \cD_\ff\text{-}\hspace*{-1mm}\lim_{N\to\infty} \,
    \biggl(\frac{1}{n^{\frac{1}{\alpha}}a_N}
           \sum_{k=1}^\nt \sum_{j=1}^N
            X^{(j)}_k
            \Bigr)
   \biggr)_{t\in\RR_+}
   = \Bigl(\cZ_t^{(\alpha)} + \frac{\alpha}{1-\alpha} t\Bigr)_{t\in\RR_+} ,
 \end{equation}
 in case of \ $\alpha = 1$, \ we have
 \begin{equation}\label{iterated_aggr_1_3}
 \begin{aligned}
  &\cD_\ff\text{-}\hspace*{-1mm}\lim_{n\to\infty} \,
   \cD_\ff\text{-}\hspace*{-1mm}\lim_{N\to\infty} \,
    \biggl(\frac{1}{n\log(n)a_N}
           \sum_{k=1}^\nt \sum_{j=1}^N
            \Bigl(X^{(j)}_k
                  - \EE\bigl(X^{(j)}_k
                             \bbone_{\{X^{(j)}_k\leq a_N\}}\bigr)
            \Bigr)
   \biggr)_{t\in\RR_+} \\
  &=\cD_\ff\text{-}\hspace*{-1mm}\lim_{n\to\infty} \,
   \cD_\ff\text{-}\hspace*{-1mm}\lim_{N\to\infty} \,
    \biggl(\frac{1}{n\log(n)a_N}
           \sum_{k=1}^\nt \sum_{j=1}^N X^{(j)}_k
           - \frac{\nt N}{n\log(n)a_N}
             \EE\bigl(X_0 \bbone_{\{X_0\leq a_N\}}\bigr)\biggr)_{t\in\RR_+} \\
  &= (t)_{t\in\RR_+} ,
 \end{aligned}
 \end{equation}
 and in case of \ $\alpha \in (1, 2)$, \ we have
 \begin{equation}\label{iterated_aggr_1_4}
 \begin{aligned}
  &\cD_\ff\text{-}\hspace*{-1mm}\lim_{n\to\infty} \,
   \cD_\ff\text{-}\hspace*{-1mm}\lim_{N\to\infty} \,
    \biggl(\frac{1}{n^{\frac{1}{\alpha}}a_N}
           \sum_{k=1}^\nt \sum_{j=1}^N
            (X^{(j)}_k - \EE(X^{(j)}_k))
            \Bigr)
   \biggr)_{t\in\RR_+} \\
  &=\cD_\ff\text{-}\hspace*{-1mm}\lim_{n\to\infty} \,
   \cD_\ff\text{-}\hspace*{-1mm}\lim_{N\to\infty} \,
    \biggl(\frac{1}{n^{\frac{1}{\alpha}}a_N}
           \sum_{k=1}^\nt \sum_{j=1}^N X^{(j)}_k
           - \frac{\nt N}{n^{\frac{1}{\alpha}}a_N} \EE(X_0)
   \biggr)_{t\in\RR_+} \\
  &= \Bigl(\cZ_t^{(\alpha)} + \frac{\alpha}{1-\alpha} t\Bigr)_{t\in\RR_+} ,
 \end{aligned}
 \end{equation}
\noindent
 where \ $\bigl(\cZ_t^{(\alpha)}\bigr)_{t\in\RR_+}$ \ is an \ $\alpha$-stable process such that the characteristic function of the distribution of \ $\cZ_1^{(\alpha)}$ \ has the form
 \[
   \EE\bigl(\ee^{\ii\vartheta\cZ_1^{(\alpha)}}\bigr)
   = \exp\biggl\{\ii b_\alpha \vartheta + \frac{1-m_\xi^\alpha}{(1-m_\xi)^\alpha}
                 \int_0^\infty
                  (\ee^{\ii\vartheta u} - 1
                   - \ii \vartheta u \bbone_{(0,1]}(u)) \alpha
                  u^{-1-\alpha} \, \dd u\biggr\} , \qquad \vartheta \in \RR ,
 \]
 where
 \[
   b_\alpha := \biggl(\frac{1-m_\xi^\alpha}{(1-m_\xi)^\alpha} - 1\biggr) \frac{\alpha}{1-\alpha} , \qquad \alpha \in (0, 1) \cup (1, 2) ,
 \]
 and \ $\bigl(\cZ_t^{(\alpha)} + \frac{\alpha}{1-\alpha} t\bigr)_{t\in\RR_+}$ \ is an \ $\alpha$-stable process such that the characteristic function of the distribution of \ $\cZ_1^{(\alpha)} + \frac{\alpha}{1-\alpha}$ \ has the form
 \begin{align*}
  &\EE\Bigl(\exp\Bigl\{\ii\vartheta\Bigl(\cZ_1^{(\alpha)}+\frac{\alpha}{1-\alpha}\Bigr)\Bigr\}\Bigr)
   = \begin{cases}
      \exp\Bigl\{\frac{1-m_\xi^\alpha}{(1-m_\xi)^\alpha}
                 \int_0^\infty
                  (\ee^{\ii\vartheta u} - 1) \alpha
                  u^{-1-\alpha} \, \dd u \Bigr\} ,
       & \text{if \ $\alpha \in (0, 1)$,} \\[3mm]
      \exp\Bigl\{\frac{1-m_\xi^\alpha}{(1-m_\xi)^\alpha}
                 \int_0^\infty
                  (\ee^{\ii\vartheta u} - 1
                   - \ii \vartheta u) \alpha
                  u^{-1-\alpha} \, \dd u \Bigr\} ,
       & \text{if \ $\alpha \in (1, 2)$,}
     \end{cases} \\
  &= \exp\biggl\{- C_\alpha \frac{1-m_\xi^\alpha}{(1-m_\xi)^\alpha}
                   |\vartheta|^\alpha
                    \left(1 - \ii \tan\left(\frac{\pi\alpha}{2}\right)
                                  \sign(\vartheta)\right)\biggr\} \quad \text{if \ $\alpha \in (0, 1) \cup (1, 2)$,}
 \end{align*}
 for \ $\vartheta \in \RR$.
\end{Thm}

\begin{Rem}\label{Rem_char_spectral_proc}
Note that, in accordance with Basrak and Segers \cite[Remark 4.8]{BasSeg} and Mikosch and Wintenberger \cite[page 171]{MikWin},
 in case of \ $\alpha \in (0, 1)$, \ we have
 \begin{align}\label{help_kar_theta}
 \begin{split}
 &\EE\Bigl(\exp\Bigl\{\ii \vartheta \Bigl(\cZ_1^{(\alpha)} + \frac{\alpha}{1-\alpha}\Bigr)\Bigr\}\Bigr) \\
 &= \exp\left\{ - \int_0^\infty \EE\left[ \exp\Big( \ii u \vartheta \sum_{\ell=1}^\infty \Theta_\ell \Big)
                               - \exp\Big( \ii u \vartheta \sum_{\ell=0}^\infty \Theta_\ell \Big) \right] \alpha u^{-\alpha-1}\,\dd u \right\}
 \end{split}
 \end{align}
 for \ $\vartheta \in \RR$, \ where \ $(\Theta_\ell)_{\ell\in\ZZ_+}$ \ is the (forward) spectral tail process of \ $(X_\ell)_{\ell\in\ZZ_+}$ \ given in \eqref{spectral_tail_process} and \eqref{spectral_tail_process0}.
Indeed, by \eqref{help_Z1},
 \begin{align*}
  &\exp\left\{ - \int_0^\infty  \EE\left[ \exp\Big( \ii u \vartheta \sum_{\ell=1}^\infty \Theta_\ell \Big)
                               - \exp\Big( \ii u \vartheta \sum_{\ell=0}^\infty \Theta_\ell \Big) \right] \alpha u^{-\alpha-1}\,\dd u \right\}  \\
   & = \exp\left\{ - \int_0^\infty  \EE\left[ \exp\Big( \ii u \vartheta \sum_{\ell=1}^\infty m_\xi^\ell \Big)
                               - \exp\Big( \ii u \vartheta \sum_{\ell=0}^\infty m_\xi^\ell \Big) \right] \alpha u^{-\alpha-1}\,\dd u \right\} \\
   & = \exp\left\{ - \int_0^\infty  \left( \exp\Big( \ii u \vartheta \frac{m_\xi}{1-m_\xi} \Big)
                               - \exp\Big( \ii u \vartheta \frac{1}{1-m_\xi} \Big) \right) \alpha u^{-\alpha-1}\,\dd u \right\} \\
   &= \exp\left\{ - \int_0^\infty  \left( \exp\Big( \ii u \frac {\vartheta m_\xi}{1-m_\xi} \Big) - 1 \right) \alpha u^{-\alpha-1}\,\dd u
                        +  \int_0^\infty  \left( \exp\Big( \ii u \frac{\vartheta}{1-m_\xi} \Big) - 1 \right) \alpha u^{-\alpha-1}\,\dd u \right\} \\
   & = \exp\Bigg\{ C_\alpha \left\vert \frac{\vartheta m_\xi}{1-m_\xi} \right\vert^\alpha
                   \left(1 - \ii \tan\left(\frac{\pi\alpha}{2}\right)
                                  \sign\left(\frac{\vartheta m_\xi}{1-m_\xi}\right)\right)  \\
   &\phantom{= \exp\Bigg\{ \;}
                   - C_\alpha \left\vert \frac{\vartheta }{1-m_\xi} \right\vert^\alpha
                   \left(1 - \ii \tan\left(\frac{\pi\alpha}{2}\right)
                                  \sign\left(\frac{\vartheta}{1-m_\xi}\right)\right) \Bigg\} \\
   & = \exp\left\{ - C_\alpha \frac{1-m_\xi^\alpha}{(1-m_\xi)^\alpha} \vert\vartheta\vert^\alpha
                        \left(1 - \ii \tan\left(\frac{\pi\alpha}{2}\right)
                                  \sign\left(\frac{\vartheta}{1-m_\xi}\right)\right) \right\} ,
 \end{align*}
 as desired.
We also remark that, using \eqref{help_Z2}, one can check that \eqref{help_kar_theta} does not hold in case of \ $\alpha\in(1,2)$, \ which is somewhat unexpected in view of page 171 in Mikosch and Wintenberger \cite{MikWin}.
\proofend
\end{Rem}

\begin{Rem}\label{cZ_properties}
If \ $\alpha \in (0, 1)$, \ then the drift of the distribution of \ $\cZ_1^{(\alpha)} + \frac{\alpha}{1-\alpha}$ \ is \ $0$, \ hence the process \ $\bigl(\cZ_t^{(\alpha)} + \frac{\alpha}{1-\alpha} t\bigr)_{t\in\RR_+}$ \ is strictly \ $\alpha$-stable, see, e.g., Sato \cite[Theorem 14.7 (iv) and Definition 13.2]{Sato}.

If \ $\alpha \in (1, 2)$, \ then the center, i.e., the expectation of
 \ $\cZ_1^{(\alpha)} + \frac{\alpha}{1-\alpha}$ \ is \ $0$, \ hence the process
 \ $\bigl(\cZ_t^{(\alpha)} + \frac{\alpha}{1-\alpha} t\bigr)_{t\in\RR_+}$ \ is strictly \ $\alpha$-stable see, e.g., Sato \cite[Theorem 14.7 (vi) and Definition 13.2]{Sato}.

All in all, the process \ $\bigl(\cZ_t^{(\alpha)} + \frac{\alpha}{1-\alpha}t\bigr)_{t\in\RR_+}$ \ is strictly \ $\alpha$-stable for any
 \ $\alpha \ne 1$.
\proofend
\end{Rem}

\section{Proofs}
\label{Proofs}

\noindent{\bf Proof of Theorem \ref{simple_aggregation1_stable_fd}.}
Let \ $k \in \ZZ_+$.
\ We are going to apply Theorem \ref{7.1} with \ $d = k + 1$ \ and
 \ $\bX_{N,j} := a_N^{-1} (X_0^{(j)}, \ldots, X_k^{(j)})^\top$, \ $N, j \in \NN$.
\ The aim of the following discussion is to check condition \eqref{(7.5)} of Theorem \ref{7.1}, namely
 \begin{equation}\label{vague_fd}
  N \PP(\bX_{N,1} \in \cdot)
  = N \PP\bigl(a_N^{-1} (X_0^{(1)}, \ldots, X_k^{(1)})^\top \in \cdot\bigr)
  \distrv \nu_{k,\alpha}(\cdot) \qquad \text{on \ $\RR_0^{k+1}$ \ as \ $N \to \infty$,}
 \end{equation}
 where \ $\nu_{k,\alpha}$ \ is a L\'evy measure on \ $\RR_0^{k+1}$.
\ For each \ $N \in \NN$ \ and \ $B \in \cB(\RR_0^{k+1})$, \ we can write
 \[
   N \PP(\bX_{N,1} \in B)
   = N \PP(X_0 > a_N)
     \frac{\PP(a_N^{-1} (X_0, \ldots, X_k)^\top \in B)}
          {\PP(X_0 > a_N)} .
 \]
By the assumption, we have \ $N \PP(X_0 > a_N) \to 1$ \ as \ $N \to \infty$, \ yielding also \ $a_N \to \infty$ \ as \ $N \to \infty$, \ consequently, it is enough to show that
 \begin{equation}\label{nu_alpha^k_old}
  \frac{\PP(x^{-1} (X_0, \ldots, X_k)^\top \in \cdot)}{\PP(X_0 > x)}
  \distrv \nu_{k,\alpha}(\cdot)  \qquad
  \text{on \ $\RR_0^{k+1}$ \ as \ $x \to \infty$,}
 \end{equation}
 where \ $\nu_{k,\alpha}$ \ is a L\'evy measure on \ $\RR_0^{k+1}$.
\ In fact, by Theorem \ref{Xtailprocess}, \ $(X_0, \ldots, X_k)^\top$ \ is regularly varying with index \ $\alpha$, \ hence, by Proposition \ref{vague}, we know that
 \begin{equation}\label{tnu_alpha^k_old}
  \frac{\PP(x^{-1} (X_0, \ldots, X_k)^\top \in \cdot)}{\PP(\|(X_0, \ldots, X_k)^\top\| > x)}
  \distrv \tnu_{k,\alpha}(\cdot)  \qquad
  \text{on \ $\RR_0^{k+1}$ \ as \ $x \to \infty$,}
 \end{equation}
 where \ $\tnu_{k,\alpha}$ \ is the so-called limit measure of \ $(X_0, \ldots, X_k)^\top$.
\ Applying Proposition \ref{Pro_mapping} for the canonical projection \ $p_0 : \RR^{k+1} \to \RR$ \ given by \ $p_0(\bx) := x_0$ \ for \ $\bx = (x_0, \ldots, x_k)^\top \in \RR^{k+1}$, \ which is continuous and positively homogeneous of degree 1, we obtain
 \[
  \frac{\PP(X_0 > x)}{\PP(\|(X_0, \ldots, X_k)^\top\| > x)}
  \to \tnu_{k,\alpha}(T_1) \qquad \text{as \ $x \to \infty$,}
 \]
 with \ $T_1 := \{\bx \in \RR^{k+1}_0 : p_0(\bx) > 1\}$, \ where we have \ $\tnu_{k,\alpha}(T_1) \in (0, 1]$.
\ Indeed, \ $\PP(X_0 > x) \leq \PP(\|(X_0, \ldots, X_k)^\top\| > x)$, \ hence \ $\tnu_{k,\alpha}(T_1) \leq 1$.
\ Moreover, by the strong stationarity of \ $(X_k)_{k\in\ZZ_+}$, \ we have
 \[
   \PP(\|(X_0, \ldots, X_k)^\top\| > x)
   \leq \sum_{j=0}^k \PP(X_j > x / \sqrt{k+1})
   = (k + 1) \PP(X_0 > x / \sqrt{k+1}) ,
 \]
 thus
 \[
   \frac{\PP(X_0 > x)}{\PP(\|(X_0, \ldots, X_k)^\top\| > x)}
   \geq \frac{\PP(X_0 > x)}{(k+1)\PP(X_0 > x / \sqrt{k+1})}
   \to (k + 1)^{-1-\frac{\alpha}{2}} \qquad \text{as \ $x \to \infty$,}
 \]
 since \ $X_0$ \ is regularly varying with index \ $\alpha$, \ hence \ $\tnu_{k,\alpha}(T_1) \in (0, 1]$, \ as desired.
Consequently, \eqref{nu_alpha^k_old} holds with \ $\nu_{k,\alpha} = \tnu_{k,\alpha}/\tnu_{k,\alpha}(T_1)$.
\ In general, one does not know whether \ $\nu_{k,\alpha}$ \ is a L\'evy measure on \ $\RR_0^{k+1}$ \ or not.
So, additional work is needed.
We will determine \ $\nu_{k,\alpha}$ \ explicitly, using a result of Planini\'{c} and Soulier \cite{PlaSou}.

The aim of the following discussion is to apply Theorem 3.1 in Planini\'{c} and Soulier \cite{PlaSou} in order to determine \ $\nu_{k,\alpha}$, \ namely, we will prove that for each Borel measurable function \ $f : \RR_0^{k+1} \to \RR_+$,
 \begin{equation}\label{fint}
  \int_{\RR_0^{k+1}} f(\bx) \, \nu_{k,\alpha}(\dd\bx)
   = (1 - m_\xi^\alpha) \sum_{j=0}^k \int_0^\infty f(u\bv_j^{(k)}) \alpha u^{-\alpha-1} \, \dd u .
 \end{equation}
Let \ $(X_\ell)_{\ell\in\ZZ}$ \ be a strongly stationary extension of \ $(X_\ell)_{\ell\in\ZZ_+}$.
\ For each \ $i, j \in \ZZ$ \ with \ $i \leq j$, \ by Theorem \ref{Xtailprocess}, \ $(X_i, \ldots, X_j)^\top$ \ is regularly varying with index \ $\alpha$, \ hence, by the strong stationarity of \ $(X_k)_{k\in\ZZ}$ \ and the discussion above, we know that
 \[
   \frac{\PP(x^{-1} (X_i, \ldots, X_j)^\top \in \cdot)}{\PP(X_0 > x)}
   = \frac{\PP(x^{-1} (X_0, \ldots, X_{j-i})^\top \in \cdot)}{\PP(X_0 > x)}
  \distrv \nu_{i,j,\alpha}(\cdot) \qquad
  \text{on \ $\RR_0^{j-i+1}$ \ as \ $x \to \infty$,}
 \]
 where \ $\nu_{i,j,\alpha} := \nu_{j-i,\alpha}$ \ is a non-null locally finite measure on \ $\RR_0^{j-i+1}$.
\ According to Basrak and Segers \cite[Theorem 2.1]{BasSeg}, there exists a sequence \ $(Y_\ell)_{\ell\in\ZZ}$ \ of random variables, called the (whole) tail process of \ $(X_\ell)_{\ell\in\ZZ}$, \ such that
 \[
   \PP(x^{-1} (X_i, \ldots, X_j)^\top \in \cdot \mid X_0 > x)
   \distrw \PP( (Y_i, \ldots, Y_j)^\top \in \cdot) \qquad \text{as \ $x \to \infty$.}
 \]
Let \ $K$ \ be a random variable with geometric distribution
 \[
   \PP(K = k) = m_\xi^{\alpha k} (1 - m_\xi^\alpha) , \qquad k \in \ZZ_+ .
 \]
Especially, if \ $m_\xi = 0$, \ then \ $\PP(K = 0) = 1$.
\ If \ $m_\xi \in (0, 1)$, \ then we have
 \begin{equation}\label{tail_process}
  Y_\ell = \begin{cases}
         m_\xi^\ell Y_0 , & \text{if \ $\ell \geq 0$,} \\
         m_\xi^\ell Y_0 \bbone_{\{K\geq-\ell\}} , & \text{if \ $\ell < 0$,}
        \end{cases}
 \end{equation}
 where \ $Y_0$ \ is a random variable independent of \ $K$ \ with Pareto distribution
 \[
   \PP(Y_0 > y) = \begin{cases}
                   y^{-\alpha} , & \text{if \ $y \in [1, \infty)$,} \\
                   1 , & \text{if \ $y \in (-\infty, 1)$.}
                  \end{cases}
 \]
Indeed, as shown in Basrak et al.~\cite[Lemma 3.1]{BasKulPal},  $(Y_\ell)_{\ell\in\ZZ_+}$ is the forward tail
 process of \ $(X_\ell)_{\ell\in\ZZ}$.
\ On the other hand, by Janssen and Segers \cite[Example 6.2]{JanSeg}, \ $(Y_\ell)_{\ell\in\ZZ}$ \ is the tail process of
 the stationary solution \ $(X_\ell')_{\ell\in \ZZ}$ \ to the stochastic recurrence equation \ $X_\ell'=\mu_A X_{\ell-1}' + B_\ell$, \ $\ell \in \ZZ$.
\ Since the distribution of the forward tail process determines the distribution of the (whole) tail process
 (see Basrak and Segers~\cite[Theorem 3.1 (ii)]{BasSeg}), it follows that \ $(Y_\ell)_{\ell\in\ZZ}$ \ represents the tail process of
 \ $(X_\ell)_{\ell\in\ZZ}$.
\ If \ $m_\xi = 0$, \ then one can easily check that
 \begin{equation}\label{tail_process0}
  Y_\ell = \begin{cases}
            Y_0 , & \text{if \ $\ell = 0$,} \\
            0 , & \text{if \ $\ell \ne 0$.}
           \end{cases}
 \end{equation}
By \eqref{tail_process} and \eqref{tail_process0}, we have \ $Y_\ell \as 0$ \ as \ $\ell \to \infty$ \ or \ $\ell \to -\infty$,
 \ hence condition (3.1) in Planini\'{c} and Soulier \cite{PlaSou} is satisfied.

Moreover, there exists a unique measure \ $\nu_\alpha$ \ on \ $\RR^\ZZ$ \ endowed with the cylindrical $\sigma$-algebra \ $\cB(\RR)^{\otimes\ZZ}$ \ such that \ $\nu_\alpha(\{\bzero\}) = 0$ \ and for each \ $i, j \in \ZZ$ \ with \ $i \leq j$, \ we have \ $\nu_\alpha \circ p_{i,j}^{-1} = \nu_{i,j,\alpha}$ \ on \ $\RR^{j-i+1}_0$, \ where \ $p_{i,j}$ \ denotes the canonical projection \ $p_{i,j} : \RR^\ZZ \to \RR^{j-i+1}$ \ given by \ $p_{i,j}(\by) := (y_i, \ldots, y_j)$ \ for \ $\by = (y_\ell)_{\ell\in\ZZ} \in \RR^\ZZ$, \ see, e.g., Planini\'{c} and Soulier \cite{PlaSou}.
The measure \ $\nu_\alpha$ \ is called the tail measure of \ $(X_\ell)_{\ell\in\ZZ}$.

If \ $m_\xi \in (0, 1)$, \ then, by \eqref{tail_process}, the (whole) spectral tail process \ $\bTheta = (\Theta_\ell)_{\ell\in\ZZ}$ \ of \ $(X_\ell)_{\ell\in\ZZ}$ \ is given by
 \begin{equation}\label{spectral_tail_process}
   \Theta_\ell := \frac{Y_\ell}{|Y_0|}
   = \begin{cases}
      m_\xi^\ell , & \text{if \ $\ell \geq 0$,} \\
      m_\xi^\ell \bbone_{\{K\geq-\ell\}} , & \text{if \ $\ell < 0$.}
     \end{cases}
 \end{equation}
If \ $m_\xi = 0$, \ then, by \eqref{tail_process0},
 \begin{equation}\label{spectral_tail_process0}
   \Theta_\ell := \frac{Y_\ell}{|Y_0|}
   = \begin{cases}
      1 , & \text{if \ $\ell = 0$,} \\
      0 , & \text{if \ $\ell \ne 0$.}
     \end{cases}
 \end{equation}
Let us introduce the so called infargmax functional \ $I : \RR^\ZZ \to \ZZ \cup \{-\infty, \infty\}$.
\ For \ $\by = (y_\ell)_{\ell\in\ZZ} \in \RR^\ZZ$, \ the value \ $I(\by)$ \ is the first time when the supremum
 \ $\sup_{\ell\in\ZZ} |y_\ell|$ \ is achieved, more precisely,
 \[
   I(\by) := \begin{cases}
              \ell \in \ZZ , & \text{if \ $\sup\limits_{m\leq\ell-1} |y_m| < |y_\ell|$ \ and \ $\sup\limits_{m\geq\ell+1} |y_m| \leq |y_\ell|$,} \\
              -\infty , & \text{if \ $\sup\limits_{m\leq\ell} |y_m| = \sup\limits_{m\in\ZZ} |y_m|$ \ for all \ $\ell \in \ZZ$,} \\
              \infty , & \text{if \ $\sup\limits_{m\leq\ell} |y_m| < \sup\limits_{m\in\ZZ} |y_m|$ \ for all \ $\ell \in \ZZ$.}
             \end{cases}
 \]
We have \ $\PP(I(\bTheta) = -K) = 1$, \ hence the condition \ $\PP(I(\bTheta) \in \ZZ) = 1$ \ of Theorem 3.1 in Planini\'{c} and Soulier \cite{PlaSou} is satisfied.

Consequently, we may apply Theorem 3.1 in Planini\'{c} and Soulier \cite{PlaSou} for the nonnegative measurable function \ $H : \RR^\ZZ \to \RR_+$ \ given by \ $H(\by) = f \circ p_{0,k}$, \ where \ $f : \RR^{k+1} \to \RR_+$ \ is a measurable function with \ $f(\bzero) = 0$.
\ By (3.2) in Planini\'{c} and Soulier \cite{PlaSou}, we obtain
 \begin{align*}
  \int_{\RR_0^{k+1}} f(\bx) \, \nu_{k,\alpha}(\dd\bx)
  &= \int_{\RR^{k+1}} f(\bx) \, \nu_{0,k,\alpha}(\dd\bx)
   = \int_{\RR^{k+1}} f(\bx) \, (\nu_\alpha \circ p_{0,k}^{-1})(\dd\bx)
   = \int_{\RR^\ZZ} f(p_{0,k}(\by)) \, \nu_\alpha(\dd\by) \\
  &= \int_{\RR^\ZZ} H(\by) \, \nu_\alpha(\dd\by)
   = \sum_{\ell\in\ZZ}
      \int_0^\infty
       \EE(H(uL^\ell(\bTheta)) \bbone_{\{I(\bTheta)=0\}}) \alpha u^{-\alpha-1} \, \dd u ,
 \end{align*}
 where \ $L$ \ denotes the backshift operator \ $L : \RR^\ZZ \to \RR^\ZZ$ \ given by \ $L(\by) = (L(\by)_k)_{k\in\ZZ} := (y_{k-1})_{k\in\ZZ}$ \ for \ $\by = (y_k)_{k\in\ZZ} \in \RR^\ZZ$.
\ Using \ $\PP(I(\bTheta) = -K)=1$, \ we obtain
 \[
   \int_{\RR_0^{k+1}} f(\bx) \, \nu_{k,\alpha}(\dd\bx)
   = \sum_{\ell\in\ZZ}
      \int_0^\infty
       \EE(f(p_{0,k}(uL^\ell(\bTheta))) \bbone_{\{K=0\}}) \alpha u^{-\alpha-1} \, \dd u .
 \]
For each \ $k \in \ZZ_+$ \ and \ $u \in \RR_+$, \ on the event \ $\{K = 0\}$, \ by \eqref{spectral_tail_process} and \eqref{spectral_tail_process0},
 we have
 \[
   p_{0,k}(uL^\ell(\bTheta))
   = \begin{cases}
      \bzero \in \RR^{k+1} , & \text{if \ $\ell > k$,} \\
      u \bv_\ell^{(k)} , & \text{if \ $\ell \in \{1, \ldots, k\}$,} \\
      (1 - m_\xi^\alpha)^{\frac{1}{\alpha}} m_\xi^{-\ell} u \bv_0^{(k)} , & \text{if \ $\ell \leq 0$,}
     \end{cases}
 \]
 hence, using \ $\PP(K = 0) = 1 - m_\xi^\alpha$, \ we obtain
 \begin{align*}
  &\int_{\RR_0^{k+1}} f(\bx) \, \nu_{k,\alpha}(\dd\bx) \\
  &= (1 - m_\xi^\alpha)
     \sum_{\ell\leq0}
      \int_0^\infty
       f((1 - m_\xi^\alpha)^{\frac{1}{\alpha}} m_\xi^{-\ell} u \bv_0^{(k)}) \alpha u^{-\alpha-1} \, \dd u
     + (1 - m_\xi^\alpha)
       \sum_{\ell=1}^k
        \int_0^\infty
         f(u \bv_\ell^{(k)}) \alpha u^{-\alpha-1} \, \dd u \\
  &= (1 - m_\xi^\alpha)^2
     \sum_{\ell\leq0}
       m_\xi^{-\ell\alpha}
       \int_0^\infty
        f(u \bv_0^{(k)}) \alpha u^{-\alpha-1} \, \dd u
     + (1 - m_\xi^\alpha)
       \sum_{\ell=1}^k
        \int_0^\infty
         f(u \bv_\ell^{(k)}) \alpha u^{-\alpha-1} \, \dd u \\
  &= (1 - m_\xi^\alpha)
      \sum_{\ell=0}^k
       \int_0^\infty
        f(u \bv_\ell^{(k)}) \alpha u^{-\alpha-1} \, \dd u .
 \end{align*}

The measure \ $\nu_{k,\alpha}$ \ is a L\'evy measure on \ $\RR_0^{k+1}$, \ since
 \eqref{fint} implies
 \begin{align*}
  &\int_{\RR_0^{k+1}} \min\{1, \|\bx\|^2\} \, \nu_{k,\alpha}(\dd\bx)
   = (1-m_\xi^\alpha)
     \sum_{j=0}^k
      \int_0^\infty \min\{1, \|u \bv_j^{(k)}\|^2\} \alpha u^{-\alpha-1} \, \dd u \\
  &= (1-m_\xi^\alpha)
     \sum_{j=0}^k
      \|\bv_j^{(k)}\|^\alpha \int_0^\infty \min\{1, w^2\} \alpha w^{-\alpha-1} \, \dd w
  = \frac{2(1-m_\xi^\alpha)}{2-\alpha} \sum_{j=0}^k \|\bv_j^{(k)}\|^\alpha
   < \infty .
 \end{align*}
Consequently, we obtain \eqref{nu_alpha^k_old}, and hence \eqref{vague_fd}, so condition \eqref{(7.5)} is satisfied.

The aim of the following discussion is to check condition \eqref{(7.6)} of Theorem \ref{7.1}, namely
 \begin{equation}\label{Res_cond}
   \lim_{\vare\downarrow0}
    \limsup_{N\to\infty}
     N \EE(a_N^{-2} (X_\ell^{(j)})^2 \bbone_{\{X_\ell^{(j)}\leq a_N\vare\}})
   = \lim_{\vare\downarrow0}
      \limsup_{N\to\infty}
       N \EE(a_N^{-2} X_0^2
             \bbone_{\{X_0\leq a_N\vare\}})
   = 0
 \end{equation}
 for each \ $j\in\NN$ \ and \ $\ell \in \{0, \ldots, k\}$.
\ By Lemma \ref{truncated_moments} with \ $\beta = 2$, \ we have
 \[
   \lim_{x\to\infty}
    \frac{x^2\PP(X_0>x)}{\EE(X_0^2\bbone_{\{X_0\leq x\}})}
   = \frac{2-\alpha}{\alpha} ,
 \]
 hence, for all \ $\vare \in \RR_{++}$, \ using again that \ $X_0$ \ is regularly varying with index \ $\alpha$, \ we have
 \begin{align*}
  N \EE(a_N^{-2} X_0^2 \bbone_{\{X_0\leq a_N\vare\}})
  = \frac{\EE\bigl(X_0^2 \bbone_{\{X_0\leq a_N\vare\}}\bigr)}
          {(a_N\vare)^2\PP(X_0>a_N\vare)}
     \frac{\PP(X_0>a_N\vare)}{\PP(X_0>a_N)}
     \vare^2 N \PP(X_0>a_N)
  \to \frac{\alpha}{2-\alpha} \vare^{2-\alpha}
 \end{align*}
 as \ $N \to \infty$, \ and, as \ $\vare \downarrow 0$, \ we conclude
 \eqref{Res_cond}.

Consequently, we may apply Theorem \ref{7.1}, and we obtain \eqref{help6_simple_aggregation1_stable_fd},
 where \ $(\bcX_t^{(k,\alpha)})_{t\in\RR_+}$ \ is an $\alpha$-stable process such that the characteristic function of the distribution \ $\mu_{k,\alpha}$ \ of \ $\bcX_1^{(k,\alpha)}$ \ has the form given in Theorem \ref{simple_aggregation1_stable_fd}.
Indeed, \eqref{fint} is valid for each Borel measurable function \ $f : \RR^{k+1}_0 \to \CC$ \ as well, for which the real and imaginary parts of the right hand side of \eqref{fint} are well defined.
Hence for all \ $\btheta \in \RR^{k+1}$, \ by \eqref{hmu},
 \begin{align*}
  &\widehat{\mu_{k,\alpha} }(\btheta)
   = \exp\biggl\{
          \int_{\RR_0^{k+1}} \biggl(\ee^{\ii \langle\btheta, \by\rangle} - 1 - \ii \sum_{\ell=1}^{k+1} \langle\be_\ell, \btheta\rangle \langle\be_\ell, \by\rangle \bbone_{(0,1]}(|\langle\be_\ell,\by\rangle|)\biggr) \nu_{k,\alpha}(\dd\by)
         \biggr\} \\
  &= \exp\biggl\{
             (1 - m_\xi^\alpha) \sum_{j=0}^k \int_0^\infty \biggl(\ee^{\ii\langle \btheta, \bv_j^{(k)}\rangle u} - 1 - \ii u \sum_{\ell=j+1}^{k+1} \langle\be_\ell, \btheta\rangle \langle\be_\ell, \bv_j^{(k)}\rangle \bbone_{(0,1]}(u\langle\be_\ell,\bv_j^{(k)}\rangle)\biggr) \alpha u^{-1-\alpha} \, \dd u
         \biggr\} ,
 \end{align*}
 since it will turn out that the real and imaginary parts of the exponent in the last expression are well defined.
If \ $\alpha \in (0, 1)$, \ then
 \[
   \int_0^\infty (\ee^{\pm\ii x} - 1) x^{-1-\alpha} \, \dd x
   = \Gamma(-\alpha) \ee^{\mp\ii\pi\alpha/2} ,
 \]
 see, e.g., (14.18) in Sato \cite{Sato} and its complex conjugate, thus for each
 \ $\vartheta \in \RR_{++}$,
 \begin{align*}
  &\int_0^\infty (\ee^{\ii\vartheta u} - 1) u^{-1-\alpha} \, \dd u
   = \vartheta^\alpha \int_0^\infty (\ee^{\ii v} - 1) v^{-1-\alpha} \, \dd v
   = \vartheta^\alpha \Gamma(-\alpha) \ee^{-\ii\pi\alpha/2}
    \\
  &= \vartheta^\alpha \frac{\Gamma(2-\alpha)}{(1-\alpha)(-\alpha)} \cos\left(\frac{\pi\alpha}{2}\right)
     \left(1 - \ii \tan\left(\frac{\pi\alpha}{2}\right)\right)
   = - \frac{C_\alpha}{\alpha} \vartheta^\alpha \left(1 - \ii \tan\left(\frac{\pi\alpha}{2}\right)\right) .
 \end{align*}
In a similar way, for each \ $\vartheta \in \RR_{--}$,
 \begin{align*}
  &\int_0^\infty (\ee^{\ii\vartheta u} - 1) u^{-1-\alpha} \, \dd u
   = (-\vartheta)^\alpha \int_0^\infty (\ee^{-\ii v} - 1) v^{-1-\alpha} \, \dd v
   = (-\vartheta)^\alpha \Gamma(-\alpha) \ee^{\ii\pi\alpha/2} \\
  &= (-\vartheta)^\alpha \frac{\Gamma(2-\alpha)}{(1-\alpha)(-\alpha)} \cos\left(\frac{\pi\alpha}{2}\right)
     \left(1 + \ii \tan\left(\frac{\pi\alpha}{2}\right)\right)
   = - \frac{C_\alpha}{\alpha} (-\vartheta)^\alpha \left(1 + \ii \tan\left(\frac{\pi\alpha}{2}\right)\right) .
 \end{align*}
Thus, for each \ $\vartheta \in \RR$,
 \begin{align}\label{help_Z1}
   \int_0^\infty
    (\ee^{\ii\vartheta u} - 1) u^{-1-\alpha}
    \, \dd u
   = - \frac{C_\alpha}{\alpha} |\vartheta|^\alpha
     \left(1 - \ii \tan\left(\frac{\pi\alpha}{2}\right) \sign(\vartheta)\right) ,
 \end{align}
 and hence, for each \ $\btheta \in \RR^{k+1}$ \ and \ $j \in \{0, \ldots, k\}$,
 \begin{align*}
  &\int_0^\infty \biggl(\ee^{\ii\langle \btheta, \bv_j^{(k)}\rangle u} - 1 - \ii u \sum_{\ell=j+1}^{k+1} \langle\be_\ell, \btheta\rangle \langle\be_\ell, \bv_j^{(k)}\rangle \bbone_{(0,1]}(u\langle\be_\ell,\bv_j^{(k)}\rangle)\biggr) \alpha u^{-1-\alpha} \, \dd u \\
  &= - C_\alpha |\langle \btheta, \bv_j^{(k)}\rangle|^\alpha
     \left(1 - \ii \tan\left(\frac{\pi\alpha}{2}\right) \sign(\langle \btheta, \bv_j^{(k)}\rangle)\right)
     - \ii \alpha \sum_{\ell=j+1}^{k+1} \langle\be_\ell, \btheta\rangle \langle\be_\ell, \bv_j^{(k)}\rangle \int_0^{1/\langle\be_\ell,\bv_j^{(k)}\rangle} u^{-\alpha} \, \dd u \\
  &= - C_\alpha |\langle \btheta, \bv_j^{(k)}\rangle|^\alpha
     \left(1 - \ii \tan\left(\frac{\pi\alpha}{2}\right) \sign(\langle \btheta, \bv_j^{(k)}\rangle)\right)
     - \ii \frac{\alpha}{1-\alpha} \sum_{\ell=j+1}^{k+1} \langle\be_\ell, \btheta\rangle \langle\be_\ell, \bv_j^{(k)}\rangle \biggl(\frac{1}{\langle\be_\ell,\bv_j^{(k)}\rangle}\biggr)^{1-\alpha} \\
  &= - C_\alpha |\langle \btheta, \bv_j^{(k)}\rangle|^\alpha
     \left(1 - \ii \tan\left(\frac{\pi\alpha}{2}\right) \sign(\langle \btheta, \bv_j^{(k)}\rangle)\right)
     - \ii \frac{\alpha}{1-\alpha} \sum_{\ell=j+1}^{k+1} \langle\be_\ell, \btheta\rangle \langle\be_\ell, \bv_j^{(k)}\rangle^\alpha .
 \end{align*}
Consequently,
 \begin{equation}\label{mu_k}
  \begin{aligned}
   \widehat{\mu_{k,\alpha}}(\btheta)
   &= \exp\Bigl\{- C_\alpha (1 - m_\xi^\alpha)
                   \sum_{j=0}^k
                    |\langle\btheta,\bv_j^{(k)}\rangle|^\alpha
                     \left(1 - \ii \tan\left(\frac{\pi\alpha}{2}\right)
                                   \sign(\langle\btheta,\bv_j^{(k)}\rangle)\right) \\
   &\phantom{= \exp\Bigl\{}
                  - \ii \frac{\alpha}{1-\alpha} (1 - m_\xi^\alpha) \sum_{j=0}^k \sum_{\ell=j+1}^{k+1} \langle\be_\ell, \btheta\rangle
                   \langle\be_\ell, \bv_j^{(k)}\rangle^\alpha\Bigr\}
  \end{aligned}
 \end{equation}
 for all \ $\btheta \in \RR^{k+1}$, \ where
 \begin{align*}
  \sum_{j=0}^k \sum_{\ell=j+1}^{k+1}
   \langle\be_\ell, \btheta\rangle \langle\be_\ell, \bv_j^{(k)}\rangle^\alpha
  = \sum_{\ell=1}^{k+1}
     \langle\be_\ell, \btheta\rangle
     \sum_{j=0}^{\ell-1}
      \langle\be_\ell, \bv_j^{(k)}\rangle^\alpha
  = \frac{1}{1-m_\xi^\alpha}
    \sum_{\ell=1}^{k+1}
     \langle\be_\ell, \btheta\rangle ,
 \end{align*}
 since \ $\langle\be_1, \bv_0^{(k)}\rangle^\alpha = (1-m_\xi^\alpha)^{-1}$, \ and, for each \ $\ell \in \{2, \ldots, k + 1\}$, \ we have
 \[
   \sum_{j=0}^{\ell-1} \langle\be_\ell, \bv_j^{(k)}\rangle^\alpha
   = \frac{m_\xi^{(\ell-1)\alpha}}{1-m_\xi^\alpha}
     + \sum_{j=1}^{\ell-1} m_\xi^{(\ell-j-1)\alpha}
   = \frac{m_\xi^{(\ell-1)\alpha}}{1-m_\xi^\alpha}
     + \frac{1-m_\xi^{(\ell-1)\alpha}}{1-m_\xi^\alpha}
   = \frac{1}{1-m_\xi^\alpha} .
 \]
Hence we obtain
 \begin{equation}\label{drift_center}
   (1 - m_\xi^\alpha)
   \sum_{j=0}^k \sum_{\ell=j+1}^{k+1}
    \langle\be_\ell, \btheta\rangle \langle\be_\ell, \bv_j^{(k)}\rangle^\alpha
   = \langle\btheta, \bone_{k+1}\rangle ,
 \end{equation}
 yielding the statement in case of \ $\alpha \in (0, 1)$.
\ Note that the above calculation shows that \eqref{drift_center} is valid for each \ $\alpha\in(0,2)$.

If \ $\alpha \in (1, 2)$, \ then
 \[
   \int_0^\infty (\ee^{\pm\ii x} - 1 \mp \ii x) x^{-1-\alpha} \, \dd x
   = \Gamma(-\alpha) \, \ee^{\mp\ii\pi\alpha/2} ,
 \]
 see, e.g., (14.19) in Sato \cite{Sato} and its complex conjugate, thus for each
 \ $\vartheta \in \RR_{++}$,
 \begin{align*}
  &\int_0^\infty
    (\ee^{\ii\vartheta u} - 1 - \ii \vartheta u) u^{-1-\alpha} \, \dd u
   = \vartheta^\alpha \int_0^\infty (\ee^{\ii x} - 1 - \ii x) x^{-1-\alpha} \, \dd x
   = \vartheta^\alpha \Gamma(-\alpha) \ee^{-\ii\pi\alpha/2} \\
  &= \vartheta^\alpha \frac{\Gamma(2-\alpha)}{(1-\alpha)(-\alpha)} \cos\left(\frac{\pi\alpha}{2}\right)
     \left(1 - \ii \tan\left(\frac{\pi\alpha}{2}\right)\right)
   = - \frac{C_\alpha}{\alpha} \vartheta^\alpha \left(1 - \ii \tan\left(\frac{\pi\alpha}{2}\right)\right) .
 \end{align*}
In a similar way, for each \ $\vartheta \in \RR_{--}$,
 \begin{align*}
  &\int_0^\infty
      (\ee^{\ii\vartheta u} - 1 - \ii \vartheta u) u^{-1-\alpha} \, \dd u
   = (-\vartheta)^\alpha
     \int_0^\infty (\ee^{-\ii x} - 1 + \ii x) x^{-1-\alpha} \, \dd x
   = (-\vartheta)^\alpha \Gamma(-\alpha) \ee^{\ii\pi\alpha/2} \\
  &= (-\vartheta)^\alpha \frac{\Gamma(2-\alpha)}{(1-\alpha)(-\alpha)} \cos\left(\frac{\pi\alpha}{2}\right)
     \left(1 + \ii \tan\left(\frac{\pi\alpha}{2}\right)\right)
   = - \frac{C_\alpha}{\alpha} (-\vartheta)^\alpha \left(1 + \ii \tan\left(\frac{\pi\alpha}{2}\right)\right) .
 \end{align*}
Thus, for each \ $\vartheta \in \RR$,
 \begin{align}\label{help_Z2}
   \int_0^\infty
    (\ee^{\ii\vartheta u} - 1 - \ii \vartheta u) u^{-1-\alpha}
    \, \dd u
   = - \frac{C_\alpha}{\alpha} |\vartheta|^\alpha
     \left(1 - \ii \tan\left(\frac{\pi\alpha}{2}\right) \sign(\vartheta)\right) ,
 \end{align}
 and hence, for each \ $\btheta \in \RR^{k+1}$ \ and \ $j \in \{0, \ldots, k\}$,
 \begin{align*}
  &\int_0^\infty \biggl(\ee^{\ii\langle \btheta, \bv_j^{(k)}\rangle u} - 1 - \ii u \sum_{\ell=j+1}^{k+1} \langle\be_\ell, \btheta\rangle \langle\be_\ell, \bv_j^{(k)}\rangle \bbone_{(0,1]}(u\langle\be_\ell,\bv_j^{(k)}\rangle)\biggr) \alpha u^{-1-\alpha} \, \dd u \\
  &= - C_\alpha |\langle \btheta, \bv_j^{(k)}\rangle|^\alpha
     \left(1 - \ii \tan\left(\frac{\pi\alpha}{2}\right) \sign(\langle \btheta, \bv_j^{(k)}\rangle)\right)
     + \ii \alpha \sum_{\ell=j+1}^{k+1} \langle\be_\ell, \btheta\rangle \langle\be_\ell, \bv_j^{(k)}\rangle \int_{1/\langle\be_\ell,\bv_j^{(k)}\rangle}^\infty u^{-\alpha} \, \dd u \\
  &= - C_\alpha |\langle \btheta, \bv_j^{(k)}\rangle|^\alpha
     \left(1 - \ii \tan\left(\frac{\pi\alpha}{2}\right) \sign(\langle \btheta, \bv_j^{(k)}\rangle)\right)
     - \ii \frac{\alpha}{1-\alpha} \sum_{\ell=j+1}^{k+1} \langle\be_\ell, \btheta\rangle \langle\be_\ell, \bv_j^{(k)}\rangle \biggl(\frac{1}{\langle\be_\ell,\bv_j^{(k)}\rangle}\biggr)^{1-\alpha} \\
  &= - C_\alpha |\langle \btheta, \bv_j^{(k)}\rangle|^\alpha
     \left(1 - \ii \tan\left(\frac{\pi\alpha}{2}\right) \sign(\langle \btheta, \bv_j^{(k)}\rangle)\right)
     + \ii \frac{\alpha}{\alpha-1} \sum_{\ell=j+1}^{k+1} \langle\be_\ell, \btheta\rangle \langle\be_\ell, \bv_j^{(k)}\rangle^\alpha .
 \end{align*}
Consequently, we obtain \eqref{mu_k} for all \ $\btheta \in \RR^{k+1}$, \ and, applying again \eqref{drift_center}, we conclude the statement in case of \ $\alpha \in (1, 2)$.

Finally, we consider the case \ $\alpha = 1$.
\ For each \ $\vartheta \in \RR_{++}$,
 \[
   \int_0^\infty
    (\ee^{\ii\vartheta u} - 1 - \ii \vartheta u \bbone_{(0,1]}(u)) u^{-2} \, \dd u
   = - \frac{\pi\vartheta}{2} - \ii \vartheta \log(\vartheta) + \ii C \vartheta ,
 \]
 where \ $C$ \ is given in \eqref{c}, see, e.g., (14.20) in Sato \cite{Sato}.
Its complex conjugate has the form
 \[
   \int_0^\infty
    (\ee^{-\ii\vartheta u} - 1 + \ii \vartheta u \bbone_{(0,1]}(u)) u^{-2} \, \dd u
   = - \frac{\pi\vartheta}{2} + \ii \vartheta \log(\vartheta) - \ii C \vartheta , \qquad
   \vartheta \in \RR_{++} ,
 \]
 thus
 \[
   \int_0^\infty
    (\ee^{\ii(-\vartheta)u} - 1 - \ii (-\vartheta) u \bbone_{(0,1]}(u)) u^{-2} \, \dd u
   = - \frac{\pi(-(-\vartheta))}{2} - \ii (-\vartheta) \log(-(-\vartheta)) + \ii C (-\vartheta)
 \]
 for \ $\vartheta \in \RR_{++}$, \ and hence
 \begin{align*}
  &\int_0^\infty
   (\ee^{\ii\vartheta u} - 1 - \ii \vartheta u \bbone_{(0,1]}(u)) u^{-2} \, \dd u
   = - \frac{\pi|\vartheta|}{2} - \ii \vartheta \log(|\vartheta|) + \ii C \vartheta \\
  &= - C_1 |\vartheta| \left(1 + \ii \frac{2}{\pi} \sign(\vartheta) \log(|\vartheta|)\right) + \ii C \vartheta , \qquad \vartheta \in \RR .
 \end{align*}
Consequently, for each \ $\btheta \in \RR^{k+1}$ \ and \ $j \in \{0, \ldots, k\}$,
 \begin{align*}
  &\int_0^\infty \biggl(\ee^{\ii\langle \btheta, \bv_j^{(k)}\rangle u} - 1 - \ii u \sum_{\ell=j+1}^{k+1} \langle\be_\ell, \btheta\rangle \langle\be_\ell, \bv_j^{(k)}\rangle \bbone_{(0,1]}(u\langle\be_\ell,\bv_j^{(k)}\rangle)\biggr) u^{-2} \, \dd u \\
  &= - C_1 |\langle \btheta, \bv_j^{(k)}\rangle| \left(1 + \ii \frac{2}{\pi} \sign(\langle \btheta, \bv_j^{(k)}\rangle) \log(|\langle \btheta, \bv_j^{(k)}\rangle|)\right) + \ii C \langle \btheta, \bv_j^{(k)}\rangle \\
  &\quad
  + \ii \sum_{\ell=j+1}^{k+1} \langle\be_\ell, \btheta\rangle \langle\be_\ell, \bv_j^{(k)}\rangle \int_{1/\langle\be_\ell,\bv_j^{(k)}\rangle}^1 u^{-1} \, \dd u \\
  &= - C_1 |\langle \btheta, \bv_j^{(k)}\rangle| \left(1 + \ii \frac{2}{\pi} \sign(\langle \btheta, \bv_j^{(k)}\rangle) \log(|\langle \btheta, \bv_j^{(k)}\rangle|)\right) + \ii C \langle \btheta, \bv_j^{(k)}\rangle \\
  &\quad
  - \ii \sum_{\ell=j+1}^{k+1} \langle\be_\ell, \btheta\rangle \langle\be_\ell, \bv_j^{(k)}\rangle \log\biggl(\frac{1}{\langle\be_\ell,\bv_j^{(k)}\rangle}\biggr) \\
  &= - C_1 |\langle \btheta, \bv_j^{(k)}\rangle| \left(1 + \ii \frac{2}{\pi} \sign(\langle \btheta, \bv_j^{(k)}\rangle) \log(|\langle \btheta, \bv_j^{(k)}\rangle|)\right) + \ii C \langle \btheta, \bv_j^{(k)}\rangle \\
  &\quad
     + \ii \sum_{\ell=j+1}^{k+1} \langle\be_\ell, \btheta\rangle \langle\be_\ell, \bv_j^{(k)}\rangle \log(\langle\be_\ell,\bv_j^{(k)}\rangle) .
 \end{align*}
Applying again \eqref{drift_center}, we have the statement in case of \ $\alpha = 1$.
\proofend

\noindent{\bf Proof of Corollary \ref{simple_aggregation1_stable_centering_fd}.}
In case of \ $\alpha \in (0, 1)$, \ by \eqref{centering01T} with \ $t = 1$, \ we have
 \begin{equation}\label{centering01}
  \lim_{N\to\infty} \frac{N}{a_N} \EE\bigl(X_0 \bbone_{\{X_0\leq a_N\}}\bigr)
  = \frac{\alpha}{1-\alpha} .
 \end{equation}
Next, we may apply Lemma \ref{Conv2Funct} with
 \begin{gather*}
  \bcU_t^{(N)} := \frac{1}{a_N}
     \sum_{j=1}^{\lfloor Nt\rfloor}
      \bigl(X^{(j)}_0, \ldots, X^{(j)}_k\bigr)^\top
     - \frac{\lfloor Nt\rfloor}{a_N} \EE\bigl(X_0 \bbone_{\{X_0\leq a_N\}}\bigr) \bone_{k+1} ,
      \qquad N \in \NN , \\
  \Phi_N(f)(t) := f(t) + \frac{\lfloor Nt\rfloor}{a_N} \EE(X_0  \bbone_{\{X_0\leq a_N\}} ) \bone_{k+1},
                  \qquad N \in \NN , \\
  \bcU_t := \bcX_t^{(k,\alpha)} , \qquad
  \Phi(f)(t) := f(t) + \frac{\alpha}{1-\alpha} t \bone_{k+1}
 \end{gather*}
 for \ $t \in \RR_+$ \ and \ $f \in \DD(\RR_+,\RR^{k+1})$.
\ Indeed, in order to show \ $\Phi_N(f_N) \to \Phi(f)$ \ in \ $\DD(\RR_+, \RR^{k+1})$ \ as \ $N \to \infty$ \
 whenever \ $f_N \to f$ \ in \ $\DD(\RR_+, \RR^{k+1})$ \ as \ $N \to \infty$ \ with \ $f, f_N \in \DD(\RR_+, \RR^{k+1})$, \ $N \in \NN$,
 \ by Propositions VI.1.17 and VI.1.23 in Jacod and Shiryaev \cite{JacShi}, it is enough to check that
 for each \ $T \in \RR_{++}$, \ we have
 \begin{align*}
  \sup_{t\in[0,T]} \biggl\|\frac{\lfloor Nt\rfloor}{a_N} \EE\bigl(X_0 \bbone_{\{X_0\leq a_N\}}\bigr) \bone_{k+1} - \frac{\alpha}{1-\alpha} t \bone_{k+1}\biggr\|
  \to 0 \qquad \text{as \ $N\to\infty$.}
 \end{align*}
This follows, since, by \eqref{centering01}, we obtain
 \begin{align*}
  &\sup_{t\in[0,T]} \biggl\|\frac{\lfloor Nt\rfloor}{a_N} \EE\bigl(X_0 \bbone_{\{X_0\leq a_N\}}\bigr) \bone_{k+1} - \frac{\alpha}{1-\alpha} t \bone_{k+1}\biggr\| \\
  &\leq \sup_{t\in[0,T]} \biggl\|\frac{\lfloor Nt\rfloor}{N} \biggl(\frac{N}{a_N} \EE\bigl(X_0 \bbone_{\{X_0\leq a_N\}}\bigr) - \frac{\alpha}{1-\alpha}\biggr) \bone_{k+1}\biggr\|
        + \sup_{t\in[0,T]} \Biggl\|\frac{\alpha}{1-\alpha} \biggl(\frac{\lfloor Nt\rfloor}{N} - t\biggr) \bone_{k+1}\Biggr\| \\
  &\leq T \sqrt{k + 1} \biggl|\frac{N}{a_N} \EE\bigl(X_0 \bbone_{\{X_0\leq a_N\}}\bigr) - \frac{\alpha}{1-\alpha}\biggr|
        + \frac{\sqrt{k+1}}{N} \frac{\alpha}{1-\alpha}
   \to 0 \qquad \text{as \ $N \to \infty$.}
 \end{align*}
Applying Lemma \ref{Conv2Funct}, we obtain
 \[
   \biggl(\frac{1}{a_N} \sum_{j=1}^{\lfloor Nt\rfloor} (X_0^{(j)}, \ldots, X_k^{(j)})\biggr)_{t\in\RR_+}
   = \Phi_N(\bcU^{(N)})
   \distr \Phi(\bcU) \qquad \text{as \ $N \to \infty$,}
 \]
 where \ $\Phi(\bcU)_t = \bcX_t^{(k,\alpha)} + \frac{\alpha}{1-\alpha} t \bone_{k+1}$,
 \ $t \in \RR_+$, \ is a $(k + 1)$-dimensional \ $\alpha$-stable process.
By Theorem \ref{simple_aggregation1_stable_fd} and Remark \ref{mu_properties2},
 the characteristic function of \ $\bcX_1^{(k,\alpha)} + \frac{\alpha}{1-\alpha} \bone_{k+1}$
 \ has the form given in the theorem, and hence we conclude the statement in case of \ $\alpha \in (0, 1)$.

In case of \ $\alpha \in (1, 2)$, \ by \eqref{centering12T} with \ $t = 1$, \ we have
 \begin{equation}\label{centering12}
  \lim_{N\to\infty} \frac{N}{a_N} \EE\bigl(X_0 \bbone_{\{X_0>a_N\}}\bigr)
  = \frac{\alpha}{\alpha-1} .
 \end{equation}
Next, we may apply Lemma \ref{Conv2Funct} with \ $\bcU$, \ $\Phi$ \ and \ $\bcU^{(N)}$, \ $N \in \NN$, \ as defined above, and with
 \[
   \Phi_N(f)(t) := f(t) - \frac{\lfloor Nt\rfloor}{a_N} \EE(X_0  \bbone_{\{X_0>a_N\}} ) \bone_{k+1},
                  \qquad N \in \NN , \quad f \in \DD(\RR_+,\RR^{k+1}) , \quad t \in \RR_+ .
 \]
Indeed, in order to show \ $\Phi_N(f_N) \to \Phi(f)$ \ in \ $\DD(\RR_+, \RR^{k+1})$ \ as \ $N \to \infty$ \
 whenever \ $f_N \to f$ \ in \ $\DD(\RR_+, \RR^{k+1})$ \ as \ $N \to \infty$ \ with \ $f, f_N \in \DD(\RR_+, \RR^{k+1})$, \ $N \in \NN$,
 \ by Propositions VI.1.17 and VI.1.23 in Jacod and Shiryaev \cite{JacShi}, it is enough to check that
 for each \ $T \in \RR_{++}$, \ we have
 \begin{align*}
  \sup_{t\in[0,T]} \biggl\|\frac{\lfloor Nt\rfloor}{a_N} \EE\bigl(X_0 \bbone_{\{X_0 > a_N\}}\bigr) \bone_{k+1} - \frac{\alpha}{\alpha-1} t \bone_{k+1}\biggr\| \to 0
     \qquad \text{as \ $N\to\infty$.}
 \end{align*}
This follows, since, by \eqref{centering12}, we obtain
 \begin{align*}
  &\sup_{t\in[0,T]} \biggl\|\frac{\lfloor Nt\rfloor}{a_N} \EE\bigl(X_0 \bbone_{\{X_0>a_N\}}\bigr) \bone_{k+1} - \frac{\alpha}{\alpha-1} t \bone_{k+1}\biggr\| \\
  &\leq \sup_{t\in[0,T]} \biggl\|\frac{\lfloor Nt\rfloor}{N} \biggl(\frac{N}{a_N} \EE\bigl(X_0 \bbone_{\{X_0>a_N\}}\bigr) - \frac{\alpha}{\alpha-1}\biggr) \bone_{k+1}\biggr\|
        + \sup_{t\in[0,T]} \Biggl\|\frac{\alpha}{\alpha-1} \biggl(\frac{\lfloor Nt\rfloor}{N} - t\biggr) \bone_{k+1}\Biggr\| \\
  &\leq T \sqrt{k + 1} \biggl|\frac{N}{a_N} \EE\bigl(X_0 \bbone_{\{X_0>a_N\}}\bigr) - \frac{\alpha}{\alpha-1}\biggr|
        + \frac{\sqrt{k+1}}{N} \frac{\alpha}{\alpha-1}
   \to 0 \qquad \text{as \ $N \to \infty$.}
 \end{align*}
Applying Lemma \ref{Conv2Funct}, we obtain
 \[
   \biggl(\frac{1}{a_N} \sum_{j=0}^{\lfloor Nt\rfloor} \bigl(X_1^{(j)}, \ldots, X_k^{(j)}\bigr) - \frac{\lfloor Nt\rfloor}{a_N} \EE(X_0) \bone_{k+1}\biggr)_{t\in\RR_+}
   = \Phi_N(\bcU^{(N)})
   \distr \Phi(\bcU) \qquad \text{as \ $N \to \infty$,}
 \]
 where \ $\Phi(\bcU)_t = \bcX_t^{(k,\alpha)} - t \frac{\alpha}{\alpha-1} \bone_{k+1}$, \ $t \in \RR_+$, \ is a $(k + 1)$-dimensional \ $\alpha$-stable process.
By Theorem \ref{simple_aggregation1_stable_fd} and Remark \ref{mu_properties2}, the characteristic function of
 \ $\bcX_1^{(k,\alpha)} - \frac{\alpha}{\alpha-1} \bone_{k+1}$ \ has the form given in the theorem,
 and hence we conclude the statement in case of \ $\alpha \in (1, 2)$ \ as well.
\proofend

\noindent{\bf Proof of Proposition \ref{Pro_AR1}.}
The sequence \ $(\tvare^{(\alpha)}_k)_{k\in\NN}$ \ consists of identically distributed random variables,
 since the strong stationarity of \ $\bigl(\cY^{(\alpha)}_k\bigr)_{k\in\ZZ_+}$ \ yields that
 \ $\bigl(\cY^{(\alpha)}_k,\cY^{(\alpha)}_{k-1} \bigr)_{k\in\NN}$ \ consists of identically distributed random
 variables.

In what follows, let \ $k\in\NN$ \ be fixed.
By \eqref{calY_def}, the characteristic function of \ $(\cY^{(\alpha)}_0, \ldots, \cY^{(\alpha)}_{k-1}, \tvare_k)$ \ has the form
 \[
   \EE\bigl(\ee^{\ii(\vartheta_0\cY^{(\alpha)}_0+\cdots+\vartheta_{k-1}\cY^{(\alpha)}_{k-1}+\vartheta_k\tvare^{(\alpha)}_k)}\bigr)
   = \EE\bigl(\ee^{\ii(\vartheta_0\cY^{(\alpha)}_0+\cdots+(\vartheta_{k-1}-m_\xi\vartheta_k)\cY^{(\alpha)}_{k-1}+\vartheta_k\cY^{(\alpha)}_k)}\bigr)
   = \EE\bigl(\ee^{\ii\langle\btheta_k,\bcX^{(k,\alpha)}_1\rangle}\bigr)
 \]
 for \ $(\vartheta_0, \ldots, \vartheta_k)^\top \in \RR^{k+1}$, \ where
 \[
   \btheta_k := (\vartheta_0, \ldots, \vartheta_{k-2}, \vartheta_{k-1} - m_\xi \vartheta_k, \vartheta_k)^\top \in \RR^{k+1} .
 \]
We can write \ $\btheta_k = \btheta_k^{(1)} + \btheta_k^{(2)}$ \ with
 \[
   \btheta_k^{(1)} := (\vartheta_0, \ldots, \vartheta_{k-2}, \vartheta_{k-1}, 0)^\top , \qquad
   \btheta_k^{(2)} := (0, \ldots, 0, - m_\xi \vartheta_k, \vartheta_k)^\top .
 \]
We have \ $\langle\btheta_k^{(2)}, \bv_j^{(k)}\rangle = 0$ \ for each \ $j \in \{0, \ldots, k - 1\}$, \ and \ $\langle\btheta_k^{(1)}, \bv_k^{(k)}\rangle = 0$, \ hence
 \[
   \langle\btheta_k,\bv_j^{(k)}\rangle
   = \begin{cases}
      \langle\btheta_k^{(1)},\bv_j^{(k)}\rangle & \text{if \ $j \in \{0, \ldots, k - 1\}$,} \\
      \langle\btheta_k^{(2)},\bv_k^{(k)}\rangle & \text{if \ $j = k$.}
     \end{cases}
 \]
In case of \ $\alpha \in (0, 1) \cup (1, 2)$, \ by Theorem \ref{simple_aggregation1_stable_fd}, we have
 \begin{align*}
  &\EE\bigl(\ee^{\ii(\vartheta_0\cY^{(\alpha)}_0+\cdots+\vartheta_{k-1}\cY^{(\alpha)}_{k-1}+\vartheta_k\tvare^{(\alpha)}_k)}\bigr) \\
  &= \exp\biggl\{- C_\alpha (1 - m_\xi^\alpha)
                   \sum_{j=0}^k
                    |\langle\btheta_k,\bv_j^{(k)}\rangle|^\alpha
                     \left(1 - \ii \tan\left(\frac{\pi\alpha}{2}\right)
                                   \sign(\langle\btheta_k,\bv_j^{(k)}\rangle)\right)
                 - \ii \frac{\alpha}{1-\alpha} \langle\btheta_k, \bone_{k+1}\rangle\biggr\} \\
  &= \exp\biggl\{- C_\alpha (1 - m_\xi^\alpha)
                   \sum_{j=0}^{k-1}
                    |\langle\btheta_k^{(1)},\bv_j^{(k)}\rangle|^\alpha
                     \left(1 - \ii \tan\left(\frac{\pi\alpha}{2}\right)
                                   \sign(\langle\btheta_k^{(1)},\bv_j^{(k)}\rangle)\right)
                 - \ii \frac{\alpha}{1-\alpha} \langle\btheta_k^{(1)}, \bone_{k+1}\rangle \\
  &\phantom{= \exp\biggl\{}
                 - C_\alpha (1 - m_\xi^\alpha)
                   |\langle\btheta_k^{(2)},\bv_k^{(k)}\rangle|^\alpha
                    \left(1 - \ii \tan\left(\frac{\pi\alpha}{2}\right)
                                  \sign(\langle\btheta_k^{(2)},\bv_k^{(k)}\rangle)\right)
                 - \ii \frac{\alpha}{1-\alpha} \langle\btheta_k^{(2)}, \bone_{k+1}\rangle \biggr\} \\
  &= \EE\bigl(\ee^{\ii(\vartheta_0\cY^{(\alpha)}_0+\cdots+\vartheta_{k-1}\cY^{(\alpha)}_{k-1})}\bigr)
  \EE\bigl(\ee^{\ii\vartheta_k\tvare^{(\alpha)}_k}\bigr) ,
 \end{align*}
where we used that
 \begin{align*}
   &\langle\btheta_k^{(1)},\bv_j^{(k)}\rangle = \langle (\vartheta_0, \ldots, \vartheta_{k-1})^\top , \bv_j^{(k-1)} \rangle, \qquad j=0,1,\ldots,k-1,\\
   & \langle\btheta_k^{(1)}, \bone_{k+1}\rangle  = \langle (\vartheta_0, \ldots, \vartheta_{k-1})^\top , \bone_k \rangle,\\
   & \langle\btheta_k^{(2)},\bv_k^{(k)}\rangle = \vartheta_k,\qquad  \langle\btheta_k^{(2)}, \bone_{k+1}\rangle = \vartheta_k,
 \end{align*}
 and \ $\vartheta_0 = \vartheta_1=\ldots = \vartheta_{k-1}=0$ \ yields \ $\btheta_k^{(1)} = \bzero\in\RR^{k+1}$.
\ Thus we obtain the independence of \ $(\cY^{(\alpha)}_0, \ldots, \cY^{(\alpha)}_{k-1})^\top$ \ and \ $\tvare^{(\alpha)}_k$, \ and the characteristic function of \ $\tvare^{(\alpha)}_k$ \ has the form
 \[
   \EE\bigl(\ee^{\ii\vartheta_k\tvare^{(\alpha)}_k}\bigr)
   = \exp\biggl\{- C_\alpha (1 - m_\xi^\alpha) |\vartheta_k|^\alpha
                   \left(1 - \ii \tan\left(\frac{\pi\alpha}{2}\right)
                                  \sign(\vartheta_k)\right)
                 - \ii \frac{\alpha}{1-\alpha} \vartheta_k\biggr\} , \qquad \vartheta_k \in \RR ,
 \]
 hence \ $\tvare^{(\alpha)}_k$ \ is \ $\alpha$-stable (see, e.g., Sato \cite[Theorem 14.10]{Sato}).
In fact, \ $\tvare^{(\alpha)}_k + \frac{\alpha}{1-\alpha} \distre (1 - m_\xi^\alpha)^{\frac{1}{\alpha}} \bigl(\cY^{(\alpha)}_0 + \frac{\alpha}{1-\alpha}\bigr)$, \ $k \in \NN$, \ since for all \ $\vartheta_k \in \RR$, \ we have
 \begin{align*}
  &\EE\Bigl(\ee^{\ii\vartheta_k(1-m_\xi^\alpha)^{\frac{1}{\alpha}}\bigl(\cY^{(\alpha)}_0+\frac{\alpha}{1-\alpha}\bigr)}\Bigr)
   = \ee^{\ii \frac{\alpha}{1-\alpha}\vartheta_k (1 - m_\xi^\alpha)^{\frac{1}{\alpha}} }
    \widehat\mu_{0,\alpha} (\vartheta_k (1 - m_\xi^\alpha)^{\frac{1}{\alpha}}) \\
  &= \exp\Bigl\{- C_\alpha (1 - m_\xi^\alpha) |\langle\vartheta_k (1 - m_\xi^\alpha)^{\frac{1}{\alpha}}, (1 - m_\xi^\alpha)^{-\frac{1}{\alpha}}\rangle|^\alpha \\
  &\phantom{= \exp\Bigg\{-}
                  \times \left(1 - \ii \tan\left(\frac{\pi\alpha}{2}\right) \sign(\langle \vartheta_k (1 - m_\xi^\alpha)^{\frac{1}{\alpha}},
                  (1 - m_\xi^\alpha)^{- \frac{1}{\alpha}}\rangle)\right)\Bigr\} \\
  &= \exp\Big\{- C_\alpha (1 - m_\xi^\alpha) |\vartheta_k|^\alpha \left(1 - \ii \tan\left(\frac{\pi\alpha}{2}\right) \sign(\vartheta_k)\right)\Big\} .
 \end{align*}
Further, \ $\bigl(\cY^{(\alpha)}_k + \frac{\alpha}{1-\alpha}\bigr)_{k\in\ZZ_+}$ \ is also a strongly stationary \ $\alpha$-stable time homogeneous Markov process.
In fact, it is a subcritical autoregressive process of order 1 with autoregressive coefficient \ $m_\xi$ \ such that the distribution of its
 innovations satisfies
 \begin{align*}
  \biggl(\cY^{(\alpha)}_k + \frac{\alpha}{1-\alpha}\biggr) - m_\xi \biggl(\cY^{(\alpha)}_{k-1} + \frac{\alpha}{1-\alpha}\biggr)
  &= \widetilde\vare_k + (1-m_\xi)\frac{\alpha}{1-\alpha} \\
  &\distre (1 - m_\xi^\alpha)^{\frac{1}{\alpha}} \biggl(\cY^{(\alpha)}_0 + \frac{\alpha}{1-\alpha}\biggr) -m_\xi\frac{\alpha}{1-\alpha} ,
  \qquad k \in \NN .
\end{align*}

In case of \ $\alpha = 1$, \ by Theorem \ref{simple_aggregation1_stable_fd}, we have
 \begin{align*}
  &\EE\bigl(\ee^{\ii(\vartheta_0\cY^{(1)}_0+\cdots+\vartheta_{k-1}\cY^{(1)}_{k-1}+\vartheta_k\tvare^{(1)}_k)}\bigr) \\
  &= \exp\biggl\{- C_1 (1 - m_\xi)
                   \sum_{j=0}^k
                    |\langle\btheta_k,\bv_j^{(k)}\rangle|
                    \Bigl(1 + \ii \frac{2}{\pi} \sign(\langle\btheta_k,\bv_j^{(k)}\rangle)
                                   \log(|\langle\btheta_k,\bv_j^{(k)}\rangle|)\Bigr) \\
  &\phantom{= \exp\biggl\{}
                 + \ii C \langle\btheta_k, \bone_{k+1}\rangle
                 + \ii (1 - m_\xi) \sum_{j=0}^k \sum_{\ell=j+1}^{k+1} \langle\be_\ell, \btheta_k\rangle
                   \langle\be_\ell, \bv_j^{(k)}\rangle \log(\langle\be_\ell, \bv_j^{(k)}\rangle)\biggr\} .
 \end{align*}

If \ $\alpha = 1$ \ and \ $m_\xi = 0$, \ then for each \ $j \in \{0, \ldots, k\}$, \ we have \ $\bv_j^{(k)} = \be_{j+1}$, \ and hence \ $\langle\btheta_k,\bv_j^{(k)}\rangle = \vartheta_j$, \ and \ $\langle\be_\ell, \bv_j^{(k)}\rangle = 1$ \ if \ $\ell = j + 1$ \ and \ $\langle\be_\ell, \bv_j^{(k)}\rangle = 0$ \ if \ $\ell \ne j + 1$.
\ Consequently,
 \begin{align*}
  &\EE\bigl(\ee^{\ii(\vartheta_0\cY^{(1)}_0+\cdots+\vartheta_{k-1}\cY^{(1)}_{k-1}+\vartheta_k\tvare^{(1)}_k)}\bigr) \\
  &= \exp\biggl\{- C_1
                   \sum_{j=0}^k
                    |\vartheta_j|
                    \Bigl(1 + \ii \frac{2}{\pi} \sign(\vartheta_j) \log(|\vartheta_j|)\Bigr)
                 + \ii C \sum_{j=0}^k \vartheta_j\biggr\} \\
  &= \EE\bigl(\ee^{\ii\vartheta_0\cY^{(\alpha)}_0}\bigr) \cdots \EE\bigl(\ee^{\ii\vartheta_0\cY^{(\alpha)}_{k-1}}\bigr)
  \EE\bigl(\ee^{\ii\vartheta_k\tvare^{(\alpha)}_k}\bigr) ,
 \end{align*}
 thus we obtain the independence of \ $\cY^{(1)}_0$, \ldots, $\cY^{(1)}_{k-1}$ \ and \ $\tvare^{(1)}_k$, \ and the characteristic function of \ $\tvare^{(1)}_k$ \ has the form
 \[
   \EE\bigl(\ee^{\ii\vartheta_k\tvare^{(1)}_k}\bigr)
   = \exp\Bigl\{- C_1 |\vartheta_k|
                   \Bigl(1 + \ii \frac{2}{\pi} \sign(\vartheta_k) \log(|\vartheta_k|)\Bigr)
                 + \ii C \vartheta_k\Bigr\} , \qquad \vartheta_k \in \RR ,
 \]
 hence \ $\tvare^{(1)}_k$ \ is \ 1-stable (see, e.g., Sato \cite[Theorem 14.10]{Sato}).
In fact, \ $\tvare^{(1)}_k \distre \cY^{(1)}_0$, \ and \ $(\cY_k^{(1)})_{k\in\ZZ_+}$ \ is a sequence of independent, identically distributed \ $1$-stable random variables.

If \ $\alpha = 1$ \ and \ $m_\xi \in (0, 1)$, \ then, by Theorem \ref{simple_aggregation1_stable_fd},
 \begin{align*}
  &\EE\bigl(\ee^{\ii(\vartheta_0\cY^{(1)}_0+\cdots+\vartheta_{k-1}\cY^{(1)}_{k-1}+\vartheta_k\tvare^{(1)}_k)}\bigr) \\
  &= \exp\biggl\{- C_1 (1 - m_\xi)
                   \sum_{j=0}^{k-1}
                    |\langle\btheta_k^{(1)},\bv_j^{(k)}\rangle|
                    \Bigl(1 + \ii \frac{2}{\pi} \sign(\langle\btheta_k^{(1)},\bv_j^{(k)}\rangle)
                                   \log(|\langle\btheta_k^{(1)},\bv_j^{(k)}\rangle|)\Bigr) \\
  &\phantom{= \exp\biggl\{}
                 + \ii C \langle\btheta_k^{(1)}, \bone_{k+1}\rangle
                 + \ii (1 - m_\xi) \sum_{j=0}^{k-1} \sum_{\ell=j+1}^k \langle\be_\ell, \btheta_k^{(1)}\rangle
                   \langle\be_\ell, \bv_j^{(k)}\rangle \log(\langle\be_\ell, \bv_j^{(k)}\rangle) \\
  &\phantom{= \exp\biggl\{}
                 - C_1 (1 - m_\xi)
                   |\langle\btheta_k^{(2)},\bv_k^{(k)}\rangle|
                   \Bigl(1 + \ii \frac{2}{\pi} \sign(\langle\btheta_k^{(2)},\bv_k^{(k)}\rangle)
                                   \log(|\langle\btheta_k^{(2)},\bv_k^{(k)}\rangle|)\Bigr) \\
  &\phantom{= \exp\biggl\{}
                 + \ii C \langle\btheta_k^{(2)}, \bone_{k+1}\rangle
                 + \ii (1 - m_\xi) \sum_{j=0}^k \sum_{\ell=j+1}^{k+1} \langle\be_\ell, \btheta_k^{(2)}\rangle
                   \langle\be_\ell, \bv_j^{(k)}\rangle \log(\langle\be_\ell, \bv_j^{(k)}\rangle)\biggr\} \\
  &= \EE\bigl(\ee^{\ii(\vartheta_0\cY^{(1)}_0+\cdots+\vartheta_{k-1}\cY^{(1)}_{k-1})}\bigr)
  \EE\bigl(\ee^{\ii\vartheta_k\tvare^{(1)}_k}\bigr) ,
 \end{align*}
 where we used that \ $\langle\be_\ell, \btheta_k^{(1)}\rangle = 0$ \ if \ $\ell = k + 1$, \ and \ $\vartheta_0 = \vartheta_1=\ldots=\vartheta_{k-1}=0$
  \ yields that \ $\btheta_k^{(1)}=\bzero\in\RR^{k+1}$.
\ Thus we obtain the independence of \ $(\cY^{(1)}_0, \ldots, \cY^{(1)}_{k-1})^\top$ \ and \ $\tvare^{(1)}_k$, \ and the characteristic function of \ $\tvare^{(1)}_k$ \ has the form
 \begin{align*}
  \EE\bigl(\ee^{\ii\vartheta_k\tvare^{(1)}_k}\bigr)
  &= \exp\biggl\{- C_1 (1 - m_\xi)
                   |\vartheta_k|
                   \Bigl(1 + \ii \frac{2}{\pi} \sign(\vartheta_k)
                                 \log(|\vartheta_k|)\Bigr)
                 + \ii C \vartheta_k  \\
  &\phantom{= \exp\biggl\{}
                 + \ii (1 - m_\xi) \sum_{j=0}^k \sum_{\ell=j+1}^{k+1} \langle\be_\ell, \btheta_k^{(2)}\rangle
                   \langle\be_\ell, \bv_j^{(k)}\rangle \log(\langle\be_\ell, \bv_j^{(k)}\rangle)\biggr\}
 \end{align*}
 \begin{align*}
  &= \exp\biggl\{- C_1 (1 - m_\xi)
                   |\vartheta_k|
                   \Bigl(1 + \ii \frac{2}{\pi} \sign(\vartheta_k)
                                 \log(|\vartheta_k|)\Bigr)
                 + \ii \vartheta_k (C + m_\xi\log(m_\xi))\biggr\}
 \end{align*}
 for \ $\vartheta_k \in \RR$, \ since
 \begin{align*}
  &\sum_{j=0}^k \sum_{\ell=j+1}^{k+1} \langle\be_\ell, \btheta_k^{(2)}\rangle
    \langle\be_\ell, \bv_j^{(k)}\rangle \log(\langle\be_\ell, \bv_j^{(k)}\rangle)
   = \sum_{\ell=1}^{k+1} \langle\be_\ell, \btheta_k^{(2)}\rangle
      \sum_{j=0}^{\ell-1} \langle\be_\ell, \bv_j^{(k)}\rangle \log(\langle\be_\ell, \bv_j^{(k)}\rangle) \\
  &= (- m_\xi \vartheta_k)
     \biggl(\frac{m_\xi^{k-1}}{1-m_\xi} \log\biggl(\frac{m_\xi^{k-1}}{1-m_\xi}\biggr)
            + \sum_{j=1}^{k-1} m_\xi^{k-j-1} \log(m_\xi^{k-j-1})\biggr) \\
  &\quad
     + \vartheta_k
       \biggl(\frac{m_\xi^k}{1-m_\xi} \log\biggl(\frac{m_\xi^k}{1-m_\xi}\biggr)
              + \sum_{j=1}^k m_\xi^{k-j} \log(m_\xi^{k-j})\biggr) \\
  &= \vartheta_k
     \biggl(\frac{m_\xi^k}{1-m_\xi} \log(m_\xi) + \sum_{j=1}^{k-1} m_\xi^{k-j} \log(m_\xi)\biggr)
   = \vartheta_k \biggl(\frac{m_\xi^k}{1-m_\xi} + \frac{m_\xi-m_\xi^k}{1-m_\xi}\biggr) \log(m_\xi) \\
  &= \vartheta_k \frac{m_\xi}{1-m_\xi} \log(m_\xi) .
 \end{align*}
Hence \ $\tvare^{(1)}_k$ \ is \ 1-stable (see, e.g., Sato \cite[Theorem 14.10]{Sato}).
In fact,
 \[
   \tvare^{(1)}_k - m_\xi(C+\log(m_\xi))  \distre (1 - m_\xi) (\cY^{(1)}_0 + \log(1 - m_\xi)) , \qquad k \in \NN,
 \]
 since for all \ $\vartheta_k \in \RR$,
 \begin{align*}
  \EE(\ee^{\ii\vartheta_k(1-m_\xi)\cY^{(1)}_0})
  &= \exp\Big\{   -C_1 (1-m_\xi) \vert\vartheta_k\vert \Bigl(1 + \ii \frac{2}{\pi} \sign(\vartheta_k)
                                 \log(|\vartheta_k|)\Bigr)  \\
  &\phantom{= \exp\Big\{\;}
                   + \ii C (1 - m_\xi) \vartheta_k  - \ii \vartheta_k (1-m_\xi) \log(1-m_\xi)
                                 \Big\} .
 \end{align*}
\proofend

\noindent{\bf Proof of Corollary \ref{aggr_copies1}.}
It follows from Theorem \ref{simple_aggregation1_stable_fd} and Corollary \ref{simple_aggregation1_stable_centering_fd}
 using the continuous mapping theorem.
\proofend

\noindent{\bf Proof of Theorem \ref{iterated_aggr_1}.}
In case of \ $\alpha \in (0, 1)$, \ by \eqref{centering01}, we have
 \[
   \lim_{N\to\infty} \frac{\nt N}{n^{\frac{1}{\alpha}}a_N} \EE(X_0 \bbone_{\{X_0\leq a_N\}})
   = \frac{\nt}{n^{\frac{1}{\alpha}}} \frac{\alpha}{1-\alpha}
   \to 0 \qquad \text{as \ $n \to \infty$,}
 \]
 hence, by Slutsky's lemma, \eqref{iterated_aggr_1_1} will be a consequence of \eqref{iterated_aggr_1_2}.

For each \ $n \in \NN$, \ by Corollary \ref{aggr_copies1} and by the continuous mapping theorem, we obtain
 \[
   \biggl(\frac{1}{n^{\frac{1}{\alpha}}a_N}
          \sum_{k=1}^\nt \sum_{j=1}^N X^{(j)}_k\biggr)_{t\in\RR_+}
   \distrf \biggl(\frac{1}{n^{\frac{1}{\alpha}}}
                 \sum_{k=1}^\nt
                  \biggl(\cY^{(\alpha)}_k + \frac{\alpha}{1-\alpha}\biggr)\biggr)_{t\in\RR_+} \qquad \text{as \ $N \to \infty$}
 \]
 in case of \ $\alpha \in (0, 1)$, \ and
 \[
   \biggl(\frac{1}{n^{\frac{1}{\alpha}}a_N}
          \sum_{k=1}^\nt \sum_{j=1}^N \bigl(X^{(j)}_k - \EE(X^{(j)}_k)\bigr)\biggr)_{t\in\RR_+}
   \distrf \biggl(\frac{1}{n^{\frac{1}{\alpha}}}
                 \sum_{k=1}^\nt
                  \biggl(\cY^{(\alpha)}_k + \frac{\alpha}{1-\alpha} \biggr)\biggr)_{t\in\RR_+} \qquad \text{as \ $N \to \infty$}
 \]
 in case of \ $\alpha \in (1, 2)$.
\ Consequently, in order to prove \eqref{iterated_aggr_1_2} and \eqref{iterated_aggr_1_4}, we need to show that for each \ $\alpha \in (0,1) \cup (1,2)$, \ we have
 \begin{equation}\label{conv_Z}
  \biggl(\frac{1}{n^{\frac{1}{\alpha}}}
                 \sum_{k=1}^\nt
                  \biggl(\cY^{(\alpha)}_k + \frac{\alpha}{1-\alpha}\biggr)\biggr)_{t\in\RR_+}
  \distrf \Bigl(\cZ^{(\alpha)}_t + \frac{\alpha}{1-\alpha} t\Bigr)_{t\in\RR_+} \qquad \text{as \ $n \to \infty$.}
 \end{equation}
For each \ $\alpha\in(0,1)\cup(1,2)$, \ $n \in \NN$, \ $d \in \NN$, \ $t_1, \ldots, t_d \in \RR_{++}$ \ with
 \ $t_1 < \ldots < t_d$ \ and \ $\vartheta_1, \ldots, \vartheta_d \in \RR$, \ we have
 \[
   \EE\biggl(\exp\biggl\{\ii \sum_{\ell=1}^d
                                \frac{\vartheta_\ell}{n^{\frac{1}{\alpha}}}
                                \sum_{k=\lfloor nt_{\ell-1}\rfloor+1}^{\lfloor nt_\ell\rfloor}
                                 \Bigl(\cY^{(\alpha)}_k + \frac{\alpha}{1-\alpha}\Bigr)\biggr\}\biggr)
   = \EE\Bigl(\exp\Bigl\{\ii \Bigl\langle n^{-\frac{1}{\alpha}} \btheta_n, \bcX_1^{(\lfloor nt_d\rfloor,\alpha)} + \frac{\alpha}{1-\alpha} \bone_{\lfloor nt_d\rfloor+1}\Bigr\rangle\Bigr\}\Bigr)
 \]
 with \ $t_0:=0$ \ and
 \[
   \btheta_n
   := \sum_{\ell=1}^d \vartheta_\ell \sum_{k=\lfloor nt_{\ell-1}\rfloor+1}^{\lfloor nt_\ell\rfloor} \be_{k+1} \in \RR^{\lfloor nt_d\rfloor+1} .
 \]
For each \ $\alpha\in(0,1)\cup(1,2)$, \ by the explicit form of the characteristic function of
 \ $\bcX_1^{(\lfloor nt_d\rfloor,\alpha)}$ \ given in Theorem \ref{simple_aggregation1_stable_fd},
 \begin{align*}
  &\EE\Bigl(\exp\Bigl\{\ii \Bigl\langle n^{-\frac{1}{\alpha}} \btheta_n, \bcX_1^{(\lfloor nt_d\rfloor,\alpha)} + \frac{\alpha}{1-\alpha} \bone_{\lfloor nt_d\rfloor+1}\Bigr\rangle\Bigr\}\Bigr) \\
  &= \exp\biggl\{- C_\alpha (1 - m_\xi^\alpha)
                  \sum_{j=0}^{\lfloor nt_d\rfloor}
                   \bigl|\bigl\langle n^{- \frac{1}{\alpha} } \btheta_n, \bv_j^{(\lfloor nt_d\rfloor)}\bigr\rangle\bigr|^\alpha
                    \left(1 - \ii \tan\left(\frac{\pi\alpha}{2}\right)
                                  \sign\bigl(\bigl\langle n^{- \frac{1}{\alpha}} \btheta_n,\bv_j^{(\lfloor nt_d\rfloor)}\bigr\rangle\bigr)\right)\biggr\} \\
  &= \exp\biggl\{- C_\alpha (1 - m_\xi^\alpha) n^{-1}
                  \sum_{j=0}^{\lfloor nt_d\rfloor}
                   \bigl|\bigl\langle\btheta_n, \bv_j^{(\lfloor nt_d\rfloor)}\bigr\rangle\bigr|^\alpha
                    \left(1 - \ii \tan\left(\frac{\pi\alpha}{2}\right)
                                  \sign\bigl(\bigl\langle\btheta_n,\bv_j^{(\lfloor nt_d\rfloor)}\bigr\rangle\bigr)\right)\biggr\} .
 \end{align*}
We have
 \begin{align*}
  \bigl\langle\btheta_n, \bv_0^{(\lfloor nt_d\rfloor)}\bigr\rangle
  &= \sum_{i=1}^d \vartheta_i \sum_{k=\lfloor nt_{i-1}\rfloor+1}^{\lfloor nt_i\rfloor} \bigl\langle\be_{k+1}, \bv_0^{(\lfloor nt_d\rfloor)}\bigr\rangle
   = \sum_{i=1}^d \vartheta_i \sum_{k=\lfloor nt_{i-1}\rfloor+1}^{\lfloor nt_i\rfloor} \frac{m_\xi^k}{(1-m_\xi^\alpha)^{\frac{1}{\alpha}}} \\
  &= \frac{1}{(1-m_\xi^\alpha)^{\frac{1}{\alpha}}(1-m_\xi)} \sum_{i=1}^d \vartheta_i (m_\xi^{\lfloor nt_{i-1}\rfloor+1} - m_\xi^{\lfloor nt_i\rfloor+1}) ,
 \end{align*}
 hence for each \ $\alpha\in(0, 1) \cup (1,2)$,
 \[
   n^{-1}
   \bigl|\bigl\langle\btheta_n, \bv_0^{(\lfloor nt_d\rfloor)}\bigr\rangle\bigr|^\alpha
   \left(1 - \ii \tan\left(\frac{\pi\alpha}{2}\right)
   \sign\bigl(\bigl\langle\btheta_n,\bv_0^{(\lfloor nt_d\rfloor)}\bigr\rangle\bigr)\right)
   \to 0 \qquad \text{as \ $n \to \infty$.}
 \]
The aim of the following discussion is to show that for each \ $\alpha\in(0, 1)\cup(1,2)$ \ and \ $\ell \in \{1, \ldots, d\}$,
 \begin{equation}\label{conv_alpha}
  n^{-1}
  \sum_{j=\lfloor nt_{\ell-1}\rfloor+1}^{\lfloor nt_\ell\rfloor}
   \bigl|\bigl\langle\btheta_n, \bv_j^{(\lfloor nt_d\rfloor)}\bigr\rangle\bigr|^\alpha
  \to \frac{(t_\ell-t_{\ell-1})|\vartheta_\ell|^\alpha}{(1-m_\xi)^\alpha}\qquad \text{as \ $n \to \infty$.}
 \end{equation}
Here, for each \ $\ell \in \{1, \ldots, d\}$ \ and \ $j \in \{\lfloor nt_{\ell-1}\rfloor + 1, \ldots, \lfloor nt_\ell\rfloor\}$,
 \begin{align}\label{help_7}
  \begin{split}
  \bigl\langle\btheta_n, \bv_j^{(\lfloor nt_d\rfloor)}\bigr\rangle
  &= \sum_{i=1}^d \vartheta_i \sum_{k=\lfloor nt_{i-1}\rfloor+1}^{\lfloor nt_i\rfloor} \bigl\langle\be_{k+1}, \bv_j^{(\lfloor nt_d\rfloor)}\bigr\rangle
   = \sum_{i=1}^d \vartheta_i \sum_{k=\lfloor nt_{i-1}\rfloor+1}^{\lfloor nt_i\rfloor} m_\xi^{k-j} \bbone_{\{k\geq j\}} \\
  &= \frac{1}{1-m_\xi} \biggl(\vartheta_\ell (1 - m_\xi^{\lfloor nt_\ell\rfloor-j+1}) + \sum_{i=\ell+1}^d \vartheta_i (m_\xi^{\lfloor nt_{i-1}\rfloor-j+1} - m_\xi^{\lfloor nt_i\rfloor-j+1})\biggr) .
  \end{split}
 \end{align}
In case of \ $\alpha \in (0, 1]$, \ we have
 \begin{align}\label{help_normlike_alpha01}
 |x|^\alpha - |y|^\alpha \leq |x+y|^\alpha \leq |x|^\alpha + |y|^\alpha , \qquad x, y \in \RR .
 \end{align}
In case of \ $\alpha \in (1, 2)$, \ by the mean value theorem and by \eqref{help_normlike_alpha01}, we have
 \[
   \bigl||x + y|^\alpha - |x|^\alpha\bigr|
   \leq \alpha |y| \max\{|x + y|^{\alpha-1}, |x|^{\alpha-1}\}
   \leq \alpha |y| (|x|^{\alpha-1} + |y|^{\alpha-1}) , \qquad x, y \in \RR .
 \]
Hence for each \ $\alpha \in (0, 2)$ \ and \ $x, y \in \RR$, \ we obtain
 \[
  |x|^\alpha - 2 |y| (|x|^{\alpha-1} + |y|^{\alpha-1})
  \leq |x+y|^\alpha
  \leq |x|^\alpha + 2 |y| (|x|^{\alpha-1} + |y|^{\alpha-1}) ,
 \]
 so, by \eqref{help_7} and the squeeze theorem, to prove \eqref{conv_alpha}, it is enough to check that
  \begin{gather}
   \frac{1}{n}
   \sum_{j=\lfloor nt_{\ell-1}\rfloor+1}^{\lfloor nt_\ell\rfloor}
   (1 - m_\xi^{\lfloor nt_\ell\rfloor-j+1})^\alpha
   \to t_\ell - t_{\ell-1} , \label{conv_alpha_02_1} \\
   \frac{1}{n}
   \sum_{j=\lfloor nt_{\ell-1}\rfloor+1}^{\lfloor nt_\ell\rfloor}
    \Biggl|\sum_{i=\ell+1}^d \vartheta_i (m_\xi^{\lfloor nt_{i-1}\rfloor-j+1} - m_\xi^{\lfloor nt_i\rfloor-j+1})\Biggr|^\alpha
   \to 0 , \label{conv_alpha_02_2} \\
   \frac{1}{n}
   \sum_{j=\lfloor nt_{\ell-1}\rfloor+1}^{\lfloor nt_\ell\rfloor}
    \Biggl|\sum_{i=\ell+1}^d \vartheta_i (m_\xi^{\lfloor nt_{i-1}\rfloor-j+1} - m_\xi^{\lfloor nt_i\rfloor-j+1})\Biggr|
    (1 - m_\xi^{\lfloor nt_\ell\rfloor-j+1})^{\alpha-1}
   \to 0 \label{conv_alpha_02_3}
  \end{gather}
  as \ $n \to \infty$.
\ Since \ $(1 - t)^\alpha = 1 - \alpha t + \oo(t)$ \ as \ $t \downarrow 0$, \ there exists \ $j_0\in\NN$ \ such that
 \ $|(1-m_\xi^j)^\alpha - 1 + \alpha m_\xi^j| \leq m_\xi^j$ \ for all \ $j\geq j_0$.
\ Hence
 \begin{align*}
  &\Biggl|\frac{1}{n}
          \sum_{j=\lfloor nt_{\ell-1}\rfloor+1}^{\lfloor nt_\ell\rfloor}
           \bigl(1 - m_\xi^{\lfloor nt_\ell\rfloor-j+1}\bigr)^\alpha
          - \frac{1}{n}
            \sum_{j=1}^{\lfloor nt_\ell\rfloor-\lfloor nt_{\ell-1}\rfloor}
             (1 - \alpha m_\xi^j)\Biggr| \\
 &= \Biggl|\frac{1}{n}
            \sum_{j=1}^{\lfloor nt_\ell\rfloor-\lfloor nt_{\ell-1}\rfloor}
             (1 - m_\xi^j)^\alpha
           - \frac{1}{n}
             \sum_{j=1}^{\lfloor nt_\ell\rfloor-\lfloor nt_{\ell-1}\rfloor}
              (1 - \alpha m_\xi^j)\Biggr|
\end{align*}
\begin{align*}
 &\leq \Biggl|\frac{1}{n}
              \sum_{j=1}^{j_0-1}
               \bigl((1 - m_\xi^j)^\alpha - 1 + \alpha m_\xi^j\bigr)\Biggr|
       + \frac{1}{n}
         \sum_{j=j_0}^{\lfloor nt_\ell\rfloor-\lfloor nt_{\ell-1}\rfloor} m_\xi^j  \\
 &\leq \Biggl|\frac{1}{n}
              \sum_{j=1}^{j_0-1}
               \bigl( (1 - m_\xi^j)^\alpha - 1 + \alpha m_\xi^j\bigr) \Biggr|
       + \frac{1}{n} \frac{m_\xi^{j_0}}{1-m_\xi}
   \to 0 \qquad \text{as \ $n \to \infty$.}
 \end{align*}
Thus
 \begin{align*}
  &\lim_{n\to\infty}
     \frac{1}{n} \sum_{j=\lfloor nt_{\ell-1}\rfloor+1}^{\lfloor nt_\ell\rfloor}
           \bigl(1 - m_\xi^{\lfloor nt_\ell\rfloor-j+1}\bigr)^\alpha
   = \lim_{n\to\infty}
         \frac{1}{n} \sum_{j=1}^{\lfloor nt_\ell\rfloor - \lfloor nt_{\ell-1}\rfloor}
                  (1 - \alpha m_\xi^j)\\
  &= \lim_{n\to\infty}
      \frac{1}{n} \Biggl(\lfloor nt_\ell\rfloor - \lfloor nt_{\ell-1}\rfloor
                          - \alpha \frac{m_\xi- m_\xi^{\lfloor nt_\ell\rfloor - \lfloor nt_{\ell-1}\rfloor+1}}{1-m_\xi}\Biggr)
   = t_\ell - t_{\ell-1} ,
 \end{align*}
 yielding \eqref{conv_alpha_02_1}.
In case of \ $\alpha \in (1, 2)$, \ for all \ $x_1, \ldots, x_k \in \RR$, \ we have
 \ $|x_1 + \cdots + x_k|^\alpha \leq k^{\alpha-1} (|x_1|^\alpha + \cdots + |x_k|^\alpha)$, \ hence, by \eqref{help_normlike_alpha01}, for each \ $\alpha \in (0, 2)$, \ we obtain
 \[
  |x_1 + \cdots + x_k|^\alpha \leq k (|x_1|^\alpha + \cdots + |x_k|^\alpha) , \qquad x_1, \ldots, x_k \in \RR .
 \]
Consequently, we have
 \begin{align*}
  &\frac{1}{n}
   \sum_{j=\lfloor nt_{\ell-1}\rfloor+1}^{\lfloor nt_\ell\rfloor}
    \Biggl| \sum_{i=\ell+1}^d \vartheta_i (m_\xi^{\lfloor nt_{i-1}\rfloor-j+1} - m_\xi^{\lfloor nt_i\rfloor-j+1})\Biggr|^\alpha \\
  &\leq \frac{1}{n}
        \sum_{j=\lfloor nt_{\ell-1}\rfloor+1}^{\lfloor nt_\ell\rfloor}
         (d - \ell)
         \sum_{i=\ell+1}^d
         |\vartheta_i|^\alpha
         \bigl(m_\xi^{\lfloor nt_{i-1}\rfloor-j+1} - m_\xi^{\lfloor nt_i\rfloor-j+1}\bigr)^\alpha \\
  &\leq \frac{d}{n}
        \sum_{i=\ell+1}^d
         |\vartheta_i|^\alpha
         \sum_{j=\lfloor nt_{\ell-1}\rfloor+1}^{\lfloor nt_\ell\rfloor}
          m_\xi^{(\lfloor nt_{i-1}\rfloor-j+1)\alpha}
   \leq \frac{d}{n} \sum_{i=\ell+1}^d |\vartheta_i|^\alpha
                     \sum_{k=0}^\infty m_\xi^{k\alpha}
  \to 0 \qquad \text{as \ $n \to \infty$,}
 \end{align*}
 yielding \eqref{conv_alpha_02_2}.
For each \ $n \in \NN$ \ and for each \ $j \in \{\lfloor nt_{\ell-1}\rfloor + 1, \ldots, \lfloor nt_\ell\rfloor\}$, \ we have
 \[
   (1 - m_\xi^{\lfloor nt_\ell\rfloor-j+1})^{\alpha-1}
   \leq \begin{cases}
         (1 - m_\xi)^{\alpha-1} & \text{if \ $\alpha \in (0, 1]$,} \\
         1 & \text{if \ $\alpha \in (1, 2)$,}
        \end{cases}
 \]
 and hence
 \begin{align*}
  &\frac{1}{n}
   \sum_{j=\lfloor nt_{\ell-1}\rfloor+1}^{\lfloor nt_\ell\rfloor}
    \Biggl| \sum_{i=\ell+1}^d \vartheta_i (m_\xi^{\lfloor nt_{i-1}\rfloor-j+1} - m_\xi^{\lfloor nt_i\rfloor-j+1})\Biggr|
    (1 - m_\xi^{\lfloor nt_\ell\rfloor-j+1})^{\alpha-1} \\
  &\leq \frac{(1 - m_\xi)^{\alpha-1} \vee 1}{n}
        \sum_{j=\lfloor nt_{\ell-1}\rfloor+1}^{\lfloor nt_\ell\rfloor}
         \sum_{i=\ell+1}^d
         |\vartheta_i|
         \bigl(m_\xi^{\lfloor nt_{i-1}\rfloor-j+1} - m_\xi^{\lfloor nt_i\rfloor-j+1}\bigr) \\
  &\leq \frac{(1 - m_\xi)^{\alpha-1} \vee 1}{n}
        \sum_{i=\ell+1}^d
         |\vartheta_i|
         \sum_{j=\lfloor nt_{\ell-1}\rfloor+1}^{\lfloor nt_\ell\rfloor}
          m_\xi^{\lfloor nt_{i-1}\rfloor-j+1}
   \leq \frac{(1 - m_\xi)^{\alpha-1} \vee 1}{n}
        \sum_{i=\ell+1}^d |\vartheta_i|
        \sum_{k=0}^\infty m_\xi^k
  \to 0
 \end{align*}
 as \ $n \to \infty$, \ yielding \eqref{conv_alpha_02_3}.
Thus we obtain \eqref{conv_alpha}.

Next we show that for each \ $\ell \in \{1, \ldots, d\}$, \ we have
 \[
   n^{-1}
   \sum_{j=\lfloor nt_{\ell-1}\rfloor+1}^{\lfloor nt_\ell\rfloor}
    \bigl|\bigl\langle\btheta_n, \bv_j^{(\lfloor nt_d\rfloor)}\bigr\rangle\bigr|^\alpha
    \sign\bigl(\bigl\langle\btheta_n, \bv_j^{(\lfloor nt_d\rfloor)}\bigr\rangle\bigr)
   \to \frac{(t_\ell-t_{\ell-1})|\vartheta_\ell|^\alpha}{(1-m_\xi)^\alpha} \sign(\vartheta_\ell)
 \]
 as \ $n \to \infty$.
\ If \ $\vartheta_\ell = 0$, \ then this readily follows from \eqref{help_7} and \eqref{conv_alpha_02_2}.
If \ $\vartheta_\ell \ne 0$, \ then we show that there exists \ $\tC_\ell \in \RR_{++}$ \ such that
 \begin{equation}\label{sign}
  \sign\bigl(\bigl\langle\btheta_n, \bv_j^{(\lfloor nt_d\rfloor)}\bigr\rangle\bigr)
  = \sign(\vartheta_\ell)
 \end{equation}
 for each \ $n \in \NN$ \ and for each \ $j \in \{\lfloor nt_{\ell-1}\rfloor + 1, \ldots, \lfloor nt_\ell\rfloor\}$ \ with \ $j < \lfloor nt_\ell\rfloor + 1 - \tC_\ell$.
\ First, observe that, by \eqref{help_7}, the inequality
 \begin{equation}\label{sign1}
   \Biggl|\sum_{i=\ell+1}^d
           \vartheta_i
           \bigl(m_\xi^{\lfloor nt_{i-1}\rfloor-j+1} - m_\xi^{\lfloor nt_i\rfloor-j+1}\bigr)\Biggr|
   < |\vartheta_\ell (1 - m_\xi^{\lfloor nt_{\ell}\rfloor - j+1})|
 \end{equation}
 implies \eqref{sign}.
Then we have
 \begin{align*}
  &\Biggl|\sum_{i=\ell+1}^d
           \vartheta_i
           \bigl(m_\xi^{\lfloor nt_{i-1}\rfloor-j+1} - m_\xi^{\lfloor nt_i\rfloor-j+1}\bigr)\Biggr|
   \leq \left(\max_{\ell+1\leq i\leq d} |\vartheta_i|\right)
        \sum_{i=\ell+1}^d
         \bigl(m_\xi^{\lfloor nt_{i-1}\rfloor-j+1} - m_\xi^{\lfloor nt_i\rfloor-j+1}\bigr) \\
  &= \left(\sum_{i=\ell+1}^d |\vartheta_i|\right)
     \bigl(m_\xi^{\lfloor nt_\ell\rfloor-j+1} - m_\xi^{\lfloor nt_d\rfloor-j+1}\bigr)
   \leq \left(\sum_{i=\ell+1}^d |\vartheta_i|\right)
        m_\xi^{\lfloor nt_\ell\rfloor-j+1} ,
 \end{align*}
 hence \eqref{sign1} is satisfied if
 \[
   \left(\sum_{i=\ell+1}^d |\vartheta_i|\right)
        m_\xi^{\lfloor nt_\ell\rfloor-j+1}
   < |\vartheta_\ell| (1 - m_\xi^{\lfloor nt_{\ell}\rfloor - j+1}) ,
 \]
 which is satisfied if
 \[
   m_\xi^{\lfloor nt_\ell\rfloor - j+1}
   < \frac{|\vartheta_\ell|}{|\vartheta_\ell|+\cdots+|\vartheta_d|} ,
 \]
 or equivalently, if
 \[
   j < \lfloor nt_\ell\rfloor + 1 - \tC_\ell \qquad
   \text{with \ $\tC_\ell :=\log\biggl(\frac{|\vartheta_\ell|}{|\vartheta_\ell|+\cdots+|\vartheta_d|}\biggr) \bigg/ \log(m_\xi) \in \RR_{++}$.}
 \]
Hence, for \ $\vartheta_\ell \ne 0$, \ $n \in \NN$, \ and \ $j \in \{\lfloor nt_{\ell-1}\rfloor + 1, \ldots, \lfloor nt_\ell\rfloor\}$
 \ with \ $j < \lfloor nt_\ell\rfloor + 1 - \tC_\ell$, \ we have \eqref{sign}.
Moreover, for each \ $n \in \NN$ \ and \ $j \in \{\lfloor nt_{\ell-1}\rfloor + 1, \ldots, \lfloor nt_\ell\rfloor\}$, \ by \eqref{help_7}, we have
 \begin{align*}
  \bigl|\bigl\langle\btheta_n, \bv_j^{(\lfloor nt_d\rfloor)}\bigr\rangle\bigr|
  &\leq \frac{1}{1-m_\xi} \biggl(|\vartheta_\ell| (1 - m_\xi^{\lfloor nt_\ell\rfloor-j+1}) + \sum_{i=\ell+1}^d |\vartheta_i| (m_\xi^{\lfloor nt_{i-1}\rfloor-j+1} - m_\xi^{\lfloor nt_i\rfloor-j+1})\biggr) \\
  &\leq \frac{1}{1-m_\xi} \sum_{i=\ell}^d |\vartheta_i| ,
 \end{align*}
 yielding that
  \begin{align*}
   & \left\vert n^{-1}
    \sum_{j=\lfloor nt_{\ell}\rfloor+1 - \widetilde C_\ell}^{\lfloor nt_\ell\rfloor}
     \bigl|\bigl\langle\btheta_n, \bv_j^{(\lfloor nt_d\rfloor)}\bigr\rangle\bigr|^\alpha
     \sign\bigl(\bigl\langle\btheta_n, \bv_j^{(\lfloor nt_d\rfloor)}\bigr\rangle\bigr)
     \right\vert
    \leq n^{-1} \sum_{j= \lfloor nt_\ell\rfloor +1- \widetilde C_\ell}^{\lfloor nt_\ell\rfloor}
                     \bigl|\bigl\langle\btheta_n, \bv_j^{(\lfloor nt_d\rfloor)}\bigr\rangle\bigr|^\alpha\\
   &\leq \frac{\widetilde C_\ell}{n(1-m_\xi)^\alpha} \left(\sum_{i=\ell}^d |\vartheta_i|\right)^\alpha
    \to 0\qquad \text{as \ $n\to\infty$.}
  \end{align*}
Consequently, by \eqref{conv_alpha},
 \begin{align*}
  &\lim_{n\to\infty}
    n^{-1}
    \sum_{j=\lfloor nt_{\ell-1}\rfloor+1}^{\lfloor nt_\ell\rfloor}
     \bigl|\bigl\langle\btheta_n, \bv_j^{(\lfloor nt_d\rfloor)}\bigr\rangle\bigr|^\alpha
     \sign\bigl(\bigl\langle\btheta_n, \bv_j^{(\lfloor nt_d\rfloor)}\bigr\rangle\bigr) \\
  &= \lim_{n\to\infty}
       n^{-1}
      \sum_{\lfloor nt_{\ell-1}\rfloor+1\leq j<\lfloor nt_\ell\rfloor+1-\tC_\ell}
       \bigl|\bigl\langle\btheta_n, \bv_j^{(\lfloor nt_d\rfloor)}\bigr\rangle\bigr|^\alpha
       \sign\bigl(\bigl\langle\btheta_n, \bv_j^{(\lfloor nt_d\rfloor)}\bigr\rangle\bigr) \\
  &= \lim_{n\to\infty}
     n^{-1}
      \sum_{\lfloor nt_{\ell-1}\rfloor+1\leq j<\lfloor nt_\ell\rfloor+1-\tC_\ell}
       \bigl|\bigl\langle\btheta_n, \bv_j^{(\lfloor nt_d\rfloor)}\bigr\rangle\bigr|^\alpha
       \sign(\vartheta_\ell) \\
  &= \lim_{n\to\infty}
      n^{-1}
      \sum_{j=\lfloor nt_{\ell-1}\rfloor+1}^{\lfloor nt_\ell\rfloor}
       \bigl|\bigl\langle\btheta_n, \bv_j^{(\lfloor nt_d\rfloor)}\bigr\rangle\bigr|^\alpha
       \sign(\vartheta_\ell)
   = \frac{(t_\ell-t_{\ell-1})|\vartheta_\ell|^\alpha}{(1-m_\xi)^\alpha} \sign(\vartheta_\ell) ,
 \end{align*}
 as desired.
We conclude for all \ $\alpha \in (0, 1) \cup (1, 2)$,
 \begin{align*}
  &\EE\biggl(\exp\biggl\{\ii \sum_{\ell=1}^d
                                \frac{\vartheta_\ell}{n^{\frac{1}{\alpha}}}
                                \sum_{k=\lfloor nt_{\ell-1}\rfloor+1}^{\lfloor nt_\ell\rfloor}
                                 \Bigl(\cY^{(\alpha)}_k + \frac{\alpha}{1-\alpha}\Bigr)\biggr\}\biggr) \\
  &=\EE\Bigl(\exp\Bigl\{\ii \Bigl\langle n^{- \frac{1}{\alpha} } \btheta_n, \bcX_1^{(\lfloor nt_d\rfloor,\alpha)} + \frac{\alpha}{1-\alpha} \bone_{\lfloor nt_d\rfloor+1}\Bigr\rangle\Bigr\}\Bigr) \\
  &\to \exp\biggl\{-C_\alpha \frac{1-m_\xi^\alpha}{(1-m_\xi)^\alpha}
                    \sum_{\ell=1}^d
                     (t_\ell - t_{\ell-1}) |\vartheta_\ell|^\alpha
                     \left(1 - \ii \tan\left(\frac{\pi\alpha}{2}\right) \sign(\vartheta_\ell)\right)\biggr\} \\
  &= \EE\biggl(\exp\biggl\{\ii \sum_{\ell=1}^d \vartheta_\ell \Bigl(\Bigl(\cZ_{t_\ell}^{(\alpha)} + \frac{\alpha}{1-\alpha} t_\ell\Bigl) - \Bigr(\cZ_{t_{\ell-1}}^{(\alpha)} + \frac{\alpha}{1-\alpha} t_{\ell-1}\Bigr)\Bigr)\biggr\}\biggr) \qquad \text{as \ $n \to \infty$.}
 \end{align*}
By the continuity theorem, we obtain for all \ $\alpha \in (0, 1) \cup (1, 2)$,
 \[
   \Biggl(\frac{1}{n^{\frac{1}{\alpha}}}
                 \sum_{k=\lfloor nt_{\ell-1}\rfloor+1}^{\lfloor nt_\ell\rfloor}
                  \Bigl(\cY^{(\alpha)}_k + \frac{\alpha}{1-\alpha}\Bigr)\Biggr)_{\ell\in\{1,\ldots,d\}}
   \distr \Bigl(\Bigl(\cZ_{t_\ell}^{(\alpha)} + \frac{\alpha}{1-\alpha} t_\ell\Bigl) - \Bigr(\cZ_{t_{\ell-1}}^{(\alpha)} + \frac{\alpha}{1-\alpha} t_{\ell-1}\Bigr)\Bigr)_{\ell\in\{1,\ldots,d\}}
 \]
 as \ $n \to \infty$, \ hence the continuous mapping theorem yields \eqref{conv_Z},
 and we finished the proofs of \eqref{iterated_aggr_1_1}, \eqref{iterated_aggr_1_2} and \eqref{iterated_aggr_1_4}.

Now we turn to prove \eqref{iterated_aggr_1_3}.
For each \ $n \in \NN$, \ by Corollary \ref{aggr_copies1} and by the continuous mapping theorem, in case of \ $\alpha = 1$, \ we obtain
 \[
   \biggl(\frac{1}{n\log(n)a_N}
          \sum_{k=1}^\nt \sum_{j=1}^N \bigl(X^{(j)}_k - \EE\bigl(X^{(j)}_k \bbone_{\{X^{(j)}_k\leq a_N\}}\bigr)\bigr)\biggr)_{t\in\RR_+}
   \distrf \biggl(\frac{1}{n\log(n)}
                 \sum_{k=1}^\nt
                  \cY^{(1)}_k \biggr)_{t\in\RR_+}
 \]
 as \ $N \to \infty$.
\ Consequently, in order to prove \eqref{iterated_aggr_1_3}, we need to show that
 \begin{equation}\label{conv_Z_1}
  \biggl(\frac{1}{n\log(n)}
                 \sum_{k=1}^\nt
                  \cY^{(1)}_k \biggr)_{t\in\RR_+}
  \distrf (t)_{t\in\RR_+} \qquad \text{as \ $n \to \infty$.}
 \end{equation}
Since the limit in \eqref{conv_Z_1} is deterministic, by van der Vaart \cite[Theorem 2.7, part (vi)]{Vaa},
 it is enough to show that for each \ $t \in \RR_+$, \ we have
 \begin{equation}\label{conv_Z_1+}
  \frac{1}{n\log(n)}
  \sum_{k=1}^\nt \cY^{(1)}_k
  \distr t \qquad \text{as \ $n \to \infty$.}
 \end{equation}
For each \ $n \in \NN$, \ $t \in \RR_{+}$, \ and \ $\vartheta \in \RR$, \ we have
 \[
   \EE\biggl(\exp\biggl\{\ii \frac{\vartheta}{n\log(n)}
                             \sum_{k=1}^\nt
                              \cY^{(1)}_k  \biggr\}\biggr)
   = \EE\biggl(\exp\biggl\{\ii \biggl\langle\frac{\vartheta\bone_{\nt+1}}{n\log(n)} , \bcX_1^{(\nt,1)}\biggr\rangle\biggr\}\biggr) .
 \]
By the explicit form of the characteristic function of \ $\bcX_1^{(\lfloor nt\rfloor,1)}$ \ given in Theorem \ref{simple_aggregation1_stable_fd},
 \begin{align*}
  &\EE\biggl(\exp\biggl\{\ii \biggl\langle\frac{\vartheta\bone_{\nt+1}}{n\log(n)} , \bcX_1^{(\nt,1)}\biggr\rangle\biggr\}\biggr) \\
  &= \exp\biggl\{- C_1 (1 - m_\xi)
                   \sum_{j=0}^\nt
                    \biggl|\biggl\langle\frac{\vartheta\bone_{\nt+1}}{n\log(n)},\bv_j^{(\nt)}\biggr\rangle\biggr|
                   \biggl(1 + \ii \frac{2}{\pi} \sign\biggl(\biggl\langle\frac{\vartheta\bone_{\nt+1}}{n\log(n)},\bv_j^{(\nt)}\biggr\rangle\biggr) \\
  &\phantom{= \exp\biggl\{- C_1 (1 - m_\xi)
                   \sum_{j=0}^\nt
                    \biggl|\biggl\langle\frac{\vartheta\bone_{\nt+1}}{n\log(n)},\bv_j^{(\nt)}\biggr\rangle\biggr|
                   \biggl(1 +=}
                                 \times
                                 \log\biggl(\biggl|\biggl\langle\frac{\vartheta\bone_{\nt+1}}{n\log(n)},\bv_j^{(\nt)}\biggr\rangle\biggr|\biggr)\biggr) \\[2mm]
  &\phantom{= \exp\biggl\{}
                 + \ii C \biggl\langle\frac{\vartheta\bone_{\nt+1}}{n\log(n)},
 \bone_{\nt+1}\biggr\rangle\\
  &\phantom{= \exp\biggl\{}
                 + \ii (1 - m_\xi) \sum_{j=0}^\nt \sum_{\ell=j+1}^{\nt+1} \biggl\langle\be_\ell, \frac{\vartheta\bone_{\nt+1}}{n\log(n)}\biggr\rangle
                   \bigl\langle\be_\ell, \bv_j^{(\nt)}\bigr\rangle \log\bigl(\bigl\langle\be_\ell, \bv_j^{(\nt)}\bigr\rangle\bigr)\biggr\}
  \end{align*}
 \begin{align*}
  &= \exp\biggl\{- \frac{C_1(1-m_\xi)|\vartheta|}{n\log(n)}
                   \sum_{j=0}^\nt
                    \bigl\langle\bone_{\nt+1}, \bv_j^{(\nt)}\bigr\rangle \\
  &\phantom{= \exp\biggl\{- \frac{C_1(1-m_\xi)|\vartheta|}{n\log(n)}
                   \sum_{j=0}^\nt\;}
                                 \times
                                 \biggl(1 + \ii \frac{2}{\pi} \sign(\vartheta)
                                 \log\biggl(\frac{|\vartheta|}{n\log(n)}
                                 \bigl\langle\bone_{\nt+1},\bv_j^{(\nt)}\bigr\rangle\biggr)\biggr) \\[2mm]
  &\phantom{= \exp\biggl\{}
                 + \ii C \frac{\vartheta}{n\log(n)} \bigl\langle\bone_{\nt+1} ,
\bone_{\nt+1}\bigr\rangle \\[2mm]
  &\phantom{= \exp\biggl\{}
                 + \ii \frac{(1 - m_\xi)\vartheta}{n\log(n)} \sum_{j=0}^\nt \sum_{\ell=j+1}^{\nt+1}
                   \bigl\langle\be_\ell, \bv_j^{(\nt)}\bigr\rangle \log\bigl(\bigl\langle\be_\ell, \bv_j^{(\nt)}\bigr\rangle\bigr)\biggr\} \\
  &\to \ee^{\ii t\vartheta}
 \end{align*}
 as \ $n \to \infty$ \ for each \ $\vartheta\in\RR$.
\ Indeed,
 \[
   \frac{1}{n\log(n)} \bigl\langle\bone_{\nt+1}, \bone_{\nt+1}\bigr\rangle
   = \frac{\nt+1}{n\log(n)} \to 0 \qquad \text{as \ $n \to \infty$,}
 \]
 and
 \begin{align*}
  &\biggl|\frac{1}{n\log(n)} \sum_{j=0}^\nt \sum_{\ell=j+1}^{\nt+1}
                   \bigl\langle\be_\ell, \bv_j^{(\nt)}\bigr\rangle
                   \log\bigl(\bigl\langle\be_\ell, \bv_j^{(\nt)}\bigr\rangle\bigr)\biggr|
   = \biggl|\frac{1}{n\log(n)} \sum_{j=0}^\nt \sum_{\ell=j+1}^{\nt+1}
                   m_\xi^{\ell-j-1}
                   \log\bigl(m_\xi^{\ell-j-1}\bigr)\biggr| \\
  &\leq \frac{|\log(m_\xi)|}{n\log(n)} \sum_{j=0}^\nt \sum_{\ell=j+1}^{\nt+1}
                   (\ell-j-1) m_\xi^{\ell-j-1}
   \leq \frac{|\log(m_\xi)|}{n\log(n)} \sum_{j=0}^\nt \sum_{\ell=j+1}^\infty
                   (\ell-j-1) m_\xi^{\ell-j-1} \\
  &= \frac{|\log(m_\xi)|(\nt+1)}{n\log(n)} \sum_{k=0}^\infty
                   k m_\xi^k
   = \frac{m_\xi|\log(m_\xi)|(\nt+1)}{(1-m_\xi)^2n\log(n)}
   \to 0 \qquad \text{as \ $n \to \infty$,}
 \end{align*}
 and
 \begin{equation}\label{help}
   \frac{1}{n\log(n)}
                   \sum_{j=0}^\nt
                    \bigl\langle\bone_{\nt+1}, \bv_j^{(\nt)}\bigr\rangle
   = \frac{1}{n\log(n)}
                   \sum_{j=0}^\nt \frac{1-m_\xi^{\nt-j+1}}{1-m_\xi}
   \leq \frac{\nt+1}{(1-m_\xi)n\log(n)}
   \to 0
 \end{equation}
 as \ $n \to \infty$, \ and
 \[
   - \frac{C_1(1-m_\xi)|\vartheta|}{n\log(n)}
                   \sum_{j=0}^\nt
                    \bigl\langle\bone_{\nt+1}, \bv_j^{(\nt)}\bigr\rangle
     \ii \frac{2}{\pi} \sign(\vartheta)
                                 \log\biggl(\frac{|\vartheta|}{n\log(n)}
                                 \bigl\langle\bone_{\nt+1},\bv_j^{(\nt)}\bigr\rangle\biggr)
   \to \ii t \vartheta
 \]
 as \ $n \to \infty$, \ since
 \begin{align*}
   &\biggl|\log\biggl(\frac{|\vartheta|}{\log(n)}
                \bigl\langle\bone_{\nt+1},\bv_j^{(\nt)}\bigr\rangle\biggr)\biggr|
    = \biggl|\log\biggl(\frac{|\vartheta|}{\log(n)}
                \frac{1-m_\xi^{\nt-j+1}}{1-m_\xi}\biggr)\biggr| \\
   &= \biggl|\log(|\vartheta|) - \log(\log(n)) +
    \log\biggl(\frac{1-m_\xi^{\nt-j+1}}{1-m_\xi}\biggr)\biggr|
   \leq |\log(|\vartheta|)| + \log(\log(n)) + |\log(1-m_\xi)| ,
 \end{align*}
 hence, by \eqref{help},
 \begin{align*}
  &\biggl|\frac{1}{n\log(n)}
   \sum_{j=0}^\nt
    \bigl\langle\bone_{\nt+1}, \bv_j^{(\nt)}\bigr\rangle
     \log\biggl(\frac{|\vartheta|}{\log(n)}
                \bigl\langle\bone_{\nt+1},\bv_j^{(\nt)}\bigr\rangle\biggr)\biggr| \\
  &\leq \frac{\nt+1}{(1-m_\xi)n\log(n)} \bigl(|\log(|\vartheta|)| + \log(\log(n)) + |\log(1-m_\xi)|\bigr)
   \to 0 \qquad \text{as \ $n \to \infty$,}
 \end{align*}
 and
 \begin{align*}
  \frac{C_1(1-m_\xi)|\vartheta|}{n\log(n)}
                   \sum_{j=0}^\nt
                    \bigl\langle\bone_{\nt+1}, \bv_j^{(\nt)}\bigr\rangle
     \ii \frac{2}{\pi} \sign(\vartheta) \log(n)
  = \ii \frac{(1-m_\xi)\vartheta}{n}
                   \sum_{j=0}^\nt \frac{1-m_\xi^{\nt-j+1}}{1-m_\xi}
  \to \ii t \vartheta
 \end{align*}
 as \ $n \to \infty$.
\ By the continuity theorem, we obtain \eqref{conv_Z_1+}, hence we finished the proof of \eqref{iterated_aggr_1_3}.
\proofend

\vspace*{5mm}

\appendix

\vspace*{5mm}

\noindent{\bf\Large Appendices}

\section{A version of the continuous mapping theorem}\label{App_cont_map_theorem}

If \ $\xi$ \ and \ $\xi_n$,
 \ $n \in \NN$, \ are random elements with values in a metric space \ $(E, d)$,
 \ then we also denote by \ $\xi_n \distr \xi$ \ the weak convergence of the
 distributions of \ $\xi_n$ \ on the space \ $(E, \cB(E))$ \ towards the
 distribution of \ $\xi$ \ on the space \ $(E, \cB(E))$ \ as \ $n \to \infty$,
 \ where \ $\cB(E)$ \ denotes the Borel \ $\sigma$-algebra on \ $E$ \ induced by
 the given metric \ $d$.

 The following version of the continuous mapping theorem can be found for
 example in Theorem 3.27 of Kallenberg \cite{Kal}.

\begin{Lem}\label{Lem_Kallenberg}
Let \ $(S, d_S)$ \ and \ $(T, d_T)$ \ be metric spaces and
 \ $(\xi_n)_{n \in \NN}$, \ $\xi$ \ be random elements with values in \ $S$
 \ such that \ $\xi_n \distr \xi$ \ as \ $n \to \infty$.
\ Let \ $f : S \to T$ \ and \ $f_n : S \to T$, \ $n \in \NN$, \ be measurable
 mappings and \ $C \in \cB(S)$ \ such that \ $\PP(\xi \in C) = 1$ \ and
 \ $\lim_{n \to \infty} d_T(f_n(s_n), f(s)) = 0$ \ if
 \ $\lim_{n \to \infty} d_S(s_n,s) = 0$ \ and \ $s \in C$, \ $s_n\in S$, \ $n\in\NN$.
\ Then \ $f_n(\xi_n) \distr f(\xi)$ \ as \ $n \to \infty$.
\end{Lem}

We will use the following corollary of this lemma several
 times.

\begin{Lem}\label{Conv2Funct}
Let \ $d \in \NN$, \ and let \ $(\bcU_t)_{t \in \RR_+}$,
 \ $(\bcU^{(n)}_t)_{t \in \RR_+}$, \ $n \in \NN$, \ be \ $\RR^d$-valued stochastic
 processes with c\`adl\`ag paths such that \ $\bcU^{(n)} \distr \bcU$ \ as $n \to \infty$.
\ Let \ $\Phi : \DD(\RR_+, \RR^d) \to \DD(\RR_+, \RR^d)$ \ and
 \ $\Phi_n : \DD(\RR_+, \RR^d) \to \DD(\RR_+, \RR^d)$, \ $n \in \NN$, \ be measurable
 mappings such that \ $\Phi_n(f_n) \to \Phi(f)$ \ in \ $\DD(\RR_+, \RR^d)$ \ as \ $n \to \infty$ \ whenever
 \ $f_n \to f$ \ in \ $\DD(\RR_+, \RR^d)$ \ as \ $n \to \infty$ \ with \ $f, f_n \in \DD(\RR_+, \RR^d)$, \ $n \in \NN$.
\ Then \ $\Phi_n(\bcU^{(n)}) \distr \Phi(\bcU)$ \ as \ $n \to \infty$.
\end{Lem}

\section{The underlying space and vague convergence}\label{App_vague}

For each \ $d \in \NN$, \ put \ $\RR_0^d := \RR^d \setminus \{\bzero\}$, \ and denote by
 \ $\cB(\RR_0^d)$ \ the Borel \ $\sigma$-algebra of \ $\RR_0^d$ \ induced by the metric \ $\varrho : \RR^d_0 \times \RR^d_0 \to \RR_+$, \ given by
 \begin{align}\label{S_metric}
  \varrho(\bx, \by) := \min\{\|\bx- \by\|, 1\} + \biggl|\frac{1}{\|\bx\|} - \frac{1}{\|\by\|}\biggr| , \qquad \bx, \by \in \RR^d_0 .
 \end{align}

\begin{Lem}\label{Lem_S_top}
The set \ $\RR^d_0$ \ furnished with the metric \ $\varrho$ \ given in \eqref{S_metric} is a complete separable metric space, and \ $B \subset \RR^d_0$ \ is bounded with respect to the metric \ $\varrho$ \ if and only if \ $B$ \ is separated from the origin \ $\bzero \in \RR^d$, \ i.e., there exists \ $\vare \in \RR_{++}$ \ such that \ $B \subset \{\bx \in \RR^d_0 : \|\bx\| > \vare\}$.
Moreover, the topology and the Borel \ $\sigma$-algebra \ $\cB(\RR_0^d )$ \ on \ $\RR_0^d$ \ induced by the metric \ $\varrho$ \ coincides with
 the topology and the Borel \ $\sigma$-algebra on \ $\RR_0^d$ \ induced by the usual metric \ $d(\bx, \by) := \|\bx-\by\|$, \ $\bx, \by \in \RR_0^d$,
 \ respectively.
\end{Lem}

\noindent{\bf Proof.}
First, we check that \ $\RR^d_0$ \ furnished with the metric \ $\varrho$ \ is a complete separable metric space.
If \ $(\bx_n)_{n\in\NN}$ \ is a Cauchy sequence in \ $\RR^d_0$, \ then for all \ $\vare \in (0, 1)$, \ there exists an \ $N_\vare\in\NN$ \ such that $\varrho(\bx_n, \bx_m)<\vare$ for $n,m\geq N_\vare$.
Hence $\|\bx_n-\bx_m\|<\vare$ and $\left\vert \frac{1}{\|\bx_n\|} - \frac{1}{\|\bx_m\|}\right\vert<\vare$ for $n,m\geq N_\vare$,
 i.e., $(\bx_n)_{n\in\NN}$ and $(1/\|\bx_n\|)_{n\in\NN}$ are Cauchy sequences in $\RR^d$ and in $\RR$, respectively.
 Consequently, there exists an $\bx\in\RR^d$ such that $\lim_{n\to\infty} \|\bx_n - \bx\| = 0$
 and $\frac{1}{\|\bx_n\|}$ being convergent as $n\to\infty$, yielding that $\bx\ne\bzero$, and so $\bx\in\RR^d_0$.
By the continuity of the norm, $\lim_{n\to\infty} \varrho(\bx_n, \bx) = 0$, as desired.
The separability of $\RR^d_0$ readily follows, since $\RR^d_0\cap \QQ^d$ is a countable everywhere dense subset of $\RR^d_0$.

\noindent
Next, we check that $B \subset \RR^d_0$ is bounded with respect to the metric $\varrho$ if and only if
 there exists $\vare\in\RR_{++}$ such that $B \subset \{\bx\in \RR^d_0 : \|\bx\|>\vare\}$.
If $B \subset \RR^d_0$ is bounded, then there exists $r>0$ such that $\varrho(\bx,\be_1)<r$, $\bx\in B$,
 yielding that $\vert \frac{1}{\|\bx\|} - 1\vert < r$, $\bx\in B$, and then $\|\bx\|>\frac{1}{1+r}$, $\bx\in B$,
 so one can choose $\vare=\frac{1}{1+r}$.
If there exists $\vare>0$ such that $B \subset \{\bx\in \RR^d_0 : \|\bx\|>\vare\}$, then
 $\varrho(\bx,\be_1) = \min\{\|\bx-\be_1\|,1\} + \vert \frac{1}{\|\bx\|} - 1\vert \leq 1+\frac{1}{\vare}+1$, $\bx\in B$.
\proofend

Since \ $\RR^d_0$ \ is locally compact, second countable and Hausdorff, one could choose a metric such that the relatively compact sets
 are precisely the bounded ones, see Kallenberg \cite[page 18]{kallenberg:2017}.
The metric \ $\varrho$ \ does not have this property, but we do not need it.

Write \ $(\RR^d_0)\hspace*{.5mm}\widehat{}$ \ for the class of bounded Borel sets with respect to the metric
 \ $\varrho$ \ given in \eqref{S_metric}.
A measure \ $\nu$ \ on \ $(\RR_0^d, \cB(\RR_0^d))$ \ is said to be locally finite if \ $\nu(B) < \infty$ \ for every \ $B \in (\RR^d_0)\hspace*{.5mm}\widehat{}$, \ and write \ $\cM(\RR^d_0)$ \ for the class of locally finite measures on \ $(\RR_0^d, \cB(\RR_0^d))$.

Write \ $\hC_{\RR^d_0}$ \ for the class of bounded, continuous functions \ $f : \RR^d_0 \to \RR_+$ \ with bounded support.
Hence, if \ $f \in \hC_{\RR^d_0}$, \ then there exist an \ $\vare \in \RR_{++}$ \ such that \ $f(\bx) = 0$ \ for all \ $\bx \in \RR^d_0$ \ with \ $\|\bx\| \leq \vare$.
\ The vague topology on \ $\cM(\RR^d_0)$ \ is constructed as in Chapter 4 in Kallenberg \cite{kallenberg:2017}.
The associated notion of vague convergence of a sequence \ $(\nu_n)_{n\in\NN}$ \ in \ $\cM(\RR^d_0)$ \ towards $\nu\in\cM(\RR^d_0)$,
 denoted by \ $\nu_n \distrv \nu$ \ as \ $n\to\infty$, \ is defined by the condition \ $\nu_n(f) \to \nu(f)$ \ as \ $n\to\infty$
 \ for all \ $f \in \hC_{\RR^d_0}$, \ where \ $\kappa(f) := \int_{\RR^d_0} f(\bx) \, \kappa(\dd\bx)$ \ for each \ $\kappa \in \cM(\RR^d_0)$.

If \ $\nu$ \ is a measure on \ $(\RR_0^d, \cB(\RR_0^d))$, \ then
 \ $B \in \cB(\RR_0^d)$ \ is called a \ $\nu$-continuity set if \ $\nu(\partial B) = 0$, \ and the class of bounded \ $\nu$-continuity sets will be denoted by \ $(\RR^d_0)_\nu\hspace*{-1.2mm}\widehat{}\hspace*{1.2mm}$.
\ The following statement is an analogue of the portmanteau theorem for vague
 convergence, see, e.g., Kallenberg \cite[15.7.2]{kallenberg:1983}.

\begin{Lem}\label{portmanteau}
Let \ $\nu, \nu_n \in \cM(\RR^d_0)$, $n\in\NN$.
\ Then the following statements are equivalent:
\renewcommand{\labelenumi}{{\rm(\roman{enumi})}}
 \begin{enumerate}
  \item
   $\nu_n \distrv \nu$ \ as \ $n \to \infty$,
  \item
   $\nu_n(B) \to \nu(B)$ \ as \ $n \to \infty$ \ for all \ $B \in (\RR^d_0)_\nu\hspace*{-1.2mm}\widehat{}\hspace*{1.2mm}$.
 \end{enumerate}
\end{Lem}

The following statement is an analogue of the continuous mapping theorem
 for vague convergence, see, e.g., Kallenberg \cite[15.7.3]{kallenberg:1983}.
Write \ $\mathrm{D}_f$ \ for the set of discontinuities of a function \ $f : \RR^d_0 \to \RR$.

\begin{Lem}\label{cmt}
Let \ $\nu, \nu_n \in \cM(\RR^d_0)$, \ $n\in\NN$, \ with \ $\nu_n \distrv \nu$ \ as \ $n \to \infty$.
\ Then \ $\nu_n(f) \to \nu(f)$ \ as \ $n \to \infty$ \ for every bounded measurable function \ $f : \RR^d_0 \to \RR_+$ \ with bounded support satisfying \ $\nu(\mathrm{D}_f) = 0$.
\end{Lem}

\section{Regularly varying distributions}\label{App_reg_var_distr}

First, we recall the notions of slowly varying and regularly varying functions, respectively.

\begin{Def}
A measurable function \ $U: \RR_{++} \to \RR_{++}$ \ is called regularly varying at infinity with
 index \ $\rho \in \RR$ \ if for all \ $c \in \RR_{++}$,
 \[
   \lim_{x\to\infty} \frac{U(cx)}{U(x)} = c^\rho .
 \]
In case of \ $\rho = 0$, \ we call \ $U$ \ slowly varying at infinity.
\end{Def}

\begin{Def}
A random variable \ $Y$ \ is called regularly varying with index \ $\alpha \in \RR_{++}$ \ if \ $\PP(|Y| > x) \in \RR_{++}$ \ for all \ $x \in \RR_{++}$, \ the function \ $\RR_{++} \ni x \mapsto \PP(|Y| > x) \in \RR_{++}$ \ is regularly varying at infinity with index \ $-\alpha$, \ and a tail-balance condition holds:
 \begin{equation}\label{TB}
   \lim_{x\to\infty} \frac{\PP(Y>x)}{\PP(|Y|>x)} = p , \qquad
   \lim_{x\to\infty} \frac{\PP(Y\leq-x)}{\PP(|Y|>x)} = q ,
 \end{equation}
 where \ $p + q = 1$.
\end{Def}

\begin{Rem}\label{tail-balance}
In the tail-balance condition \eqref{TB}, the second convergence can be replaced by
 \begin{equation}\label{tail-balance+}
  \lim_{x\to\infty} \frac{\PP(Y<-x)}{\PP(|Y|>x)} = q .
 \end{equation}
Indeed, if \ $Y$ \ is regularly varying with index \ $\alpha \in \RR_{++}$, \ then
 \[
   \limsup_{x\to\infty} \frac{\PP(Y<-x)}{\PP(|Y|>x)}
   \leq \lim_{x\to\infty} \frac{\PP(Y\leq-x)}{\PP(|Y|>x)} = q ,
 \]
 and
 \begin{align*}
  \liminf_{x\to\infty} \frac{\PP(Y<-x)}{\PP(|Y|>x)}
  &\geq \liminf_{x\to\infty} \frac{\PP(Y\leq-x-1)}{\PP(|Y|>x)} \\
  &= \liminf_{x\to\infty} \frac{\PP(Y\leq-x-1)}{\PP(|Y|>x+1)} \frac{\PP(|Y|>x(1+1/x))}{\PP(|Y|>x)}
   = q ,
 \end{align*}
 since, by the uniform convergence theorem for regularly varying functions (see, e.g., Bingham et al.\ \cite[Theorem 1.5.2]{BinGolTeu}) together with the fact that \ $1 + 1/x \in [1, 2]$ \ for \ $x \in [1, \infty)$, \ we obtain
 \[
   \lim_{x\to\infty} \frac{\PP(|Y|>x(1+1/x))}{\PP(|Y|>x)} = 1 ,
 \]
 and hence, we conclude \eqref{tail-balance+}.

On the other hand, if \ $Y$ \ is a random variable such that \ $\PP(|Y| > x) \in \RR_{++}$ \ for all \ $x \in \RR_{++}$, \ the function \ $\RR_{++} \ni x \mapsto \PP(|Y| > x) \in \RR_{++}$ \ is regularly varying at infinity with index \ $-\alpha$, \ and \eqref{tail-balance+} holds, then the second convergence in the tail-balance condition \eqref{TB} can be derived in a similar way.
\proofend
\end{Rem}

\begin{Lem}\label{hom}
\renewcommand{\labelenumi}{{\rm(\roman{enumi})}}
 \begin{enumerate}
  \item
   A non-negative random variable \ $Y$ \ is regularly varying with index \ $\alpha \in \RR_{++}$ \ if and only if \ $\PP(Y > x) \in \RR_{++}$ \ for all \ $x \in \RR_{++}$, \ and the function \ $\RR_{++} \ni x \mapsto \PP(Y > x) \in \RR_{++}$ \ is regularly varying at infinity with index \ $-\alpha$.
  \item
   If \ $Y$ \ is a regularly varying random variable with index \ $\alpha \in \RR_{++}$, \ then for each \ $\beta \in \RR_{++}$, \ $|Y|^\beta$ \ is regularly varying with index \ $\alpha/\beta$.
 \end{enumerate}
\end{Lem}

\begin{Lem}\label{a_n}
If \ $Y$ \ is a regularly varying random variable with index \ $\alpha \in \RR_{++}$, \ then there exists a sequence \ $(a_n)_{n\in\NN}$ \ in \ $\RR_{++}$ \ such that \ $n \PP(|Y| > a_n) \to 1$ \ as \ $n \to \infty$.
\ If \ $(a_n)_{n\in\NN}$ \ is such a sequence, then \ $a_n \to \infty$ \ as \ $n \to \infty$.
\end{Lem}

\noindent{\bf Proof.}
We are going to show that one can choose \ $a_n := \max\{\ta_n, 1\}$, \ $n \in \NN$, \ where \ $\ta_n$ \ denotes the \ $1 - \frac{1}{n}$ \ lower quantile of \ $|Y|$, \ namely,
 \[
   \ta_n := \inf\biggl\{x \in \RR : 1 - \frac{1}{n} \leq \PP(|Y| \leq x)\biggr\}
   = \inf\biggl\{x \in \RR : \PP(|Y| > x) \leq \frac{1}{n}\biggr\} , \qquad n \in \NN .
 \]
For each \ $n \in \NN$, \ by the definition of the infimum, there exists a sequence \ $(x_m)_{m\in\NN}$ \ in \ $\RR$ \ such that \ $x_m \downarrow \ta_n$ \ as \ $m \to \infty$ \ and \ $\PP(|Y| > x_m) \leq \frac{1}{n}$, \ $m \in \NN$.
\ Letting \ $m \to \infty$, \ using that the distribution function of \ $|Y|$ \ is right-continuous, we obtain \ $\PP(|Y| > \ta_n) \leq \frac{1}{n}$, \ thus \ $n \PP(|Y| > \ta_n) \leq 1$, \ and hence
 \begin{equation}\label{limsup_a_n}
  \limsup_{n\to\infty} n \PP(|Y| > \ta_n) \leq 1 .
 \end{equation}
Moreover, for each \ $n \in \NN$, \ again by the definition of the infimum, we have \ $\frac{1}{n} < \PP(|Y| > \ta_n - 1)$, \ thus \ $n \PP(|Y| > \ta_n - 1) > 1$, \ and hence
 \begin{equation}\label{liminf_a_n}
  \liminf_{n\to\infty} n \PP(|Y| > \ta_n - 1) \geq 1 .
 \end{equation}
We have \ $\ta_n \to \infty$ \ as \ $n \to \infty$, \ since \ $|Y|$ \ is regularly variable with index \ $\alpha \in \RR_+$ \ (see part (ii) of Lemma \ref{hom}), yielding that \ $|Y|$ \ is unbounded.
Thus for each \ $q \in (0, 1)$ \ and for sufficiently large \ $n \in \NN$, \ we have \ $\ta_n \geq \frac{1}{1-q}$, \ and then \ $\ta_n - 1 \geq q \ta_n$, \ and hence \ $\PP(|Y| > \ta_n - 1) \leq \PP(|Y| > q \widetilde a_n)$.
\ Consequently, for each \ $q \in (0, 1)$, \ using \eqref{liminf_a_n} and that \ $|Y|$ \ is regularly varying with index \ $\alpha \in \RR_{++}$, \ we obtain
 \begin{align*}
  1 &\leq \liminf_{n\to\infty} n \PP(|Y| > \ta_n - 1)
  \leq \liminf_{n\to\infty} n \PP(|Y| > q \ta_n) \\
  &= \liminf_{n\to\infty} \frac{\PP(|Y| > q \ta_n)}{\PP(|Y| > \ta_n)} n \PP(|Y| > \ta_n)
   = q^{-\alpha} \liminf_{n\to\infty} n \PP(|Y| > \ta_n) .
 \end{align*}
Hence for each \ $q \in (0, 1)$, \ we have \ $\liminf_{n\to\infty} n \PP(|Y| > \ta_n) \geq q^\alpha$.
\ Letting \ $q \uparrow 1$, \ we get \ $\liminf_{n\to\infty} n \PP(|Y| > \ta_n) \geq 1$, \ and hence by \eqref{limsup_a_n}, we conclude \ $\lim_{n\to\infty} n \PP(|Y| > \ta_n) = 1$.

If \ $(a_n)_{n\in\NN}$ \ is a sequence in \ $\RR_{++}$ \ such that \ $n \PP(|Y| > a_n) \to 1$ \ as \ $n \to \infty$, \ then \ $a_n \to \infty$ \ as \ $n \to \infty$, \ since \ $|Y|$ \ is unbounded.
\proofend

\begin{Lem}[Karamata's theorem for truncated moments]\label{truncated_moments}
Consider a non-negative regularly varying random variable \ $Y$ \ with index \ $\alpha \in \RR_{++}$.
\ Then
 \begin{align*}
  \lim_{x\to\infty}
   \frac{x^\beta\PP(Y>x)}{\EE(Y^\beta\bbone_{\{Y\leq x\}})}
  &= \frac{\beta-\alpha}{\alpha} \qquad \text{for \ $\beta \in [\alpha, \infty)$,} \\
  \lim_{x\to\infty}
   \frac{x^\beta\PP(Y>x)}{\EE(Y^\beta\bbone_{\{Y>x\}})}
  &= \frac{\alpha-\beta}{\alpha} \qquad \text{for \ $\beta \in (-\infty, \alpha)$.}
 \end{align*}
\end{Lem}

For Lemma \ref{truncated_moments}, see, e.g., Bingham et al.\ \cite[pages 26-27]{BinGolTeu} or Buraczewski et al.\ \cite[Appendix B.4]{BurDamMik}.

Next, based on Buraczewski et al.\ \cite[Appendix C]{BurDamMik}, we recall the definition
 and some properties of regularly varying random vectors.

\begin{Def}
A \ $d$-dimensional random vector \ $\bY$ \ and its distribution are called regularly
 varying with index \ $\alpha \in \RR_{++}$
 \ if there exists a probability measure \ $\psi$ \ on \ $\SSS^{d-1}$ \ such that for all \ $c \in \RR_{++}$,
 \[
   \frac{\PP\bigl(\|\bY\| > c x, \, \frac{\bY}{\|\bY\|} \in \cdot\bigr)}
        {\PP(\|\bY\| > x)}
   \distrw c^{-\alpha} \psi(\cdot) \qquad \text{as \ $x \to \infty$,}
 \]
 where \ $\distrw$ \ denotes the weak convergence of finite measures on \ $\SSS^{d-1}$.
\ The probability measure \ $\psi$ \ is called the spectral measure of \ $\bY$.
\end{Def}

The following equivalent characterization of multivariate regular variation can be derived, e.g., from Resnick \cite[page 69]{Res0}.

\begin{Pro}\label{vague}
A $d$-dimensional random vector \ $\bY$ \ is regularly varying with some index \ $\alpha \in \RR_{++}$ \
 if and only if there exists a non-null locally finite measure \ $\mu$ \ on \ $\RR^d_0$ \ satisfying the limit relation
 \begin{equation}\label{vague_Kallenberg}
   \mu_x(\cdot)
   := \frac{\PP(x^{-1} \bY \in \cdot)}{\PP(\|\bY\| > x)}
   \distrv \mu(\cdot) \qquad \text{as \ $x \to \infty$,}
 \end{equation}
 where \ $\distrv$ \ denotes vague convergence of locally finite measures on \ $\RR^d_0$ \ (see Appendix \ref{App_vague}
  for the notion \ $\distrv$).
\ Further, \ $\mu$ \ satisfies the property \ $\mu(cB) = c^{-\alpha}\mu(B)$ \ for any \ $c\in\RR_{++}$ \ and
 \ $B\in\cB(\RR_0^d)$ \ (see, e.g., Theorem 1.14 and 1.15 and Remark 1.16 in Lindskog \cite{Lin}).
\end{Pro}

The measure \ $\mu$ \ in Proposition \ref{vague} is called the limit measure of \ $\bY$.

\noindent{\bf Proof of Proposition \ref{vague}}.
Recall that a \ $d$-dimensional random vector \ $\bY$ \ is regularly varying with some index
 \ $\alpha\in\RR_{++}$ \ if and only if on \ $(\oRR_0^d, \cB(\oRR_0^d))$, \ furnished with an appropriate metric \ $\ovarrho$ \
 (see, e.g., Kallenberg \cite[page 125]{kallenberg:2017}), the vague convergence \ $\mu_x \distrv \omu$ \ as \ $x \to \infty$ \ holds
 with some non-null locally finite measure \ $\omu$ \ with \ $\omu(\oRR_0^d\setminus\RR_0^d) = 0$, \ where \ $\oRR_0^d := \oRR^d \setminus \{\bzero\}$ \
 with \ $\oRR := \RR \cup \{-\infty, \infty\}$, \ see, e.g., Resnick \cite[page 69]{Res0}.
It remains to check that \ $\mu_x \distrv \omu$ \ as \ $x\to\infty$ \ on \ $(\oRR_0^d, \cB(\oRR_0^d))$ \ holds if and only if
 \ $\mu_x \distrv \mu$ \ as \ $x\to\infty$ \ on \ $(\RR_0^d, \cB(\RR_0^d))$ \ with \ $\mu:=\omu\big\vert_{\RR_0^d}$.
\ By Lemma \ref{portmanteau}, \ $\mu_x(B\cap\RR_0^d)=\mu_x(B)\to \omu(B)=\omu(B\cap\RR_0^d)$ \ as \ $x\to\infty$ \ for any bounded \ $\omu$-continuity Borel set \ $B$ \ of \ $\oRR_0^d$.
\ By Kallenberg \cite[page 125]{kallenberg:2017} and Lemma \ref{Lem_S_top}, a subset \ $B$ \ of \ $\oRR_0^d$ \ is bounded with respect to the metric
 \ $\overline{\varrho}$ \ if and only if \ $B\cap \RR_0^d$ \ (as a subset of \ $\RR_0^d$) \ is bounded with respect to the metric \ $\varrho$.
\ Further, for any \ $B\in\cB(\oRR_0^d)$, \ $(\partial_{\oRR_0^d}B) \cap \RR_0^d = \partial_{\RR_0^d}(B\cap\RR_0^d)$, \ where
\ $\partial_{\oRR_0^d}B$ \ and \ $\partial_{\RR_0^d}(B\cap\RR_0^d)$ \ denotes the boundary of \ $B$ \ in \ $\oRR_0^d$ \ and that of
 \ $(B\cap\RR_0^d)$ \ in \ $\RR_0^d$, \ respetively, since a set \ $G \subset \oRR_0^d$ \ is open with respect to \ $\ovarrho$ \ if and only
 if \ $G \cap \RR^d_0$ \ is open with respect to \ $\varrho$.
\ Thus \ $\omu(\partial_{\oRR_0^d}B)=\omu((\partial_{\oRR_0^d}B) \cap \RR_0^d)=0$ \ if and only if \ $\mu(\partial_{\RR_0^d}(B \cap \RR_0^d))=0$.
\ Hence \ $\mu_x(B)\to \omu(B)$ \ as \ $x\to\infty$ \ for any bounded \ $\omu$-continuity set \ $B$ \ of \ $\oRR_0^d$ \ if and only if \ $\mu_x(B) \to \mu(B)$ \ as \ $x\to\infty$ \ for any bounded \ $\mu$-continuity set \ $B$ \ of \ $\RR_0^d$.
\ Consequently, by Lemma \ref{portmanteau}, \ $\mu_x \distrv \omu$ \ as \ $x \to \infty$ \ on \ $\oRR^d_0$ \ if and only if
 \ $\mu_x\distrv \mu$ \ as \ $x\to\infty$ \ on \ $\RR_0^d$.
\proofend

The next statement follows, e.g., from part (i) in Lemma C.3.1 in Buraczewski et al.\ \cite{BurDamMik}.

\begin{Lem}\label{Lem_shift}
If \ $\bY$ \ is a regularly varying \ $d$-dimensional random vector with index \ $\alpha \in \RR_{++}$, \ then for each \ $\bc \in \RR^d$, \ the random vector \ $\bY - \bc$ \ is regularly varying with index \ $\alpha$.
\end{Lem}

Recall that if \ $\bY$ \ is a regularly varying \ $d$-dimensional random vector with index \ $\alpha \in \RR_{++}$ \ and with limit measure \ $\mu$ \ given in \eqref{vague_Kallenberg}, and \ $f: \RR^d \to \RR$ \ is a continuous function with \ $f^{-1}(\{0\}) = \{\bzero\}$ \ and it is positively homogeneous of degree \ $\beta \in \RR_{++}$ \ (i.e., \ $f(c \bv) = c^\beta f(\bv)$ \ for every \ $c \in \RR_{++}$ \ and \ $\bv \in \RR^d)$, \ then \ $f(\bY)$ \ is regularly varying with index \ $\frac{\alpha}{\beta}$ \ and with limit measure \ $\mu(f^{-1}(\cdot))$, \ see, e.g., Buraczewski et al.\ \cite[page 282]{BurDamMik}.
Next we describe the tail behavior of \ $f(\bY)$ \ for appropriate positively homogeneous functions \ $f: \RR^d \to \RR$.

\begin{Pro}\label{Pro_mapping}
Let \ $\bY$ \ be a regularly varying \ $d$-dimensional random vector with index \ $\alpha \in \RR_{++}$ \ and let \ $f: \RR^d \to \RR$ \ be a measurable function which is positively homogeneous of degree \ $\beta \in \RR_{++}$, \ continuous at \ $\bzero$ \ and \ $\mu(D_f) = 0$, \ where \ $\mu$ \ is the limit measure of \ $\bY$ \ given in \eqref{vague_Kallenberg} and \ $D_f$ \ denotes the set of discontinuities of \ $f$.
\ Then \ $\mu(\partial_{\RR_0^d}(f^{-1}((1, \infty)))) = 0$, \ where
\ $\partial_{\RR_0^d}(f^{-1}((1, \infty)))$ \ denotes the boundary of \ $f^{-1}((1, \infty))$ \ in \ $\RR_0^d$.
\ Consequently,
 \[
   \lim_{x\to\infty} \frac{\PP(f(\bY)>x)}{\PP(\|\bY\|^\beta>x)}
   = \mu(f^{-1}((1,\infty))) ,
 \]
 and \ $f(\bY)$ \ is regularly varying with tail index \ $\frac{\alpha}{\beta}$.
\end{Pro}

\noindent{\bf Proof.}
For all \ $x \in \RR_{++}$, \ we have
 \[
   \frac{\PP(f(\bY)>x)}{\PP(\|\bY\|^\beta>x)}
   = \frac{\PP(x^{-1}f(\bY)>1)}{\PP(\|\bY\|>x^{1/\beta})}
   = \frac{\PP(f(x^{-1/\beta}\bY)>1)}{\PP(\|\bY\|>x^{1/\beta})}
   = \frac{\PP(x^{-1/\beta}\bY\in f^{-1}((1,\infty)))}{\PP(\|\bY\|>x^{1/\beta})} .
 \]
Next, we check that \ $f^{-1}((1,\infty))$ \ is a \ $\mu$-continuity set being bounded with respect to the metric \ $\varrho$ \ given in \eqref{S_metric}.
Since \ $f(\bzero) = 0$ \ (following from the positive homogeneity of \ $f$), \ we have \ $f^{-1}((1,\infty)) \in \cB(\RR_0^d)$.
\ The continuity of \ $f$ \ at \ $\bzero$ \ implies the existence of an \ $\vare \in \RR_{++}$ \ such that for all \ $\bx \in \RR^d$ \ with \ $\|\bx\| \leq \vare$ \ we have \ $|f(\bx)| \leq 1$, \ thus \ $\bx \notin f^{-1}((1, \infty))$, \ hence \ $f^{-1}((1, \infty)) \subset \{\bx \in \RR_0^d : \|\bx\| > \vare\}$, \ i.e., \ $f^{-1}((1, \infty))$ \ is separated from \ $\bzero$, \ and hence, by Lemma \ref{Lem_S_top}, \ $f^{-1}((1, \infty))$ \ is bounded in \ $\RR_0^d$ \ with respect to the metric \ $\varrho$.
\ Further, we have
 \[
   \partial_{\RR_0^d}(f^{-1}((1, \infty))) \subset f^{-1}(\partial_{\RR}((1,\infty))) \cup D_f = f^{-1}(\{1\}) \cup D_f ,
 \]
 and hence
 \[
   \mu(\partial_{\RR_0^d}(f^{-1}((1, \infty))))
   \leq \mu(f^{-1}(\{1\}))  +  \mu(D_f) = \mu(f^{-1}(\{1\})) .
 \]
Here \ $\mu(f^{-1}(\{1\})) = 0$, \ since if, on the contrary, we suppose that \ $\mu(f^{-1}(\{1\})) \in (0,
\infty]$, \ then for all \ $u, v \in \RR_{++}$ \ with \ $u < v$, \ we have
 \begin{align*}
  \mu(f^{-1}((u,v))) \geq \mu\left(\bigcup_{q\in\QQ\cap(u,v)} f^{-1}(\{q\})\right)
  &= \sum_{q\in\QQ\cap(u,v)} \mu(f^{-1}(\{q\}))
   = \sum_{q\in\QQ\cap(u,v)} \mu(q^{\frac{1}{\beta}} f^{-1}(\{1\})) \\
  &= \sum_{q\in\QQ\cap(u,v)} q^{-\frac{\alpha}{\beta}} \mu(f^{-1}(\{1\}))
   = \infty ,
 \end{align*}
 where we used that \ $\mu(cB) = c^{-\alpha} \mu(B)$, \ $c \in \RR_{++}$, \ $B \in \cB(\RR_0^d)$ \ (see Proposition \ref{vague}), and that
 \begin{align*}
  f^{-1}(\{q\}) &= \{\bx \in \RR_0^d : f(\bx) = q\}
                 = \{\bx \in \RR_0^d : f(q^{-\frac{1}{\beta}}\bx) = 1\} \\
                &= \{q^{\frac{1}{\beta}} \by \in \RR_0^d : f(\by) = 1\}
                 = q^{\frac{1}{\beta}} f^{-1}(\{1\}) , \qquad q \in \RR_{++} .
 \end{align*}
This leads us to a contradiction, since \ $f^{-1}((u, v))$ \ is separated from \ $\bzero$ \ (can be seen similarly as for \ $f^{-1}((1, \infty))$), so, by Lemma \ref{Lem_S_top}, it is bounded  with respect to the metric \ $\varrho$, \ and hence \ $\mu(f^{-1}((u, v))) < \infty$ \ due to the local finiteness of \ $\mu$.
\ Hence \ $\mu(\partial_{\RR_0^d}(f^{-1}((1, \infty)))) = 0$, \ as desired.

Consequently, by portmanteau theorem for vague convergence (see Lemma \ref{portmanteau}), we have
 \[
   \frac{\PP(x^{-1/\beta}\bY\in f^{-1}((1,\infty)))}{\PP(\|\bY\|>x^{1/\beta})}
   \to \mu(f^{-1}((1, \infty))) \qquad \text{as \ $x \to \infty$,}
 \]
 as desired.
\proofend

\section{Weak convergence of partial sum processes towards L\'evy processes}\label{App_Resnick_gen}

We formulate a slight modification of Theorem 7.1 in Resnick \cite{Res} with a different centering.

\begin{Thm}\label{7.1}
Suppose that for each \ $N \in \NN$, \ $\bX_{N,j}$, \ $j \in \NN$, \ are independent and identically distributed \ $d$-dimensional random vectors such that
 \begin{equation}\label{(7.5)}
  N \PP(\bX_{N,1} \in \cdot) \distrv \nu(\cdot) \qquad \text{on \ $\RR_0^d$ \ as \ $N \to \infty$,}
 \end{equation}
 where \ $\nu$ \ is a L\'evy measure on \ $\RR_0^d$ \ such that \ $\nu(\{\bx \in \RR^d_0 : |\langle\be_\ell, \bx\rangle| = 1\}) = 0$ \ for every \ $\ell \in \{1, \ldots, d\}$, \ and that
 \begin{equation}\label{(7.6)}
  \lim_{\vare\downarrow0} \limsup_{N\to\infty} N \EE\bigl(\langle\be_\ell,\bX_{N,1}\rangle^2 \bbone_{\{|\langle\be_\ell,\bX_{N,1}\rangle|\leq\vare\}}\bigr) = 0 ,
   \qquad \ell \in \{1, \ldots, d\} .
 \end{equation}
Then we have
 \[
   \Biggl(\sum_{j=1}^{\lfloor Nt\rfloor} \biggl(\bX_{N,j} - \sum_{\ell=1}^d \EE\bigl(\langle\be_\ell, \bX_{N,j}\rangle \bbone_{\{|\langle\be_\ell,\bX_{N,j}\rangle|\leq1\}}\bigr) \be_\ell\biggr)\Biggr)_{t\in\RR_+}
   \distr (\bcX_t)_{t\in\RR_+} \qquad
   \text{as \ $N \to \infty$,}
 \]
 where \ $(\bcX_t)_{t\in\RR_+}$ \ is a L\'evy process such that the characteristic function of
 the distribution \ $\mu$ \ of \ $\bcX_1$ \ has the form
 \begin{equation}\label{hmu}
  \hmu(\btheta)
   = \exp\biggl\{\int_{\RR^d_0} \biggl(\ee^{\ii\langle\btheta,\bx\rangle} - 1 - \ii \sum_{\ell=1}^d \langle\be_\ell, \btheta\rangle \langle\be_\ell,
   \bx\rangle \bbone_{(0,1]}(|\langle\be_\ell,\bx\rangle|)\biggr) \nu(\dd\bx)\biggr\} , \quad \btheta \in \RR^d .
 \end{equation}
\end{Thm}

\noindent{\bf Proof.}
There exists \ $r \in \RR_{++}$ \ such that \ $\nu(\{\bx \in \RR^d_0 : \|\bx\| = r\}) = 0$, \ since the function \ $\RR_{++} \ni t \mapsto \nu(\{\bx \in \RR^d_0 : \|\bx\| > t\})$ \ is decreasing.
By an appropriate modification of Theorem 7.1 in Resnick \cite{Res}, we obtain
 \[
   \Biggl(\sum_{j=1}^{\lfloor Nt\rfloor} \bigl(\bX_{N,j}
       - \EE\bigl(\bX_{N,j} \bbone_{\{\Vert \bX_{N,j} \Vert\leq r\}}\bigr)\bigr)\Biggr)_{t\in\RR_+}
   \distr (\tbcX_t)_{t\in\RR_+} \qquad
   \text{as \ $N \to \infty$,}
 \]
 where \ $(\tbcX_t)_{t\in\RR_+}$ \ is a L\'evy process such that the characteristic function of \ $\tbcX_1$ \ has the form
 \begin{equation*}
  \EE(\ee^{\ii \langle \btheta, \tbcX_1 \rangle})
     = \exp\biggl\{\int_{\RR^d_0} \biggl(\ee^{\ii\langle\btheta,\bx\rangle} - 1
           - \ii \langle \btheta,\bx \rangle \bbone_{(0,r]}(\Vert \bx \Vert)\biggr) \nu(\dd\bx)\biggr\} , \qquad \btheta \in \RR^d .
 \end{equation*}
Let us consider the decomposition
 \begin{align*}
   &\sum_{j=1}^{\lfloor Nt\rfloor} \biggl(\bX_{N,j} - \sum_{\ell=1}^d \EE\bigl(\langle\be_\ell, \bX_{N,j}\rangle \bbone_{\{|\langle\be_\ell,\bX_{N,j}\rangle|\leq1\}}\bigr) \be_\ell\biggr) \\
   &= \sum_{j=1}^{\lfloor Nt\rfloor} \bigl(\bX_{N,j}
         - \EE\bigl(\bX_{N,j} \bbone_{\{\Vert \bX_{N,j} \Vert\leq r\}}\bigr)\bigr)
         + \sum_{\ell=1}^d \sum_{j=1}^{\lfloor Nt\rfloor} \EE\bigl(\langle\be_\ell, \bX_{N,j}\rangle  ( \bbone_{\{ \Vert \bX_{N,j} \Vert\leq r\}}
         - \bbone_{\{|\langle\be_\ell,\bX_{N,j}\rangle|\leq1\}}   ) \bigr) \be_\ell
 \end{align*}
 for each \ $t\in\RR_{++}$.
Here for each \ $\ell\in\{1,\ldots,d \}$, \ we have
 \begin{align*}
  &\sum_{j=1}^{\lfloor Nt\rfloor} \EE\bigl(\langle\be_\ell, \bX_{N,j}\rangle  ( \bbone_{\{ \Vert \bX_{N,j} \Vert\leq r\}}
         - \bbone_{\{|\langle\be_\ell,\bX_{N,j}\rangle|\leq1\}}   ) \bigr) \\
  &= \lfloor Nt\rfloor \EE\bigl(\langle\be_\ell, \bX_{N,1}\rangle  ( \bbone_{\{ \Vert \bX_{N,1} \Vert\leq r\}}
         - \bbone_{\{|\langle\be_\ell,\bX_{N,1}\rangle|\leq1\}}   ) \bigr)
          = \frac{\lfloor Nt\rfloor }{N} N \EE(g_\ell(\bX_{N,1})),
 \end{align*}
 where \ $g_\ell:\RR^d\to\RR$, \ $g_\ell(\bx):=x_\ell( \bbone_{\{ \Vert \bx \Vert\leq r\}}
                               - \bbone_{\{|x_\ell|\leq1\}} )$, \ $\bx=(x_1,\ldots,x_d)^\top\in\RR^d$.
\ For each \ $\ell\in\{1,\ldots,d\}$, \ the positive and negative parts \ $g_\ell^+$ \ and \ $g_\ell^{-}$ \ of the function
 \ $g_\ell$ \ are bounded, measurable with a bounded support (following from Lemma \ref{Lem_S_top}), and, due to
 \ $\nu(\{\bx \in \RR^d_0 : |\langle\be_\ell, \bx\rangle| = 1\}) = 0$, \ $\ell \in \{1, \ldots, d\}$, \
 and \ $\nu(\{\bx \in \RR^d_0 : \Vert \bx \Vert = r\})=0$, \ the sets of discontinuity points \ $D_{g^+_\ell}$ \ and
 \ $D_{g^-_\ell}$ \ have \ $\nu$-measure \ $0$, \ i.e., \ $\nu(D_{g^+_\ell})=\nu(D_{g^-_\ell})=0$.
\ Consequently, by \eqref{(7.5)} and Lemma \ref{cmt}, we have
 \[
   N \EE(g_\ell(\bX_{N,1})) = N \EE(g^+_\ell(\bX_{N,1})) - N \EE(g^-_\ell(\bX_{N,1}))
     \to  \nu(g^+_\ell) - \nu(g^-_\ell) = \nu(g_\ell) \in \RR
 \]
 as \ $N\to\infty$, \ since \ $\nu(g^+_\ell), \nu(g^-_\ell) \in \RR_+$ \ due to the fact that \ $\nu$ \ is a L\'evy measure.
Next, we may apply Lemma \ref{Conv2Funct} with
 \begin{gather*}
  \bcU_t^{(N)} :=
     \sum_{j=1}^{\lfloor Nt\rfloor}
            \bigl(\bX_{N,j} - \EE\bigl(\bX_{N,j} \bbone_{\{ \Vert \bX_{N,j}\Vert\leq r\}}\bigr) \bigr) ,
      \qquad N \in \NN , \\
  \Phi_N(f)(t) := f(t) + \lfloor Nt\rfloor \sum_{\ell=1}^d \EE(g_\ell(\bX_{N,1}))\be_\ell,
                  \qquad N \in \NN , \\
  \bcU_t := \tbcX_t , \qquad
  \Phi(f)(t) := f(t) +t \sum_{\ell=1}^d \nu(g_\ell)\be_\ell
 \end{gather*}
 for \ $t \in \RR_+$ \ and \ $f \in \DD(\RR_+,\RR^d)$.
\ Indeed, in order to show \ $\Phi_N(f_N) \to \Phi(f)$ \ in \ $\DD(\RR_+, \RR^{k+1})$ \ as \ $N \to \infty$ \
 whenever \ $f_N \to f$ \ in \ $\DD(\RR_+, \RR^{k+1})$ \ as \ $N \to \infty$ \ with \ $f, f_N \in \DD(\RR_+, \RR^{k+1})$, \ $N \in \NN$,
 \ by Propositions VI.1.17 and VI.1.23 in Jacod and Shiryaev \cite{JacShi}, it is enough to check that
 for each \ $T \in \RR_{++}$, \ we have
 \begin{align*}
  \sup_{t\in[0,T]} \sum_{\ell=1}^d \biggl| \lfloor Nt\rfloor \EE(g_\ell(\bX_{N,1})) - t\nu(g_\ell) \biggr| \to 0
   \qquad \text{as \ $N\to\infty$.}
 \end{align*}
This follows, since for each \ $\ell\in\{1,\ldots,d\}$,
 \begin{align*}
  &\sup_{t\in[0,T]} \biggl|   \lfloor Nt\rfloor \EE(g_\ell(\bX_{N,1})) - t\nu(g_\ell) \biggr| \\
  &\leq \sup_{t\in[0,T]} \biggl|\frac{\lfloor Nt\rfloor}{N} \bigl( N\EE(g_\ell(\bX_{N,1})) - \nu(g_\ell) \bigr)\biggr|
        + \sup_{t\in[0,T]} \Biggl| \nu(g_\ell)\biggl( \frac{\lfloor Nt\rfloor}{N} -t  \biggr)  \Biggr| \\
  &\leq T \vert N\EE(g_\ell(\bX_{N,1})) - \nu(g_\ell) \vert + \frac{ \nu(g_\ell)}{N}
   \to 0 \qquad \text{as \ $N \to \infty$.}
 \end{align*}
Applying Lemma \ref{Conv2Funct}, we obtain
 \[
   \biggl(\sum_{j=1}^{\lfloor Nt\rfloor}
         \biggl(\bX_{N,j} - \sum_{\ell=1}^d \EE\bigl(\langle\be_\ell, \bX_{N,j}\rangle \bbone_{\{|\langle\be_\ell,\bX_{N,j}\rangle|\leq1\}}\bigr) \be_\ell\biggr)
   \biggr)_{t\in\RR_+}
   = \Phi_N(\bcU^{(N)})
   \distr \Phi(\bcU) \qquad \text{as \ $N \to \infty$,}
 \]
 where \ $\Phi(\bcU)_t = \tbcX_t + t \sum_{\ell=1}^d \nu(g_\ell)\be_\ell = \bcX_t$,
 \ $t \in \RR_+$, \ is a $d$-dimensional L\'evy process, since
 \begin{align*}
   \EE\big( \ee^{\ii \langle \btheta, \tbcX_1 + \sum_{\ell=1}^d \nu(g_\ell)\be_\ell \rangle}\big)
    & = \exp\biggl\{\int_{\RR^d_0} \biggl(\ee^{\ii\langle\btheta,\bx\rangle} - 1
            - \ii \langle \btheta,\bx \rangle {\bbone_{(0,r]}}(\Vert \bx \Vert)\biggr) \nu(\dd\bx)
            + \ii \sum_{\ell=1}^d \langle \btheta,\be_\ell \rangle \nu(g_\ell)
           \biggr\} \\
    & = \exp\biggl\{\int_{\RR^d_0} \biggl(\ee^{\ii\langle\btheta,\bx\rangle} - 1
            - \ii \langle \btheta,\bx \rangle {\bbone_{(0,r]}}(\Vert \bx \Vert)\biggr) \nu(\dd\bx) \\
    &\phantom{= \exp\biggl\{}
            + \ii \sum_{\ell=1}^d \langle \btheta,\be_\ell \rangle
              \int_{\RR_0^d} \langle \be_\ell,\bx \rangle (  \bbone_{\{\Vert \bx \Vert\leq r \}}
                               - \bbone_{\{ \vert \langle \be_\ell,\bx \rangle \vert \leq1\}})\,\nu(\dd \bx)
           \biggr\},
 \end{align*}
 yielding \eqref{hmu}.
\proofend

\section{Tail behavior of \ $(X_k)_{k\in\ZZ_+}$}\label{App_tail}

Due to Basrak et al.\ \cite[Theorem 2.1.1]{BasKulPal}, we have the following tail
 behavior.

\begin{Thm}\label{Xtail}
We have
 \[
   \lim_{x\to\infty} \frac{\pi((x, \infty))}{\PP(\vare > x)}
   = \sum_{i=0}^\infty m_\xi^{i\alpha}
   = \frac{1}{1-m_\xi^\alpha} ,
 \]
 where \ $\pi$ \ denotes the unique stationary distribution of the Markov chain \ $(X_k)_{k\in\ZZ_+}$, \ and consequently, \ $\pi$ \ is also regularly varying with index \ $\alpha$.
\end{Thm}

Note that in case of \ $\alpha = 1$ \ and \ $m_\vare = \infty$ \ Basrak et al.\
 \cite[Theorem 2.1.1]{BasKulPal} assume additionally that \ $\vare$ \ is consistently varying (or
 in other words intermediate varying), but, eventually, it follows from the fact that \ $\vare$ \ is regularly
 varying.

Let \ $(X_k)_{k\in\ZZ}$ \ be a strongly stationary extension of \ $(X_k)_{k\in\ZZ_+}$.
\ Basrak et al. \cite[Lemma 3.1]{BasKulPal} described the so-called forward tail process of the
 strongly stationary process \ $(X_k)_{k\in\ZZ}$, \ and hence, due to Basrak and
 Segers \cite[Theorem 2.1]{BasSeg}, the strongly stationary process
 \ $(X_k)_{k\in\ZZ}$ \ is jointly regularly varying.

\begin{Thm}\label{Xtailprocess}
The finite dimensional conditional distributions of \ $(x^{-1} X_k)_{k\in\ZZ_+}$ \ with respect
 to the condition \ $X_0 > x$ \ converge weakly to the corresponding finite dimensional distributions of
 \ $(m_\xi^k Y)_{k\in\ZZ_+}$ \ as \ $x \to \infty$, \ where \ $Y$ \ is a random
 variable with Pareto distribution
 \ $\PP(Y \leq y) = (1 - y^{-\alpha}) \bbone_{[1,\infty)}(y)$, \ $y \in \RR$.
\ Consequently, the strongly stationary process \ $(X_k)_{k\in\ZZ}$ \ is jointly
 regularly varying with index \ $\alpha$, \ i.e., all its finite dimensional
 distributions are regularly varying with index \ $\alpha$.
\ The process \ $(m_\xi^k Y)_{k\in\ZZ_+}$ \ is the so called forward tail process of $(X_k)_{k\in\ZZ}$.
\ Moreover, there exists a (whole) tail process of $(X_k)_{k\in\ZZ}$ as well.
\end{Thm}

By the proof of Theorem \ref{simple_aggregation1_stable_fd} and Proposition \ref{Pro_mapping}, we obtain the following results.

\begin{Pro}\label{Pro_limit_meaure}
For each \ $k \in \ZZ_+$,
\renewcommand{\labelenumi}{{\rm(\roman{enumi})}}
 \begin{enumerate}
  \item
   the limit measure \ $\tnu_{k,\alpha}$ \ of \ $(X_0, \ldots, X_k)^\top$ \ given in \eqref{tnu_alpha^k_old} takes the form
    \[
      \tnu_{k,\alpha} = \frac{\nu_{k,\alpha}}{\nu_{k,\alpha}(\{\bx\in\RR_0^{k+1}:\|\bx\|>1\})} ,
    \]
    where \ $\nu_{k,\alpha}$ \ is given by \eqref{fint} and
    \[
      \nu_{k,\alpha}(\{\bx\in\RR_0^{k+1}:\|\bx\|>1\})
      = \frac{1-m_\xi^\alpha}{(1-m_\xi^2)^{\alpha/2}}
        \biggl(\frac{(1-m_\xi^{2(k+1)})^{\alpha/2}}{1-m_\xi^\alpha}
               + \sum_{j=1}^k (1-m_\xi^{2(k-j+1)})^{\alpha/2}\biggr) ;
    \]
  \item
   the tail behavior of \ $X_0 + \cdots + X_k$ \ is given by
    \[
      \lim_{x\to\infty} \frac{\PP(X_0+\cdots+X_k>x)}{\PP(X_0>x)}
      = \frac{1-m_\xi^\alpha}{(1-m_\xi)^\alpha}
        \biggl(\frac{ (1-m_\xi^{k+1})^\alpha }{1-m_\xi^\alpha}
               + \sum_{j=1}^k (1-m_\xi^{k-j+1})^\alpha\biggr) .
    \]
 \end{enumerate}
\end{Pro}

\noindent{\bf Proof.}
(i). \ In the proof of Theorem \ref{simple_aggregation1_stable_fd}, we derived \ $\nu_{k,\alpha} = \tnu_{k,\alpha}/\tnu_{k,\alpha}(\{\bx\in\RR_0^{k+1}:x_0>1\})$.
\ Consequently,
 \[
   \tnu_{k,\alpha}(\{\bx\in\RR_0^{k+1}:x_0>1\})
   = \frac{\tnu_{k,\alpha}(\{\bx\in\RR_0^{k+1}:\|\bx\|>1\})}{\nu_{k,\alpha}(\{\bx\in\RR_0^{k+1}:\|\bx\|>1\})} ,
 \]
 where, using Proposition \ref{Pro_mapping} with the 1-homogeneous function \ $\RR^{k+1} \ni \bx \mapsto \|\bx\|$, \ we have
 \[
   \tnu_{k,\alpha}(\{\bx\in\RR_0^{k+1}:\|\bx\|>1\})
   = \lim_{x\to\infty} \frac{\PP(\|(X_0,\ldots,X_k)^\top\|>x)}{\PP(\|(X_0,\ldots,X_k)^\top\|>x)} = 1 ,
 \]
 and, by \eqref{fint},
 \begin{align*}
  &\nu_{k,\alpha}(\{\bx \in \RR_0^{k+1} : \|\bx\| > 1\})
   = \int_{\RR_0^{k+1}} \bbone_{\{\|\bx\|>1\}} \, \nu_{k,\alpha}(\dd\bx) \\
  &= (1 - m_\xi^\alpha)
     \sum_{j=0}^k
      \int_0^\infty \bbone_{\{\|u\bv_{j}^{(k)}\|>1\}} \alpha u^{-\alpha-1} \, \dd u \\
  &= (1 - m_\xi^\alpha)
     \sum_{j=0}^k
       \int_{\|\bv_{j}^{(k)}\|^{-1}}^\infty \alpha u^{-\alpha-1} \, \dd u
   = (1 - m_\xi^\alpha) \sum_{j=0}^k \|\bv_j^{(k)}\|^\alpha \\
  &= (1 - m_\xi^\alpha) \biggl((1-m_\xi^\alpha)^{-1}(1+m_\xi^2 + \cdots + m_\xi^{2k})^{\alpha/2}
                                + \sum_{j=1}^k (1+m_\xi^2 + \cdots + m_\xi^{2(k-j)})^{\alpha/2}\biggr) \\
  &= \frac{1-m_\xi^\alpha}{(1-m_\xi^2)^{\alpha/2}}
     \biggl(\frac{(1-m_\xi^{2(k+1)})^{\alpha/2}}{1-m_\xi^\alpha}
            + \sum_{j=1}^k (1-m_\xi^{2(k-j+1)})^{\alpha/2}\biggr) .
 \end{align*}
(ii). \ Applying Proposition \ref{Pro_mapping} for the 1-homogeneous functions \ $\RR^{k+1} \ni \bx \mapsto x_0$ \ and \ $\RR^{k+1} \ni \bx \mapsto x_0 + \cdots + x_k$ \ and formula \eqref{fint}, we obtain
 \begin{align*}
  &\lim_{x\to\infty} \frac{\PP(X_0+\cdots+X_k>x)}{\PP(X_0>x)}
   = \lim_{x\to\infty}
      \frac{\PP(\|(X_0,\ldots,X_k)^\top\|>x)}{\PP(X_0>x)}
      \frac{\PP(X_0+\cdots+X_k>x)}{\PP(\|(X_0,\ldots,X_k)^\top\|>x)} \\
  &= \frac{\tnu_{k,\alpha}(\{\bx\in\RR_0^{k+1}:x_0+\cdots+x_k>1\})}{\tnu_{k,\alpha}(\{\bx\in\RR_0^{k+1}:x_0>1\})}
   = \nu_{k,\alpha}(\{\bx \in \RR_0^{k+1} : x_0 + \cdots + x_k > 1\}) \\
  &= \nu_{k,\alpha}(\{\bx \in \RR_0^{k+1} : \langle\bone_{k+1}, \bx\rangle > 1\})
   = \int_{\RR_0^{k+1}} \bbone_{\{\langle\bone_{k+1}, \bx\rangle > 1\}} \, \nu_{k,\alpha}(\dd\bx) \\
  &= (1 - m_\xi^\alpha)
     \sum_{j=0}^k
      \int_0^\infty \bbone_{\{\langle\bone_{k+1}, u\bv_{j}^{(k)}\rangle > 1\}} \alpha u^{-\alpha-1} \, \dd u \\
  &= (1 - m_\xi^\alpha)
     \sum_{j=0}^k
       \int_{\langle\bone_{k+1}, \bv_{j}^{(k)}\rangle^{-1}}^\infty \alpha u^{-\alpha-1} \, \dd u
   = (1 - m_\xi^\alpha) \sum_{j=0}^k \langle \bone_{k+1},  \bv_j^{(k)} \rangle^\alpha \\
  &= (1 - m_\xi^\alpha) \biggl((1-m_\xi^\alpha)^{-1}(1+m_\xi + \cdots + m_\xi^k)^\alpha
                                + \sum_{j=1}^k (1+m_\xi + \cdots + m_\xi^{k-j})^\alpha\biggr) \\
  & = \frac{1-m_\xi^\alpha}{(1-m_\xi)^\alpha} \biggl(\frac{ (1-m_\xi^{k+1})^\alpha }{1-m_\xi^\alpha}  + \sum_{j=1}^k (1-m_\xi^{k-j+1})^\alpha\biggr) ,
 \end{align*}
 as desired.
\proofend

\bibliographystyle{plain}
\bibliography{reg_var_aggr_3}

\end{document}